\documentclass{eptcs}

\newif\ifignore 
\ignorefalse
\newcommand{\auxproof}[1]{
\ifignore\mbox{}\newline
\textbf{PROOF:} \dotfill\newline
{\it #1}\mbox{}\newline
\textbf{ENDPROOF}\dotfill
\fi}

\usepackage{amsmath}
\usepackage{amssymb}
\usepackage{amsthm}
\usepackage{amsfonts}
\usepackage{xspace}
\usepackage{xy}
\xyoption{all}
\usepackage{hyperref}
\usepackage{stmaryrd}

\newdimen\proofrulebreadth \proofrulebreadth=.05em
\newdimen\proofdotseparation \proofdotseparation=1.25ex
\newdimen\proofrulebaseline \proofrulebaseline=2ex
\newcount\proofdotnumber \proofdotnumber=3
\let\then\relax
\def\hfi{\hskip0pt plus.0001fil}
\mathchardef\squigto="3A3B
\newif\ifinsideprooftree\insideprooftreefalse
\newif\ifonleftofproofrule\onleftofproofrulefalse
\newif\ifproofdots\proofdotsfalse
\newif\ifdoubleproof\doubleprooffalse
\let\wereinproofbit\relax
\newdimen\shortenproofleft
\newdimen\shortenproofright
\newdimen\proofbelowshift
\newbox\proofabove
\newbox\proofbelow
\newbox\proofrulename
\def\shiftproofbelow{\let\next\relax\afterassignment\setshiftproofbelow\dimen0 }
\def\shiftproofbelowneg{\def\next{\multiply\dimen0 by-1 }%
\afterassignment\setshiftproofbelow\dimen0 }
\def\setshiftproofbelow{\next\proofbelowshift=\dimen0 }
\def\setproofrulebreadth{\proofrulebreadth}

\def\prooftree{
\ifnum  \lastpenalty=1
\then   \unpenalty
\else   \onleftofproofrulefalse
\fi
\ifonleftofproofrule
\else   \ifinsideprooftree
        \then   \hskip.5em plus1fil
        \fi
\fi
\bgroup
\setbox\proofbelow=\hbox{}\setbox\proofrulename=\hbox{}%
\let\justifies\proofover\let\leadsto\proofoverdots\let\Justifies\proofoverdbl
\let\using\proofusing\let\[\prooftree
\ifinsideprooftree\let\]\endprooftree\fi
\proofdotsfalse\doubleprooffalse
\let\thickness\setproofrulebreadth
\let\shiftright\shiftproofbelow \let\shift\shiftproofbelow
\let\shiftleft\shiftproofbelowneg
\let\ifwasinsideprooftree\ifinsideprooftree
\insideprooftreetrue
\setbox\proofabove=\hbox\bgroup$\displaystyle 
\let\wereinproofbit\prooftree
\shortenproofleft=0pt \shortenproofright=0pt \proofbelowshift=0pt
\onleftofproofruletrue\penalty1
}

\def\eproofbit{
\ifx    \wereinproofbit\prooftree
\then   \ifcase \lastpenalty
        \then   \shortenproofright=0pt  
        \or     \unpenalty\hfil         
        \or     \unpenalty\unskip       
        \else   \shortenproofright=0pt  
        \fi
\fi
\global\dimen0=\shortenproofleft
\global\dimen1=\shortenproofright
\global\dimen2=\proofrulebreadth
\global\dimen3=\proofbelowshift
\global\dimen4=\proofdotseparation
\global\count255=\proofdotnumber
$\egroup  
\shortenproofleft=\dimen0
\shortenproofright=\dimen1
\proofrulebreadth=\dimen2
\proofbelowshift=\dimen3
\proofdotseparation=\dimen4
\proofdotnumber=\count255
}

\def\proofover{
\eproofbit 
\setbox\proofbelow=\hbox\bgroup 
\let\wereinproofbit\proofover
$\displaystyle
}%
\def\proofoverdbl{
\eproofbit 
\doubleprooftrue
\setbox\proofbelow=\hbox\bgroup 
\let\wereinproofbit\proofoverdbl
$\displaystyle
}%
\def\proofoverdots{
\eproofbit 
\proofdotstrue
\setbox\proofbelow=\hbox\bgroup 
\let\wereinproofbit\proofoverdots
$\displaystyle
}%
\def\proofusing{
\eproofbit 
\setbox\proofrulename=\hbox\bgroup 
\let\wereinproofbit\proofusing
\kern0.3em$
}

\def\endprooftree{
\eproofbit 
  \dimen5 =0pt
\dimen0=\wd\proofabove \advance\dimen0-\shortenproofleft
\advance\dimen0-\shortenproofright
\dimen1=.5\dimen0 \advance\dimen1-.5\wd\proofbelow
\dimen4=\dimen1
\advance\dimen1\proofbelowshift \advance\dimen4-\proofbelowshift
\ifdim  \dimen1<0pt
\then   \advance\shortenproofleft\dimen1
        \advance\dimen0-\dimen1
        \dimen1=0pt
        \ifdim  \shortenproofleft<0pt
        \then   \setbox\proofabove=\hbox{%
                        \kern-\shortenproofleft\unhbox\proofabove}%
                \shortenproofleft=0pt
        \fi
\fi
\ifdim  \dimen4<0pt
\then   \advance\shortenproofright\dimen4
        \advance\dimen0-\dimen4
        \dimen4=0pt
\fi
\ifdim  \shortenproofright<\wd\proofrulename
\then   \shortenproofright=\wd\proofrulename
\fi
\dimen2=\shortenproofleft \advance\dimen2 by\dimen1
\dimen3=\shortenproofright\advance\dimen3 by\dimen4
\ifproofdots
\then
        \dimen6=\shortenproofleft \advance\dimen6 .5\dimen0
        \setbox1=\vbox to\proofdotseparation{\vss\hbox{$\cdot$}\vss}%
        \setbox0=\hbox{%
                \advance\dimen6-.5\wd1
                \kern\dimen6
                $\vcenter to\proofdotnumber\proofdotseparation
                        {\leaders\box1\vfill}$%
                \unhbox\proofrulename}%
\else   \dimen6=\fontdimen22\the\textfont2 
        \dimen7=\dimen6
        \advance\dimen6by.5\proofrulebreadth
        \advance\dimen7by-.5\proofrulebreadth
        \setbox0=\hbox{%
                \kern\shortenproofleft
                \ifdoubleproof
                \then   \hbox to\dimen0{%
                        $\mathsurround0pt\mathord=\mkern-6mu%
                        \cleaders\hbox{$\mkern-2mu=\mkern-2mu$}\hfill
                        \mkern-6mu\mathord=$}%
                \else   \vrule height\dimen6 depth-\dimen7 width\dimen0
                \fi
                \unhbox\proofrulename}%
        \ht0=\dimen6 \dp0=-\dimen7
\fi
\let\doll\relax
\ifwasinsideprooftree
\then   \let\VBOX\vbox
\else   \ifmmode\else$\let\doll=$\fi
        \let\VBOX\vcenter
\fi
\VBOX   {\baselineskip\proofrulebaseline \lineskip.2ex
        \expandafter\lineskiplimit\ifproofdots0ex\else-0.6ex\fi
        \hbox   spread\dimen5   {\hfi\unhbox\proofabove\hfi}%
        \hbox{\box0}%
        \hbox   {\kern\dimen2 \box\proofbelow}}\doll%
\global\dimen2=\dimen2
\global\dimen3=\dimen3
\egroup 
\ifonleftofproofrule
\then   \shortenproofleft=\dimen2
\fi
\shortenproofright=\dimen3
\onleftofproofrulefalse
\ifinsideprooftree
\then   \hskip.5em plus 1fil \penalty2
\fi
}

\newenvironment{myproof}[1][Proof]%
   { \begin{trivlist}%
     \item[\hskip \labelsep {\bfseries #1}]%
   }%
   { \end{trivlist}%
   }

\makeatother
 \newdir{ >}{{}*!/-7.5pt/@{>}}
 \newdir{|>}{!/4.5pt/@{|}*:(1,-.2)@^{>}*:(1,+.2)@_{>}}
 \newdir{ |>}{{}*!/-3pt/@{|}*!/-7.5pt/:(1,-.2)@^{>}*!/-7.5pt/:(1,+.2)@_{>}}
\newcommand{\xyline}[2][]{\ensuremath{\smash{\xymatrix@1#1{#2}}}}
\newcommand{\xyinline}[2][]{\ensuremath{\smash{\xymatrix@1#1{#2}}}^{\rule[8.5pt]{0pt}{0pt}}}
\newcommand{\filter}{\raisebox{5.5pt}{$\xymatrix@=6pt@H=0pt@M=0pt@W=4pt{\\ \ar@{>->}[u]}$}}
\newcommand{\ideal}{\raisebox{1pt}{$\xymatrix@=5pt@H=0pt@M=0pt@W=4pt{\ar@{>->}[d] \\ \mbox{}}$}}
\makeatletter

\newcommand{\QEDbox}{\square}
\newcommand{\QED}{\hspace*{\fill}$\QEDbox$}

\newcommand{\after}{\mathrel{\circ}}
\newcommand{\scalar}{\mathrel{\bullet}}
\newcommand{\cat}[1]{\ensuremath{\mathbf{#1}}}
\newcommand{\Cat}[1]{\ensuremath{\mathbf{#1}}}
\newcommand{\op}{\ensuremath{^{\mathrm{op}}}}
\newcommand{\idmap}[1][]{\ensuremath{\mathrm{id}_{#1}}}
\newcommand{\orthogonal}{\mathrel{\bot}}

\newcommand{\Sets}{\Cat{Sets}\xspace}

\newcommand{\Conv}{\Cat{Conv}\xspace}
\newcommand{\EA}{\Cat{EA}\xspace}
\newcommand{\EMod}{\Cat{EMod}\xspace}
\newcommand{\AEMod}{\Cat{AEMod}\xspace}
\newcommand{\BEMod}{\Cat{BEMod}\xspace}
\newcommand{\BCM}{\Cat{BCM}\xspace}
\newcommand{\BA}{\Cat{BA}\xspace}
\newcommand{\CH}{\Cat{CH}\xspace}
\newcommand{\CCH}{\Cat{CCH}\xspace}
\newcommand{\CCHobs}{\CCH_{\mathrm{obs}}}
\newcommand{\Algobs}{\Alg_{\mathrm{obs}}}
\newcommand{\poVectu}{\Cat{poVectu}\xspace}
\newcommand{\OUS}{\Cat{OUS}\xspace}
\newcommand{\BOUS}{\Cat{BOUS}\xspace}
\newcommand{\reals}{\ensuremath{\mathbb{R}}}
\newcommand{\Pow}{\mathcal{P}}
\newcommand{\Dst}{\mathcal{D}}
\newcommand{\Exp}{\mathcal{E}}
\newcommand{\UF}{\ensuremath{\mathcal{U}{\kern-.75ex}\mathcal{F}}}
\newcommand{\Mlt}{\mathcal{M}}

\newcommand{\Ef}{\ensuremath{\mathrm{Ef}}}
\newcommand{\pr}{\ensuremath{\mathrm{Pr}}}
\newcommand{\DM}{\ensuremath{\mathrm{DM}}}
\newcommand{\tr}{\ensuremath{\mathrm{tr}}\xspace}
\newcommand{\ev}{\ensuremath{\mathrm{ev}}\xspace}
\newcommand{\cv}{\ensuremath{\mathrm{cv}}\xspace} 
\newcommand{\ch}{\ensuremath{\mathrm{ch}}\xspace} 
\newcommand{\cc}{\ensuremath{\mathrm{cc}}\xspace} 
\newcommand{\hs}[1][]{{\frak{hs}_{#1}}}
\newcommand{\Alg}{\textsl{Alg}\xspace}
\newcommand{\NNO}{\mathbb{N}}

\renewcommand{\H}{\ensuremath{\mathcal{H}}}
\newcommand{\calF}{\mathcal{F}}
\newcommand{\supp}{\textsl{supp}}
\newcommand{\open}{\ensuremath{\mathcal{O}}}
\newcommand{\closed}{\ensuremath{\mathcal{C}{\kern-.45ex}\ell}}
\newcommand{\ket}[1]{\ensuremath{|{\kern.1em}#1{\kern.1em}\rangle}}
\newcommand{\sotimes}{\mathrel{\raisebox{.05pc}{$\scriptstyle \otimes$}}}
\newcommand{\Kl}{\mathcal{K}{\kern-.5ex}\ell}
\newcommand{\set}[2]{\{#1\;|\;#2\}}
\newcommand{\setin}[3]{\{#1\in#2\;|\;#3\}}
\newcommand{\conjun}{\mathrel{\wedge}}
\newcommand{\disjun}{\mathrel{\vee}}
\newcommand{\all}[2]{\forall{#1}.\,#2}
\newcommand{\allin}[3]{\forall{#1\in#2}.\,#3}
\newcommand{\ex}[2]{\exists{#1}.\,#2}
\newcommand{\exin}[3]{\exists{#1\in#2}.\,#3}
\newcommand{\lam}[2]{\lambda#1.\,#2}
\newcommand{\lamin}[3]{\lambda#1\in#2.\,#3}
\newcommand{\tuple}[1]{\langle#1\rangle}
\newcommand{\charac}[1]{\mathbf{1}_{#1}}
\newcommand{\downset}{\mathop{\downarrow}}
\newcommand{\upset}{\mathop{\uparrow}}

\newcommand{\st}{\ensuremath{\mathsf{st}}}
\newcommand{\leftScottint}{[{\kern-.3ex}[}
\newcommand{\rightScottint}{]{\kern-.3ex}]}
\newcommand{\Scottint}[1]{\leftScottint\,#1\,\rightScottint}

\newcommand{\conglongrightarrow}{\mathrel{\stackrel{
           \raisebox{.5ex}{$\scriptstyle\cong\,$}}{
           \raisebox{0ex}[0ex][0ex]{$\longrightarrow$}}}}
\newcommand{\partot}{\ensuremath{\mathcal{T}\hspace{-3pt}{\scriptstyle{o}}}}
\newcommand{\totpar}{\ensuremath{\mathcal{P}\hspace{-3pt}{\scriptstyle{a}}}}
\newcommand{\partotscript}{\ensuremath{\mathcal{T}\hspace{-3pt}{\scriptscriptstyle{o}}}}
\newcommand{\totparscript}{\ensuremath{\mathcal{P}\hspace{-2pt}{\scriptscriptstyle{a}}}}
\newcommand{\pth}{\ensuremath{\hat{\mathcal{T}}\hspace{-3pt}{\scriptstyle{o}}}}
\newcommand{\tph}{\ensuremath{\hat{\mathcal{P}}\hspace{-3pt}{\scriptstyle{a}}}}
\newcommand{\ojoin}{\ovee}
\newcommand{\eps}{\varepsilon}

\theoremstyle{plain}
\newtheorem{theorem}{Theorem}
\newtheorem{lemma}{Lemma}
\newtheorem{proposition}[lemma]{Proposition}
\newtheorem{corollary}[lemma]{Corollary}

\theoremstyle{definition}
\newtheorem{example}{Example}
\newtheorem{remark}{Remark}

\providecommand{\event}{QPL 2011} 
\title{The Expectation Monad \\ in Quantum Foundations}
\author{Bart Jacobs and Jorik Mandemaker 
\institute{Institute for Computing and Information Sciences (iCIS), \\
Radboud University Nijmegen, The Netherlands.}
\email{\{B.Jacobs,J.Mandemaker\}@cs.ru.nl} 
}

\def\titlerunning{The Expectation Monad}
\def\authorrunning{B. Jacobs \& J. Mandemaker}

\date{\small \today}

\begin{document}
\maketitle


\begin{abstract}
The expectation monad is introduced abstractly via two composable
adjunctions, but concretely captures measures. It turns out to sit in
between known monads: on the one hand the distribution and ultrafilter
monad, and on the other hand the continuation monad.  This expectation
monad is used in two probabilistic analogues of fundamental results of
Manes and Gelfand for the ultrafilter monad: algebras of the
expectation monad are convex compact Hausdorff spaces, and are dually
equivalent to so-called Banach effect algebras. These structures
capture states and effects in quantum foundations, and also the
duality between them.  Moreover, the approach leads to a new
re-formulation of Gleason's theorem, expressing that effects on a
Hilbert space are free effect modules on projections, obtained via
tensoring with the unit interval.
\end{abstract}
\renewcommand{\arraystretch}{1.3}
\setlength{\arraycolsep}{3pt}

\section{Introduction}\label{IntroSec}

Techniques that have been developed over the last decades for the
semantics of programming languages and programming logics gain wider
significance. In this way a new interdisciplinary area has emerged
where researchers from mathematics, (theoretical) physics and
(theoretical) computer science collaborate, notably on quantum
computation and quantum foundations. The article~\cite{BaezS10} uses
the phrase ``Rosetta Stone'' for the language and concepts of category
theory that form an integral part of this common area.

The present article is also part of this new field. It uses results
from programming semantics, topology and (convex) analysis, category
theory (esp.\ monads), logic and probability, and quantum foundations.
The origin of this article is an illustration of the connections
involved. Previously, the authors have worked on effect algebras and
effect modules~\cite{JacobsM12a,Jacobs10e,Jacobs11c} from quantum
logic, which are fairly general structures incorporating both logic
(Boolean and orthomodular lattices) and probability (the unit interval
$[0,1]$ and fuzzy predicates). By reading completely different work,
on formal methods in computer security (in particular the
thesis~\cite{ZanellaBeguelin10}), the expectation monad was noticed.
The monad is used in~\cite{ZanellaBeguelin10,BartheGZ09} to give
semantics to a probabilistic programming language that helps to
formalize (complexity) reduction arguments from security proofs in a
theorem prover.  In~\cite{ZanellaBeguelin10} (see
also~\cite{AudebaudP09,RamseyP02}) the expectation monad is defined in
a somewhat \textit{ad hoc} manner (see Section~\ref{SemanticsSec} for
details). Soon it was realized that a more systematic definition of
this expectation monad could be given via the (dual) adjunction
between convex sets and effect modules (elaborated in
Subsection~\ref{MndAdjSubsec}). Subsequently the two main parts of the
present paper emerged.
\begin{enumerate}
\item The expectation monad turns out to be related to several known
  monads as described in the following diagram.
\begin{equation}
\label{MonadsOverviewDiag}
\vcenter{\xymatrix@R-2pc{
\big(\mbox{distribution }\Dst\big)\ar@{ >->}[dr] \\
& \big(\mbox{expectation }\Exp\big)\ar@{ >->}[r] & 
   \big(\mbox{continuation }\mathcal{C}\big) \\
\big(\mbox{ultrafilter }\UF\big)\ar@{ >->}[ur]
}}
\end{equation}

\noindent The continuation monad $\mathcal{C}$ also comes from
programming semantics. But here we are more interested in the
connection with the distribution and ultrafilter monads $\Dst$ and
$\UF$. Since the algebras of the distribution monad are convex sets
and the algebras of the ultrafilter monad are compact Hausdorff spaces
(a result known as Manes theorem) it follows that the algebras of the
expectation monad must be some subcategory of convex compact Hausdorff
spaces. One of the main results in this paper,
Theorem~\ref{AlgobsCCHobsThm}, makes this connection precise. It can
be seen as a probabilistic version of Manes theorem. It uses basic
notions from Choquet theory, notably barycenters of measures.

\item The adjunction that gives rise to the expectation monad $\Exp$
  yields a (dual) adjunction between the category $\Alg(\Exp)$ of
  algebras and the category of effect modules. By suitable restriction
  this adjunction gives rise to an equivalence between
  ``observable'' $\Exp$-algebras and ``Banach'' (complete) effect
  modules, see Theorem~\ref{EModAlgExpDualityThm}.
\end{enumerate}

\noindent These two parts of the paper may be summarized as
follows. There are classical results:
$$\begin{array}{rcccl}
\Alg(\UF)
& \smash{\stackrel{\mbox{\scriptsize[Manes]}}{\simeq}} &
\big(\text{compact Hausdorff spaces}\big)
& \smash{\stackrel{\mbox{\scriptsize[Gelfand]}}{\simeq}} &
\big(\text{commutative } C^*\text{-algebras}\big)\op
\end{array}$$

\noindent Here we give the following ``probabilistic'' analogues:
$$\begin{array}{rcccl}
\Algobs(\Exp)
& \simeq &
\big(\text{convex compact Hausdorff spaces}\big)_{\mathrm{obs}}
& \simeq &
\big(\text{Banach effect modules}\big)\op
\end{array}$$

\noindent The subscript `obs' refers to a suitable observability
condition, see Section~\ref{ExpAlgSec}. The role played by the
two-element set $\{0,1\}$ in these classical
results---\textit{e.g.}~as ``schizophrenic'' object---is played in our
probabilistic analogues by the unit interval $[0,1]$.

Quantum mechanics is notoriously non-intuitive. Hence a proper
mathematical understanding of the relevant phenomena is important,
certainly within the emerging field of quantum computation. It seems
fair to say that such an all-encompassing understanding of quantum
mechanics does not exist yet. For instance, the categorical analysis
in~\cite{AbramskyC04,AbramskyC09} describes some of the basic
underlying structure in terms of monoidal categories, daggers, and
compact closure. However, an integrated view of logic and probability
is still missing. Here we certainly do not provide this integrated
view, but possibly we do contribute a bit. The states of a Hilbert
space $\H$, described as density matrices $\DM(\H)$, fit within the
category of convex compact Hausdorff spaces investigated here.  Also,
the effects $\Ef(\H)$ of the space fit in the associated dual category
of Banach Hausdorff spaces. The duality we obtain between convex
compact Hausdorff spaces and Banach effect algebras precisely captures
the translations back and forth between states and effects, as
expressed by the isomorphisms:
$$\begin{array}{rclcrcl}
\mathrm{Hom}\big(\Ef(\H), [0,1]) & \cong & \DM(H)
& \qquad &
\mathrm{Hom}\big(\DM(\H), [0,1]) & \cong & \Ef(H).
\end{array}$$

\noindent These isomorphisms (implicitly) form the basis for the
quantum weakest precondition calculus described in~\cite{dHondtP06}.

In this context we shed a bit more light on the relation between
quantum logic---as expressed by the projections $\pr(\H)$ on a Hilbert
space---and quantum probability---via its effects $\Ef(\H)$. In
Section~\ref{GleasonSec} it will be shown that Gleason's famous
theorem, expressing that states are probability measures, can
equivalently be expressed as an isomorphism relating projections
and effects:
$$\begin{array}{rcl}
[0,1]\otimes \pr(\H) & \cong & \Ef(\H).
\end{array}$$

\noindent This means that the effects form the free effect module on
projections, via the free functor $[0,1]\otimes(-)$. More loosely
formulated: quantum probabilities are freely obtained from quantum
predicates.

We briefly describe the organization of the paper. It starts with a
quick recap on monads in Section~\ref{MonadSec}, including
descriptions of the monads relevant in the rest of the paper.
Section~\ref{EModSec} gives a brief introduction to effect algebras
and effect modules. It also establishes equivalences between (Banach)
order unit spaces and (Banach) Archimedean effect modules.  In
Section~\ref{ExpMndSec} we give several descriptions of the
expectation monad in terms of effect algebras and effect modules. We
also describe the map between the expectation monad and the
continuation monad here. Sections~\ref{ExpUFSec} and~\ref{ExpDstSec}
deal with the construction of the other two monad maps from
Diagram~\eqref{MonadsOverviewDiag}: those from the ultrafilter and
distribution monads to the expectation monad. Here we also explore
some of the implications of these maps. Next, in
Section~\ref{ExpAlgSec}, we study the algebras of the expectation
monad.  We prove that the category of $\Exp$-algebras is equivalent to
the category compact convex sets with continuous affine mappings.  In
Section~\ref{ExpEModSec} we establish a dual adjunction between
$\Exp$-algebras and effect modules. We prove that when restricted to
so-called observable $\Exp$-algebras and Banach effect modules this
adjunction becomes an equivalence. In Section~\ref{GleasonSec} we
apply this duality to quantum logic. We prove that the isomorphism
$[0,1]\otimes\Pr(\H)\cong\Ef(\H)$ is an algebraic reformulation of
Gleason's theorem.  Finally in Section~\ref{SemanticsSec} we examine
how the expectation monad has appeared in earlier work on programming
semantics. We also suggest how it might be used to capture both
non-deterministic and probabilistic computation simultaneously,
although the details of this are left for future work.



\section{A recap on monads}\label{MonadSec}

This section recalls the basics of the theory of monads, as needed
here. For more information, see
\textit{e.g.}~\cite{MacLane71,BarrW85,Manes74,Borceux94}. Some
specific examples will be elaborated later on. 

A monad is a functor $T\colon\cat{C} \rightarrow \cat{C}$ together
with two natural transformations: a unit $\eta\colon \idmap[\cat{C}]
\Rightarrow T$ and multiplication $\mu \colon T^{2} \Rightarrow
T$. These are required to make the following diagrams commute, for
$X\in\cat{C}$.
$$\xymatrix@C-.5pc@R-.5pc{
T(X)\ar[rr]^-{\eta_{T(X)}}\ar@{=}[drr] & & T^{2}(X)\ar[d]^{\mu_X} & &
   T(X)\ar[ll]_{T(\eta_{X})}\ar@{=}[dll] 
&  & T^{3}(X)\ar[rr]^-{\mu_{T(X)}}\ar[d]_{T(\mu_{X})} 
   & & T^{2}(X)\ar[d]^{\mu_X} \\
& & T(X) & &
& &
T^{2}\ar[rr]_{\mu_X} & & T(X)
}$$

\noindent Each adjunction $F\dashv G$ gives rise to a monad $GF$.

Given a monad $T$ one can form a category $\Alg(T)$ of so-called
(Eilenberg-Moore) algebras. Objects of this category are maps of
the form $a\colon T(X) \rightarrow X$, making the first two squares
below commute.
$$\xymatrix@R-.5pc{
X\ar@{=}[dr]\ar[r]^-{\eta} & TX\ar[d]^{a}
& 
T^{2}X\ar[d]_{\mu}\ar[r]^-{T(a)} & TX\ar[d]^{a}
& &
TX\ar[d]_{a}\ar[r]^-{T(f)} & TY\ar[d]^{b} \\
& X 
&
TX\ar[r]_-{a} & X
& &
X\ar[r]_-{f} & Y
}$$

\noindent A homomorphism of algebras $(X,a) \rightarrow (Y,b)$ is a
map $f\colon X\rightarrow Y$ in $\cat{C}$ between the underlying
objects making the diagram above on the right commute. The diagram in
the middle thus says that the map $a$ is a homomorphism $\mu
\rightarrow a$. The forgetful functor $U\colon \Alg(T) \rightarrow
\cat{C}$ has a left adjoint, mapping an object $X\in\cat{X}$ to the
(free) algebra $\mu_{X}\colon T^{2}(X) \rightarrow T(X)$ with carrier
$T(X)$.

Each category $\Alg(T)$ inherits limits from the category $\cat{C}$.
In the special case where $\cat{C}=\Sets$, the category of sets and
functions (our standard universe), the category $\Alg(T)$ is not only
complete but also cocomplete (see~\cite[\S~9.3, Prop.~4]{BarrW85}).

A map of monads $\sigma\colon T\Rightarrow S$ is a natural
transformation that commutes with the units and multiplications, as in:
\begin{equation}
\label{MndMapEqn}
\vcenter{\xymatrix@R-.5pc@C-.2pc{
X\ar[d]_{\eta_X}\ar@{=}[r] & X\ar[d]^{\eta_X}
& & 
T^{2}(X)\ar[d]_{\mu_X}\ar[r]^-{\sigma_{TX}} & 
   S(T(X))\ar[r]^-{S(\sigma_{X})} &
   S^{2}(X)\ar[d]^{\mu_{X}} \\
T(X)\ar[r]_-{\sigma_X} & S(X)
& &
T(X)\ar[rr]_-{\sigma_X} & & S(X)
}}
\end{equation}

\noindent Such a map of monads $\sigma\colon T\Rightarrow S$ induces a
functor $(-)\after\sigma \colon \Alg(S) \rightarrow \Alg(T)$ between
categories of algebras that commutes with the forgetful functors.

\begin{lemma}
\label{MndMapAlgLem}
Assume a map of monads $\sigma\colon T\Rightarrow S$.
\begin{enumerate}
\item There is a functor $(-)\after\sigma \colon \Alg(S) \rightarrow
  \Alg(T)$ that commutes with the forgetful functors.

\item If the category $\Alg(S)$ has sufficiently many
  coequalizers---like when the underlying category is $\Sets$---this
  functor has a left adjoint $\Alg(T)\rightarrow\Alg(S)$; it maps an
  algebra $a\colon T(X)\rightarrow X$ to the following coequalizer
  $a_{\sigma}$ in $\Alg(S)$.
$$\vcenter{\xymatrix@C-.5pc{
\ensuremath{\left(\xy
(0,4)*{S^{2}(TX)};
(0,-4)*{S(TX)};
{\ar^{\mu} (0,2); (0,-2)};
\endxy\right)}
\ar@<.5pc>[rr]^-{\mu \after S(\sigma)}
\ar@<-.5pc>[rr]_-{S(a)}
& &
\ensuremath{\left(\xy
(0,4)*{S^{2}(X)};
(0,-4)*{S(X)};
{\ar^{\mu} (0,2); (0,-2)};
\endxy\right)}
\ar@{->>}[rr]^-{c} & & 
\ensuremath{\left(\xy
(0,4)*{S(X_{\sigma})};
(0,-4)*{X_{\sigma}};
{\ar^{a_{\sigma}} (0,2); (0,-2)};
\endxy\right)}
}}\eqno{\QEDbox}$$
\end{enumerate}
\end{lemma}

\begin{myproof}
We need to establish a bijective correspondence between algebra maps:
$$\begin{prooftree}
{\xymatrix@C+1pc{
\ensuremath{\left(\xy
(0,4)*{S(X_{\sigma})};
(0,-4)*{X_{\sigma}};
{\ar^{a_{\sigma}} (0,2); (0,-2)};
\endxy\right)}
\ar[r]^-{f} &
\ensuremath{\left(\xy
(0,4)*{S(Y)};
(0,-4)*{Y};
{\ar^{b} (0,2); (0,-2)};
\endxy\right)}}}
\Justifies
{\xymatrix@C+1pc{
\ensuremath{\left(\xy
(0,4)*{T(X)};
(0,-4)*{X};
{\ar^{a} (0,2); (0,-2)};
\endxy\right)}
\ar[r]_-{g} &
\ensuremath{\left(\xy
(0,4)*{T(Y)};
(0,-4)*{Y};
{\ar^{b\after\sigma} (0,2); (0,-2)};
\endxy\right)}}}
\end{prooftree}$$

\noindent This works as follows. Given $f$, one takes $\overline{f} =
f \after c \after \eta\colon X\rightarrow Y$. And given $g$ one
obtains $\overline{g}\colon X_{\sigma} \rightarrow Y$ because $b
\after T(g)\colon S(X)\rightarrow Y$ coequalizes the above parallel
pair $\mu \after S(\sigma)$ and $S(a)$. Remaining details are left to
the interested reader. \QED

\auxproof{
First of all, $\overline{f} = f \after c \after \eta$ is a map of
$T$-algebras:
$$\begin{array}{rcl}
b \after \sigma \after T(\overline{f})
& = &
b \after \sigma \after T(f \after c \after \eta) \\
& = &
b \after S(f) \after S(c) \after S(\eta) \after \sigma \\
& = &
f \after a_{\sigma} \after S(c) \after S(\eta) \after \sigma \\
& = &
f \after c \after \mu \after S(\eta) \after \sigma \\
& = &
f \after c \after \mu \after \eta \after \sigma \\
& = &
f \after c \after \mu \after S(\sigma) \after \eta \\
& = &
f \after c \after S(a) \after \eta \\
& = &
f \after c \after \eta \after a.
\end{array}$$

\noindent Starting from $g$ we have a map of algebras $b \after S(g)$
from $\mu$ to $b$, since:
$$\begin{array}{rcccl}
b \after S(b \after S(g))
& = &
b \after \mu \after S^{2}(g)
& = &
b \after S(g) \after \mu.
\end{array}$$

\noindent Moreover, this map $b \after S(g)$ coequalizes the above
parallel pair $\mu \after S(\sigma)$ and $S(a)$ since:
$$\begin{array}{rcl}
b \after S(g) \after \mu \after S(\sigma)
& = &
b \after \mu \after S^{2}(g) \after S(\sigma) \\
& = &
b \after \mu \after S(\sigma) \after S(T(g)) \\
& = &
b \after S(b \after \sigma) \after S(T(g)) \\
& = &
b \after S(g) \after S(a).
\end{array}$$

\noindent Hence there is a unique map of algebras $\overline{g}
\colon a_{\sigma} \rightarrow b$ with $\overline{g} \after c = 
b \after S(g)$.

Then:
$$\begin{array}{rcl}
\overline{\overline{f}} \after c
& = &
b \after S(\overline{f}) \\
& = &
b \after S(f) \after S(c) \after S(\eta) \\
& = &
f \after a_{\sigma} \after S(c) \after S(\eta) \\
& = &
f \after c \after \mu \after S(\eta) \\
& = &
f \after c \\
\overline{\overline{g}}
& = &
\overline{g} \after c \after \eta \\
& = &
b \after S(g) \after \eta \\
& = &
b \after \eta \after g \\
& = &
g.
\end{array}$$

The unit of this adjunction is the unique map of algebras:
$$\xymatrix@C+1pc{
\ensuremath{\left(\xy
(0,4)*{T(X)};
(0,-4)*{X};
{\ar^{a} (0,2); (0,-2)};
\endxy\right)}
\ar[r]^-{c \after \eta} &
\ensuremath{\left(\xy
(0,4)*{T(X_{\sigma})};
(0,-4)*{X_{\sigma}};
{\ar^{a_{\sigma}\after\sigma} (0,2); (0,-2)};
\endxy\right)}}$$

\noindent It is a map of algebras since:
$$\begin{array}{rcl}
a_{\sigma} \after \sigma \after T(c \after \eta) 
& = &
a_{\sigma} \after S(c) \after S(\eta) \after \sigma \\
& = &
c \after \mu \after S(\eta) \after \sigma \\
& = &
c \after \mu \after \eta \after \sigma \\
& = &
c \after \mu \after S(\sigma) \after \eta \\
& = &
c \after S(a) \after \eta \\
& = &
c \after \eta \after a.
\end{array}$$
}

\end{myproof}

\subsection{The Distribution monad}\label{DstSubsec}

We shall write $\Dst$ for the discrete probability distribution monad
on $\Sets$. It maps a set $X$ to the set of formal convex combinations
$r_{1}x_{1} + \cdots + r_{n}x_{n}$, where $x_{i}\in X$ and
$r_{i}\in[0,1]$ with $\sum_{i}r_{1} = 1$. Alternatively,
$$\begin{array}{rcl}
\Dst(X)
& = &
\set{\varphi\colon X\rightarrow [0,1]}{\supp(\varphi) 
   \mbox{ is finite, and } \sum_{x}\varphi(x) = 1},
\end{array}$$

\noindent where $\supp(\varphi)\subseteq X$ is the support of
$\varphi$, containing all $x$ with $\varphi(x) \neq 0$.  The functor
$\Dst\colon\Sets\rightarrow\Sets$ forms a monad with the Dirac
function as unit in:
$$\xymatrix@R-2pc{
X\ar[r]^-{\eta} & \Dst{X}
&
\Dst{\Dst{X}}\ar[r]^-{\mu} & \Dst{X} \\
x\ar@{|->}[r] & 1x = 
   {\lam{y}{\left\{\begin{array}{ll}
      1 & \mbox{if }y=x \\ 0 & \mbox{if }y\neq x \end{array}\right.}}
&
\Psi\ar@{|->}[r] & 
   \lam{y}{\sum_{\varphi\in\Dst{X}}\Psi(\varphi)\cdot\varphi(y)}.
}$$

\noindent [Here we use the ``lambda'' notation from the lambda
calculus~\cite{Barendregt84}: the expression $\lam{x}{\cdots}$ is used
for the function $x \mapsto \cdots$. We also use the associated
application rule $(\lam{x}{f(x)})(y) = f(y)$.]

Objects of the category $\Alg(\Dst)$ of (Eilenberg-Moore)
algebras of this monad $\Dst$ can be identified as \emph{convex sets},
in which sums $\sum_{i}r_{i}x_{i}$ of convex combinations exists.
Morphisms are so-called affine functions, preserving such convex sums,
see~\cite{Jacobs10e}. Hence we also write $\Alg(\Dst) = \Conv$, where
$\Conv$ is the category of convex sets and affine functions.

The prime example of a convex set is the unit interval $[0,1]\subseteq
\mathbb{R}$ of probabilities. Also, for an arbitrary set $X$, the set
of functions $[0,1]^{X}$, or fuzzy predicates on $X$, is a convex set,
via pointwise convex sums.

\subsection{The ultrafilter monad}\label{UFSubsec}

A particular monad that plays an important role in this paper is
the ultrafilter monad $\UF\colon\Sets\rightarrow\Sets$, given by:
\begin{equation}
\label{UFDefEqn}
\hspace*{-.5em}\begin{array}{rcl}
\UF(X)
& = &
\set{\calF \subseteq \Pow(X)}{\calF 
   \mbox{ is an ultrafilter}} \\
& \cong &
\set{f\colon \Pow(X) \rightarrow \{0,1\}}{f 
   \mbox{ is a homomorphism of Boolean algebras}}
\end{array}
\end{equation}

\noindent Such an ultrafilter $\calF\subseteq \Pow(X)$ satisfies,
by definition, the following three properties.
\begin{itemize}
\item It is an upset: $V\supseteq U\in\calF \Rightarrow V\in\calF$;

\item It is closed under finite intersections: $X\in\cal F$ and 
$U,V\in\calF \Rightarrow U\cap V\in\calF$;

\item For each set $U$ either $U\in\calF$ or $\neg U =
  \setin{x}{X}{x\not\in U} \in\calF$, but not both. As a consequence,
  $\emptyset \not\in\calF$.
\end{itemize}

\noindent For a function $f\colon X\rightarrow Y$ one obtains
$\UF(f)\colon \UF(X) \rightarrow \UF(Y)$ by:
$$\begin{array}{rcl}
\UF(f)(\calF)
& = &
\set{V\subseteq Y}{f^{-1}(V)\in\calF}.
\end{array}$$

\noindent Taking ultrafilters is a monad, with unit $\eta\colon 
X\rightarrow \UF(X)$ given by so-called principal ultrafilters:
$$\begin{array}{rcl}
\eta(x)
& = &
\set{U\subseteq X}{x\in U}.
\end{array}$$

\noindent The multiplication $\mu\colon\UF^{2}(X)\rightarrow\UF(X)$ is:
$$\begin{array}{rclcrcl}
\mu(\mathcal{A})
& = &
\set{U\subseteq X}{D(U)\in\mathcal{A}}
& \quad\mbox{where}\quad
D(U)
& = &
\setin{\calF}{\UF(X)}{U\in\calF}.
\end{array}$$

\auxproof{
$$\begin{array}{rcl}
\big(\mu \after \eta\big)(\cal F)
& = &
\mu\big(\set{V\subseteq \UF(X)}{\calF\in V}\big) \\
& = &
\set{U\subseteq X}{\calF \in 
   \setin{\mathcal{G}}{\UF(X)}{U\in\mathcal{G}}} \\
& = &
\set{U\subseteq X}{U\in\calF} \\
& = &
\calF \\
\big(\mu \after \UF(\eta)\big)(\calF) 
& = &
\mu\big(\set{V\subseteq\UF(X)}{\eta^{-1}(V)\in\calF}\big) \\
& = &
\set{U\subseteq X}{\eta^{-1}(\setin{\mathcal{G}}{\UF(X)}{U\in\mathcal{G}})
   \in\calF} \\
& = &
\set{U\subseteq X}{\set{x}{U\in\eta(x)}\in\calF} \\
& = &
\set{U\subseteq X}{\set{x}{x\in U}\in\calF} \\
& = &
\set{U\subseteq X}{U\in\calF} \\
& = &
\calF \\
\big(\mu \after \UF(\mu)\big)(\mathcal{B}) 
& = &
\mu\big(\set{V\subseteq\UF(X)}{\mu^{-1}(V)\in\mathcal{B}}\big) \\
& = &
\set{U\subseteq X}{\setin{\mathcal{G}}{\UF(X)}{U\in\mathcal{G}} \in
   \set{V\subseteq\UF(X)}{\mu^{-1}(V)\in\mathcal{B}}} \\
& = &
\set{U\subseteq X}{\mu^{-1}(\setin{\mathcal{G}}{\UF(X)}{U\in\mathcal{G}}) 
   \in \mathcal{B}} \\
& = &
\set{U\subseteq X}{\setin{\mathcal{A}}{\UF^{2}(X)}{
   U\in\mu(\mathcal{A})}\in \mathcal{B}} \\
& = &
\set{U\subseteq X}{\setin{\mathcal{A}}{\UF^{2}(X)}
   {\setin{\mathcal{G}}{\UF(X)}{U\in\mathcal{G}}
      \in\mathcal{A}}\in \mathcal{B}} \\
& = &
\mu\big(\set{V\subseteq\UF(X)}{
   \setin{\mathcal{A}}{\UF^{2}(X)}{V\in\mathcal{A}}\in \mathcal{B}} \\
& = &
\big(\mu \after \mu\big)(\mathcal{B}).
\end{array}$$
}

The set $\UF(X)$ of ultrafilters on a set $X$ is a topological space
with basic (compact) clopens given by subsets $D(U) =
\setin{\calF}{\UF(X)}{U\in\calF}$, for $U\subseteq X$. This makes
$\UF(X)$ into a compact Hausdorff space. The unit $\eta\colon
X\rightarrow \UF(X)$ is a dense embedding.

\auxproof{ 
In order to show that $A\subseteq Y$ is dense,
  \textit{i.e.}~$\overline{A} = Y$, one can prove that each non-empty
open $U\subseteq Y$ has non-empty intersection with $A$. If not,
\textit{i.e.}~$A\cap U = \emptyset$, then $A\subseteq \neg U$,
so that $X = \overline{A} \subseteq \neg U$, yielding $\neg U = X$ 
and thus $U=\emptyset$.

Hence, take a non-empty open $D(U)\subseteq \UF(X)$. Then
$U\neq\emptyset$, since $D(\emptyset) = \emptyset$. Pick $x\in U$;
then $U\in\eta(x)$, so $D(U) \cap \set{\eta(x)}{x\in X} \neq
\emptyset$.

Maps $\UF(f) \colon \UF(X) \rightarrow UF(Y)$ are continuous, since
for $V\subseteq Y$ one has $\UF(f)^{-1}(D(V)) = D(f^{-1}(V))$.
}

The following result shows the importance of the ultrafilter monad, see
\textit{e.g.}~\cite{Manes69}, \cite[III.2]{Johnstone82},
or~\cite[Vol.~2, Prop.~4.6.6]{Borceux94}.

\begin{theorem}[Manes]
\label{ManesThm}
$\Alg(\UF) \simeq \CH$, \textit{i.e.} the category of algebras of the
  ultrafilter monad is equivalent to the category $\CH$ of compact
  Hausdorff spaces and continuous maps.
\end{theorem}

The proof is complicated and will not be reproduced here. We only
extract the basic constructions. For a compact Hausdorff space $Y$ one
uses denseness of the unit $\eta$ to define a unique continuous
extensions $f^\#$ as in:
\begin{equation}
\label{UFextensionDiag}
\vcenter{\xymatrix@R-.5pc{
X\ar@{ >->}[rr]^-{\eta}\ar[drr]_{f} & & \UF(X)\ar@{-->}[d]^{f^\#} \\
& & Y
}}
\end{equation}

\noindent One defines $f^{\#}(\calF)$ to be the unique element in
$\bigcap\set{\overline{V}}{V\subseteq Y \mbox{ with }f^{-1}(V) \in
  \calF}$. This intersection is a singleton precisely because $Y$ is a
compact Hausdorff space. In such a way one obtains an algebra
$\UF(Y)\rightarrow Y$ as extension of the identity.

\auxproof{
Recall the basic extension results about dense embeddings.
\begin{enumerate}
\item If $i\colon A\hookrightarrow X$ is a dense embedding and
  continuous maps $f,g\colon X \rightarrow Y$, where $Y$ is Hausdorff,
  satisfy $f \after i = g \after i$, then $f=g$.

Suppose not: $f(x)\neq g(x)$. Find disjoint opens $U\ni f(x)$ and
$V\ni g(x)$. Then $f^{-1}(U) \cap g^{-1}(V) \subseteq X$ is non-empty
(containing $x$) and open, so $f^{-1}(U) \cap g^{-1}(V)\cap A \neq
\emptyset$, say with $z \in f^{-1}(U) \cap g^{-1}(V)\cap A$. But then
$f(z) = g(z) \in U \cap V$, contradicting the assumption.

\item Define $f^{\#}(\calF)$ to be the unique element in
  $\bigcap\set{\overline{V}}{V\subseteq Y \mbox{ with }f^{-1}(V) \in
  \calF}$. We first convince ourselves that this intersection is
  non-empty. Suppose not,
  \textit{i.e.}~$\bigcap\set{\overline{V}}{V\subseteq Y \mbox{ with
    }f^{-1}(V) \in \calF} = \emptyset$. Then its complement must be
  the whole space, \textit{i.e.} $Y =
  \bigcup\set{\neg\overline{V}}{V\subseteq Y \mbox{ with }f^{-1}(V)
    \in \calF}$. Since $Y$ is compact there are finitely many $V_{1},
  \ldots, V_{n}\subseteq Y$ with $f^{-1}(V_{i}) \in \calF$ and
  $Y = \neg \overline{V_{1}} \cup \cdots \cup
  \neg\overline{V_n}$. Thus $V_{1} \cap \cdots \cap V_{n} \subseteq
  \overline{V_{1}} \cap \cdots \cap \overline{V_{n}} = \neg\big(\neg
  \overline{V_{1}} \cup \cdots \cup \neg\overline{V_n}\big) = \neg Y =
  \emptyset$. But then $\emptyset = f^{-1}(\emptyset) = f^{-1}(V_{1}
  \cap \cdots \cap V_{n}) = f^{-1}(V_{1}) \cap \cdots \cap
  f^{-1}(V_{n}) \in \calF$. This is impossible.

Next suppose the intersection $\bigcap\set{\overline{V}}{V\subseteq Y
  \mbox{ with }f^{-1}(V) \in \calF}$ contains two (distinct) elements
$y_{1},y_{2}$ Since $Y$ is Hausdorff there are two disjoint opens
$V_{i}\subseteq Y$ with $y_{i}\in V_{i}$. We have either
$f^{-1}(V_{1})\in\calF$ or $f^{-1}(V_{1})\not\in\calF$.
\begin{itemize}
\item If $f^{-1}(V_{1})\in\calF$ then
  $y_{1},y_{2}\in\overline{V_{1}}$.  But since $V_{1}\cap V_{2} =
  \emptyset$, we have $V_{1} \subseteq \neg V_{2}$, which is closed,
  and thus $\overline{V_{1}} \subseteq \neg V_{2}$. Now we get $y_{2}
  \in \neg V_{2}$, which is impossible.

\item If $f^{-1}(V_{1})\not\in\calF$ then $\neg f^{-1}(V_{1}) =
  f^{-1}(\neg V_{1}) \in\calF$. Hence $y_{1},y_{2}\in \overline{\neg
    V_1} = \neg V_{1}$, which is impossible.
\end{itemize}

\noindent We get $f^{\#} \after \eta = f$ since:
$$\begin{array}{rcl}
f(x)=y
& \Longleftrightarrow &
y\in \bigcap\set{\overline{V}}{V\subseteq Y \mbox{ with }x\in f^{-1}(V)}.
\end{array}$$

\noindent Clearly, $f(x)$ is in the intersection on the right. This
gives the implication $(\Rightarrow)$. But it also gives $(\Leftarrow)$
because there is precisely one element in this intersection.

What remains is to show that $f^{\#} \colon \UF(X) \rightarrow Y$
is continuous. So assume $W\subseteq Y$ is open. We claim:
$$\begin{array}{rcl}
\big(f^{\#}\big)^{-1}(W)
& = &
\setin{\calF}{\UF(X)}{\ex{V\subseteq Y}{\overline{V} \subseteq W
   \mbox{ and }f^{-1}(V) \in \calF}}.
\end{array}$$

\noindent Then we are done, since the right-hand-side can be written
as union of opens: $\bigcup_{\overline{V}\subseteq W}D(f^{-1}(V))$.
\begin{itemize}
\item[$(\subseteq)$] Assume $\calF \in (f^{\#})^{-1}(W)$. We
  use regularity of the compact Hausdorff space $Y$, in the form: each
  open $W\subseteq Y$ satisfies $W =
  \bigcup\set{V\in\open(Y)}{\overline{V} \subseteq W}$. Thus
  $f^{\#}(\calF) \in V$ for some open $V$ with $\overline{V}\subseteq
  W$.  We wish to show $f^{-1}(V)\in\calF$. Suppose not; then
  $f^{-1}(\neg V) \in \calF$ and thus $f^{\#}(\calF) \in
  \overline{\neg V} = \neg V$. But this contradicts the assumption
  $f^{\#}(\calF) \in V$.

\item[$(\supseteq)$] Assume we have an $\calF$ with a subset
  $V\subseteq Y$ satisfying $\overline{V} \subseteq W$ and $f^{-1}(V)
  \in \calF$. But then $f^{\#}(\calF) \in \overline{V} \subseteq W$,
  so $\calF\in (f^{\#})^{-1}(W)$.
\end{itemize}
\end{enumerate}
}



Conversely, assuming an algebra $\ch_{X}\colon\UF(X)\rightarrow X$ one
defines $U\subseteq X$ to be closed if for all $\calF\in\UF(X)$,
$U\in\calF$ implies $\ch(\calF)\in U$. 
This yields a topology on $X$ which is Hausdorff and
compact. There can be at most one such algebra structure
$\ch_{X}\colon\UF(X)\rightarrow X$ on a set $X$, corresponding to a
compact Hausdorff topology, because of the following standard
result.

\begin{lemma}
\label{UniqueCHLem}
Assume a set $X$ carries two topologies $\open_{1}(X), \open_{2}(X)
\subseteq \Pow(X)$ with $\open_{1}(X) \subseteq \open_{2}(X)$,
$\open_{1}(X)$ is Hausdorff and $\open_{2}(X)$ is compact, then
$\open_{1}(X) = \open_{2}(X)$. \QED
\end{lemma}

\begin{myproof}
If $U$ is closed in $\open_{2}(X)$, then it is compact, and, because
$\open_{1}(X) \subseteq \open_{2}(X)$, also compact in
$\open_{1}(X)$. Hence it is closed there. \QED
\end{myproof}

\auxproof{
Because:
$$\begin{array}{rcl}
U\in\open_{2}(X)
& \Longrightarrow &
\neg U \mbox{ is closed in } \open_{2}(X) \\
& \smash{\stackrel{(1)}{\Longrightarrow}} &
\neg U \mbox{ is compact in } \open_{2}(X) \\
& \smash{\stackrel{(2)}{\Longrightarrow}} &
\neg U \mbox{ is compact in } \open_{1}(X) \\
& \smash{\stackrel{(3)}{\Longrightarrow}} &
\neg U \mbox{ is closed in } \open_{1}(X) \\
& \Longrightarrow &
U\in\open_{1}(X).
\end{array}$$

\noindent We briefly review these steps.
\begin{enumerate}
\item Standard: an open covering $V \subseteq \bigcup_{i}W_{i}$ of a
  closed set $V$ yields an open cover $\neg V \cup \bigcup_{i}W_{i}$
  of the whole space, and thus a finite subcover.

\item Suppose $V$ is compact in $\open_{2}(X)$ and $V \subseteq
  \bigcup_{i}W_{i}$ with $W_{i}\in\open_{1}(X)$. Then $W_{i}\subseteq
  \open_{2}(X)$. Thus there is a finite subcover.

\item Suppose $V$ is compact and $x\not\in V$. Then there is an an
  open set $U_{x}$ with $x\in U_{x}$ and $V\cap U_{x} = \emptyset$.
  (Subproof: for each $y\in V$ there are disjoint opens $V_{y}\ni y$
  and $W_{y}\ni x$. Then $V \subseteq \bigcup_{y\in V}V_{y}$. Hence $V
  \subseteq V_{y_{1}} \cup \cdots \cup V_{y_n}$. And $x\in W_{y_1}
  \cap \cdots \cap W_{y_n}$, where $W_{y_1} \cap \cdots \cap W_{y_n}$
  is open and disjoint from $V$.) Thus $\neg V = \bigcup_{x\in\neg
    V}U_{x}$ is a union of opens.
\end{enumerate}
}

We can apply this result to the space $\UF(X)$ of ultrafilters: as
described before Theorem~\ref{ManesThm}, $\UF(X)$ carries a compact
Hausdorff topology with sets $D(U) = \setin{\calF}{\UF(X)}{U\in\calF}$
as clopens. Also, it carries a compact Hausdorff topology via the
(free) algebra $\mu_{X} \colon \UF^{2}(X) \rightarrow \UF(X)$.  It is
not hard to see that the subsets $D(U)$ are closed in the latter
topology, so the two topologies on $\UF(X)$ coincide by
Lemma~\ref{UniqueCHLem}. Later we shall use a similar argument.

\auxproof{ 
The closure $\overline{V}$ of a subset $V\subseteq\Exp(X)$ is
$\overline{V} = \set{\mu(\mathcal{A})}{\mathcal{A}\in\UF^{2}(X) \mbox{
    with }V\in\mathcal{A}}$. Then, subsets $D(U)\subseteq\Exp(X)$, for
$U\subseteq X$, are closed in this topology, since:
$$\begin{array}{rcl}
\overline{D(U)}
& = &
\set{\mu(\mathcal{A})}{D(U)\in\mathcal{A}} \\
& = &
\set{\mu(\mathcal{A})}{U\in\mu(\mathcal{A})} \\
& = &
\setin{\calF}{\UF(X)}{U\in\calF},
\end{array}$$

\noindent where the inclusion $(\subseteq)$ of the last equality
is obvious; for $(\supseteq)$ one takes $\mathcal{A} = \eta(\calF)$,
so that $\mu(\mathcal{A}) = \calF$.
}

\begin{example}
\label{UnitIntervalUFEx}
The unit interval $[0,1]\subseteq\mathbb{R}$ is a standard example of
a compact Hausdorff space. Its Eilenberg-Moore algebra
$\ch\colon\UF([0,1]) \rightarrow [0,1]$ can be described concretely on
$\calF\in\UF([0,1])$ as:
\begin{equation}
\label{UnitIntervalUFEqn}
\begin{array}{rcl}
\ch(\calF)
& = &
\inf\setin{s}{[0,1]}{[0,s]\in\calF}.
\end{array}
\end{equation}

\noindent For the proof, recall that $\ch(\calF)$ is the (sole)
element of the intersection $\bigcap\set{\overline{V}}{V\in\calF}$.
Hence if $[0,s]\in\calF$, then $\ch(\calF)\in\overline{[0,s]} =
[0,s]$, so $\ch(\calF) \leq s$. This establishes the $(\leq)$-part
of~\eqref{UnitIntervalUFEqn}. Assume next that $\ch(\calF) <
\inf\set{s}{[0,s]\in\calF}$. Then there is some $r\in[0,1]$ with
$\ch(\calF) < r < \inf\set{s}{[0,s]\in\calF}$. Then $[0,r]$ is not in
$\calF$, so that $\neg[0,r] = (r,1]\in\calF$.  But this means
  $\ch(\calF) \in \overline{(r,1)} = [r,1]$, which is impossible.

Notice that~\eqref{UnitIntervalUFEqn} can be strengthened to:
$\ch(\calF) = \inf\setin{s}{[0,1]\cap\mathbb{Q}}{[0,s]\in\calF}$.
\end{example}

The second important result about compact Hausdorff spaces is
as follows.

\begin{theorem}[Gelfand]
\label{GelfandThm}
$\CH \simeq C^{*}\mbox{-}\Cat{Alg}\op$, \textit{i.e.} the category
$\CH$ of compact Hausdorff spaces is equivalent to the opposite
of the category of commutative $C^*$-algebras.
\end{theorem}

This paper presents probabilistic analogues of these two basic results
(Theorems~\ref{ManesThm} and~\ref{GelfandThm}), involving
\emph{convex} compact Hausdorff spaces (see
Theorem~\ref{EModAlgExpDualityThm}).

\auxproof{
\begin{center}
\textbf{Elaboration of Manes result, from a different context.}
\end{center}

The following easy result is worth making explicit.

\begin{lemma}
\label{UltraFilterEqLem}
For maps of Boolean algebras $\varphi,\psi\colon A\rightarrow B$, in
order to prove $\varphi=\psi$ it suffices to prove $\varphi\leq \psi$.
\end{lemma}

\begin{myproof}
Assume $\varphi\leq\psi$, \textit{i.e.}~$\varphi(a) \leq \psi(a)$ for
each $a\in A$. Then also: $\neg\varphi(a) = \varphi(\neg a) \leq
\psi(\neg a) = \neg \psi(a)$, so that $\varphi(a) \geq \psi(a)$. Hence
$\varphi(a) = \psi(a)$. \QED
\end{myproof}

The following ``Ultrafilter Lemma'' follows from the Axiom of Choice,
see~\cite[I, 2.3 and 2.4]{Johnstone82}. It is crucial in this context.

\begin{lemma}
\label{UltraFilterLem}
Let $B$ be a Boolean algebra containing a filter $F\subseteq B$ that
is disjoint from an ideal $I\subseteq B$. Then there is a homomorphism
of Boolean algebras $\varphi\colon B\rightarrow 2$ with $F\subseteq
\varphi^{-1}(1)$ and $I\subseteq \varphi^{-1}(0) =
\ker(\varphi)$. \QED
\end{lemma}

For an arbitrary subset $U\subseteq B$ of a Boolean algebra $B$ we
write
$$\begin{array}{rcl}
\ideal U
& = &
\downset\set{\bigvee V}{V\mbox{ is a finite subset of }U} \\
\filter U
& = &
\upset\set{\bigwedge V}{V\mbox{ is a finite subset of }U}
\end{array}$$

\noindent for the least ideal and filter generated by $U$.

We define a functor $\UF\colon\Sets\rightarrow\Sets$ as:
$$\begin{array}{rcl}
\UF(X)
& = &
\BA(\Pow(X), 2) \\
\UF\big(X\stackrel{f}{\rightarrow} Y\big)(\varphi)
& = &
\Big(\Pow(Y)\stackrel{f^{-1}}{\rightarrow} 
   \Pow(X)\stackrel{\varphi}{\rightarrow} 2\Big).
\end{array}$$

\noindent Later we shall see that $\UF$ is a monad, see
Theorem~\ref{CHAdjThm}. At this stage we only define a unit
$\eta\colon \idmap\Rightarrow \UF$ as:
\begin{equation}
\label{UFunitEqn}
\begin{array}{rcl}
\eta_{X}(x)(U)
& = &
\left\{\begin{array}{ll}
1 & \mbox{if }x\in U \\
0 & \mbox{otherwise.}
\end{array}\right.
\end{array}
\end{equation}

\noindent It is not hard to see that it is natural in $X$.

\auxproof{
This $\eta$ is a natural map: for $f\colon X\rightarrow Y$ in \Sets we
have for $x\in X$ and $V\in\Pow(Y)$,
$$\begin{array}{rcl}
\big(T(f) \after \eta_{X}\big)(x)(V) = 1
& \Longleftrightarrow &
T(f)(\eta_{X}(x))(V) = 1 \\
& \Longleftrightarrow &
(\eta_{X}(x) \after f^{-1})(V) = 1 \\
& \Longleftrightarrow &
\eta_{X}(x)(f^{-1}(V)) = 1 \\
& \Longleftrightarrow &
x\in f^{-1}(V) \\
& \Longleftrightarrow &
f(x) \in V \\
& \Longleftrightarrow &
\eta_{Y}(f(x))(V) = 1 \\
& \Longleftrightarrow &
\big(\eta_{Y} \after f)(x)(V) = 1.
\end{array}$$
}

Next we review that $\UF(X)$ is a topological space. The basic opens
are the sets $D(U) = \setin{\varphi}{\UF(X)}{\varphi(U) \neq 0}$, for
$U\in\Pow(X)$. Then:
\begin{equation}
\label{UltraFilterOpensEqn}
\begin{array}{rcl}
D(0)
& = &
\emptyset \\
D(1)
& = &
\UF(X) \\
D(U\cap V)
& = &
D(U) \cap D(V) \\
D(U \cup V)
& = &
D(U) \cup D(V) \\
D(\neg U)
& = &
\neg D(U).
\end{array}
\end{equation}

\noindent This makes $D\colon \Pow(X)\rightarrow \Pow(\UF(X))$ a
map of Boolean algebras. Notice that the last point implies that
subsets of the form $D(U)$ are also closed and thus clopen.

Further, 
$$\neg D(U)
=
\set{\varphi}{\varphi(U) = 0}
=
\set{\varphi}{U\in\ker(\varphi)},$$

\noindent where $\ker(\varphi)$ is the kernel of $\varphi$. Then, for
$A\subseteq \Pow(X)$,
$$\begin{array}{rcl}
\neg \bigcup_{U\in A}D(U)
& = &
\bigcap_{U\in A}\neg D(U) \\
& = &
\bigcap_{U\in A}\set{\varphi}{U\in\ker(\varphi)} \\
& = &
\set{\varphi}{\allin{U}{A}{U\in\ker(\varphi)}} \\
& = &
\set{\varphi}{A\subseteq \ker(\varphi)} \\
& = &
\set{\varphi}{\ideal A\subseteq \ker(\varphi)} \\
\end{array}$$

\noindent This allows us to show that $D(V)$ are compact:
$$\begin{array}{rcl}
D(V) \subseteq \bigcup_{U\in A}D(U)
& \Longrightarrow &
\neg\bigcup_{U\in A}D(U) \subseteq \neg D(V) \\
& \Longrightarrow &
\all{\varphi}{\ideal A\subseteq \ker(\varphi) \Rightarrow V\in\ker(\varphi)} \\
& \smash{\stackrel{(*)}{\Longrightarrow}} &
V\in\ideal A \\
& \Longrightarrow &
V \subseteq U_{1} \cup \cdots \cup U_{n}, \mbox{ for some }
   U_{1}, \ldots, U_{n} \in A \\
& \Longrightarrow &
D(V) \subseteq D(U_{1} \cup \cdots \cup U_{n}) \subseteq 
   D(U_{1}) \cup \cdots \cup D(U_{n}).
\end{array}$$

\noindent The marked implication
$\smash{\stackrel{(*)}{\Longrightarrow}}$ is a crucial step. Suppose
$V\not\in \ideal A$. Then the filter $\upset V$ and the ideal $\ideal
A$ are disjoint: $\upset V \cap \ideal A = \emptyset$.  By the
Ultrafilter Lemma~\ref{UltraFilterLem} there is a $\varphi\in\UF(X)$
with $\upset V \subseteq \varphi^{-1}(1)$ and $\ideal A \subseteq
\varphi^{-1}(0) = \ker(\varphi)$ and thus by assumption $V\in
\ker(\varphi)$. However, this contradicts $V\in\upset V\subseteq
\varphi^{-1}(1)$. Hence $V\in\ideal A$.

In particular, $D(1) = \UF(X)$ is compact, so that $\UF(X)$ is a
compact topological space. It is also Hausdorff: if we have
$\varphi,\psi \in \UF(X)$ with $\varphi\neq\psi$, then there must be a
point $U\in\Pow(X)$ with $\varphi(U)\neq\psi(U)$. Without loss of
generality, assume $\varphi(U) = 1$ and $\psi(U) = 0$. Then
$\varphi\in D(U)$, and $\psi\in D(\neg U)$, since $\psi(\neg U) = \neg
\psi(U) = \neg 0 = 1$, using that $\psi$ is a morphism of Boolean
algebras. Hence there are disjoint opens, namely $D(U)$ and $D(\neg
U)$ that separate $\varphi$ and $\psi$.

It is not hard to see that the maps $\UF(f) = - \after f^{-1}$ are
continuous, since $\UF(f)^{-1}(D(V)) = D(f^{-1}(V))$. We may thus
conclude that we have a functor $\UF\colon\Sets \rightarrow \CH$, from
sets to compact Hausdorff spaces and continuous maps between them.

\auxproof{
For a basic open $D(V) = \setin{\psi}{\UF(Y)}{\psi(V)\neq 0}$ we need
to show that we get an open subset:
$$\begin{array}{rcl}
\UF(f)^{-1}(D(V))
& = &
\setin{\varphi}{\UF(X)}{\UF(f)(\varphi)\in D(V)} \\
& = &
\setin{\varphi}{\UF(X)}{\varphi(f^{-1}(V)) \neq 0} \\
& = &
D(f^{-1}(V)).
\end{array}$$
}

\subsection{Compact Hausdorff spaces}\label{CHSubsec}

Assume now that $X\in\CH$ is an arbitrary compact Hausdorff
space. For an element $x\in X$ we define its \emph{neighborhood filter}
$\mathcal{N}(x)\subseteq \Pow(X)$ as:
$$\begin{array}{rcl}
\mathcal{N}(x)
& = &
\setin{U}{\Pow(X)}{\exin{U'}{\open(X)}{x\in U' \subseteq U}}.
\end{array}$$

\noindent Clearly, this is a proper filter:
\begin{itemize}
\item $\mathcal{N}(x)$ is an upset: if $V\supseteq U\in
  \mathcal{N}(x)$, say via $x\in U'\subseteq U$, then $x\in
  U'\subseteq V$ so that $V\in \mathcal{N}(x)$.

\item $\mathcal{N}(x)$ contains the top element $X\in\Pow(X)$, via
  $x\in X\subseteq X$.

\item $\mathcal{N}(x)$ is closed under meets: if $U,V\in
  \mathcal{N}(x)$, say via $x\in U'\subseteq U$ and $x\in V'\subseteq
  V$, then $U'\cap V'$ is open and satisfies $x\in U'\cap V' \subseteq
  U\cap V$. Hence $U\cap V\in \mathcal{N}(x)$.

\item the bottom element $\emptyset$ is not in $\mathcal{N}(x)$, since
  there is the only open $U\subseteq \emptyset$, namely $U=\emptyset$,
  and it does not contain $x$.
\end{itemize}

\noindent For $\varphi\in\UF(X)$ and $x\in X$ one defines the
\emph{convergence} relation:
$$\begin{array}{rcl}
\varphi\downset x
& \Longleftrightarrow &
\mathcal{N}(x) \subseteq \varphi^{-1}(1) \\
& \Longleftrightarrow &
\allin{U}{\mathcal{N}(x)}{\varphi(U) = 1} \\
& \Longleftrightarrow &
\allin{U}{\Pow(X)}{\allin{U'}{\open(X)}{x\in U'\subseteq U
   \Rightarrow \varphi(U) = 1}}.
\end{array}$$

\noindent Then we have the following basic result, see~\cite[III,
  2.2]{Johnstone82}, expressing basic topological properties in terms
of convergence.

\begin{proposition}
\label{ConvergenceProp}
For a topological space $X$:
\begin{enumerate}
\item $X$ is Hausdorff if and only if
$\allin{\varphi}{\UF(X)}{\allin{x,y}{X}{\varphi\downset x \conjun
   \varphi\downset y \Rightarrow x=y}}$.

\item $X$ is compact if and only if
$\allin{\varphi}{\UF(X)}{\exin{x}{X}{\varphi\downset x}}$.
\end{enumerate}
\end{proposition}

\begin{myproof}
{\em(1)}~First assume $X$ is Hausdorff and $\varphi\downset x$ and
$\varphi\downset y$, so that $\mathcal{N}(x) \subseteq
\varphi^{-1}(1)$ and $\mathcal{N}(y) \subseteq \varphi^{-1}(1)$. If
$x\neq y$ then there are $U,V\in\open(X)$ with $U\cap V =
\emptyset$ and $x\in U$, $y\in V$. But then $U\in\mathcal{N}(x)
\subseteq \varphi^{-1}(1)$ and $V\in\mathcal{N}(y) \subseteq
\varphi^{-1}(1)$, so that $\varphi(U) = 1 = \varphi(V)$. This yields a
contradiction: $1 = 1\conjun 1 = \varphi(U) \conjun \varphi(V) =
\varphi(U\cap V) = \varphi(\emptyset) = 0$. Hence $x=y$.

In the reverse direction, assume
$\allin{\varphi}{\UF(X)}{\allin{x,y}{X}{\varphi\downset x \conjun
    \varphi\downset y \Rightarrow x=y}}$. In order to prove that $X$
is Hausdorff, assume $x\neq y$ for $x,y\in X$, and consider the set
$$\begin{array}{rcl}
\mathcal{N}(x)\& \mathcal{N}(y)
& = &
\set{U\cap V}{U\in\mathcal{N}(x), V\in\mathcal{N}(y)}.
\end{array}$$

\noindent It contains both $\mathcal{N}(x)$ and $\mathcal{N}(y)$ as
subsets. Furthermore, it is a filter:
\begin{itemize}
\item $\mathcal{N}(x)\& \mathcal{N}(y)$ is an upset: if $W\supseteq U\cap V$
for $U\in\mathcal{N}(x)$ and $V\in\mathcal{N}(y)$, then 
$W\cup U\in\mathcal{N}(x)$ and $W\cup V\in\mathcal{N}(y)$ because they
are upsets. Hence:
$$W 
= 
W\cup (U\cap V)
=
(W\cup U) \cap (W\cup V) \in \mathcal{N}(x)\& \mathcal{N}(y).$$

\item $\mathcal{N}(x)\& \mathcal{N}(y)$ is obviously closed under
  finite intersections.
\end{itemize}

\noindent Now we distinguish two options:
\begin{itemize}
\item $\mathcal{N}(x)\& \mathcal{N}(y)$ is not a proper filter,
  \textit{i.e.}~it does contain $\emptyset$: $U\cap V = \emptyset$,
  say for $U\in\mathcal{N}(x)$ and $V\in\mathcal{N}(y)$, say with
  $x\in U'\subseteq U$ and $y\in V'\subseteq V$, for
  $U',V'\in\open(X)$. Then $U'\cap V' \subseteq U\cap V = \emptyset$,
  so $x,y$ are separated by disjoint opens.

\item $\mathcal{N}(x)\& \mathcal{N}(y)$ is a proper filter.  By the
  Lemma~\ref{UltraFilterLem} there is now a $\varphi\in\UF(X)$ with
  $\mathcal{N}(x)\& \mathcal{N}(y) \subseteq \varphi^{-1}(1)$.  But
  then $\mathcal{N}(x), \mathcal{N}(y) \subseteq \varphi^{-1}(1)$.
  This means $\varphi\downset x$ and $\varphi\downset y$, and thus
  $x=y$ by assumption. This is impossible.
\end{itemize}

{\em(2)}~Now assume that $X$ is compact. For $\varphi\in\UF(X)$
consider the set $U =
\bigcap\setin{V}{\closed(X)}{\varphi(V)=1}$. This set $U$ is
non-empty, because if $U=\emptyset$, then $X = \neg U =
\bigcup\set{\neg V}{V\in\closed(X), \varphi(V) = 1}$, and so $X
\subseteq \neg V_{1} \cup \cdots \cup \neg V_{n}$ for certain closed
$V_{i}$ with $\varphi(V_{i}) = 1$; but this yields a contradiction: $1
= \varphi(X) \leq \varphi(\neg V_{1} \cup \cdots \cup \neg V_{n}) =
\neg \varphi(V_{1}) \disjun \cdots \disjun \neg \varphi(V_{n}) = 0
\disjun \cdots \disjun 0 = 0$.

Hence we may assume an element $x\in U$. We aim to show
$\varphi\downset x$, \textit{i.e.}~$\mathcal{N}(x) \subseteq
\varphi^{-1}(1)$. So assume $V\in\mathcal{N}(x)$, say via $x\in
V'\subseteq V$, where $V'$ is open. We need to show $\varphi(V) =
1$. Assume otherwise, \textit{i.e.}~$\varphi(V) = 0$; then also
$\varphi(V') \leq \varphi(V) = 0$ and thus $\varphi(\neg V') =
1$. Since $\neg V'$ is closed we get $U\subseteq \neg V'$, and thus
$x\in U\subseteq \neg V'$. This is impossible since $x\in V'$. Hence
$\varphi(V) = 1$.

In the reverse direction, in order to prove that $X$ is compact assume
$X \subseteq \bigcup_{i}U_{i}$ for $U_{i}\in\open(X)$. Consider the
generated filter $F = \filter\set{\neg U_{i}}{i\in I} \subseteq
\Pow(X)$. Distinguish:
\begin{itemize}
\item $F$ is not proper. This means $\emptyset \in U$, and thus that
  $\emptyset \supseteq \neg U_{i_1} \cap \cdots \cap \neg U_{i_n}$,
  for certain $i_{j}\in I$. Hence $X \subseteq \neg (\neg U_{i_1} \cap
  \cdots \cap \neg U_{i_n}) = U_{i_1} \cup \cdots \cup U_{i_n}$. This
  means we are done.

\item $F$ is proper. Then we can apply the Lemma~\ref{UltraFilterLem}
  and obtain a $\varphi\in\UF(X)$ with $F\subseteq
  \varphi^{-1}(1)$. Since $\neg U_{i}\in F$ this means $\varphi(U_{i})
  = \neg \varphi(\neg U_{i}) = \neg 1 = 0$. By assumption there is an
  $x\in X = \bigcup_{i\in I}U_{i}$ with $\varphi\downset x$, and thus
  $\mathcal{N}(x) \subseteq \varphi^{-1}(1)$. Assume $x\in
  U_{i}\in\open(X)$; then $U_{i}\in\mathcal{N}(x) \subseteq
  \varphi^{-1}(1)$, so that $\varphi(U_{i}) = 1$. This second option
  is thus impossible. \QED
\end{itemize}
\end{myproof}

With this result one can define, for a compact Hausdorff space $X$, a
map $\varepsilon_{X}\colon \UF(X)\rightarrow X$ by:
\begin{equation}
\label{UFcounitEqn}
\begin{array}{rcccl}
\varepsilon_{X}(\varphi) = x
& \Longleftrightarrow &
\varphi\downset x
& \Longleftrightarrow &
\mathcal{N}(x) \subseteq \varphi^{-1}(1).
\end{array}
\end{equation}

\noindent In particular, $\mathcal{N}(\varepsilon_{X}(\varphi))
\subseteq \varphi^{-1}(1)$.

\begin{lemma}
\label{UFcounitLem}
The functions $\varepsilon_{X}\colon \UF(X)\rightarrow X$ defined
in~(\ref{UFcounitEqn}) satisfy the following properties.
\begin{enumerate}
\item Naturality in $X$.

\item $\{\varepsilon_{X}(\varphi)\} =
  \bigcap\set{\overline{U}}{U\in\Pow(X), \varphi(U) = 1}$.

\item The function $\varepsilon_X$ is continuous, with for $U\in\open(X)$,
$$\begin{array}{rcl}
\varepsilon_{X}^{-1}(U)
& = &
\setin{\varphi}{\UF(X)}{\exin{V}{\Pow(X)}{\overline{V} \subseteq U
   \conjun \varphi(V)=1}} \\
& = &
\bigcup_{\overline{V}\subseteq U}D(V).
\end{array}$$
\end{enumerate}
\end{lemma}

\begin{myproof}
{\em(1)}~For a continuous function $f\colon X\rightarrow Y$ we show
$f \after \varepsilon_{X} = \varepsilon_{Y} \after \UF(f)$. That is, for
$\varphi\in\UF(X)$ and $y\in Y$,
$$\begin{array}{rcl}
f(\varepsilon_{X}(\varphi)) = y
& \Longleftrightarrow &
\varepsilon_{Y}(\UF(f)(\varphi)) = y \\
& \Longleftrightarrow &
\varepsilon_{Y}(\varphi \after f^{-1}) = y \\
& \Longleftrightarrow &
\mathcal{N}(y) \subseteq \setin{V}{\Pow(X)}{\varphi(f^{-1}(V)) = 1}.
\end{array}$$

\noindent For the implication $(\Rightarrow)$ assume
$f(\varepsilon_{X}(\varphi)) = y$ and $V\in\mathcal{N}(y)$, say via
$y\in V'\subseteq V$ with $V'\in\open(Y)$. Then
$\varepsilon_{X}(\varphi) \in f^{-1}(V')$, so that $f^{-1}(V') \in
\mathcal{N}(\varepsilon_{X}(\varphi)) \subseteq
\varphi^{-1}(1)$. Hence $\varphi(f^{-1}(V')) = 1$, then thus, since
$V'\subseteq V$, also $\varphi(f^{-1}(V)) = 1$, as required.

For the reverse implication $(\Leftarrow)$ assume
$f(\varepsilon_{X}(\varphi)) \neq y$. Then there are disjoint
$V,W\in\open(Y)$ with $f(\varepsilon_{X}(\varphi)) \in V$ and $y\in
W$. One has $\varepsilon_{X}(\varphi) \in f^{-1}(V)$, so
that $f^{-1}(V) \in \mathcal{N}(\varepsilon_{X}(\varphi)) \subseteq
\varphi^{-1}(1)$, and thus $\varphi(f^{-1}(V)) = 1$. Similarly, since
$W\in\mathcal{N}(y)$ we get $\varphi(f^{-1}(W)) = 1$ by
assumption. But now we have a contradiction: $1 = 1\conjun 1 =
\varphi(f^{-1}(V)) \conjun \varphi(f^{-1}(W)) = \varphi(f^{-1}(V) \cap
f^{-1}(W)) = \varphi(f^{-1}(V\cap W)) = \varphi(f^{-1}(\emptyset)) =
\varphi(\emptyset) = 0$.

{\em(2)}~For the inclusion $(\subseteq)$ assume $x = \varepsilon_{X}(\varphi)$
and $\varphi(U) = 1$, for $U\in\Pow(X)$. We need to show
$x\in\overline{U}$. If not then $x\in\neg\overline{U}\in\open(X)$, so
$\neg\overline{U}\in\mathcal{N}(x) \subseteq \varphi^{-1}(1)$ and thus
$\varphi(\neg\overline{U}) = 1$, and so $\varphi(U) \leq
\varphi(\overline{U}) = 0$. This is impossible.

For $(\supseteq)$ assume $x\in\overline{U}$ for all $U\in\Pow(X)$ with
$\varphi(U) = 1$. We need to show $\varepsilon_{X}(\varphi) = x$,
\textit{i.e.}~$\mathcal{N}(x) \subseteq \varphi^{-1}(1)$. Hence assume
$V\in\mathcal{N}(x)$, via $x\in V'\subseteq V$, with
$V'\in\open(X)$. We need to show $\varphi(V) = 1$. If not,
$\varphi(\neg V) = 1$, so $x\in\overline{\neg V}$. From $V'\subseteq
V$ we get $\neg V \subseteq \neg V'\in\closed(X)$, and so
$x\in\overline{\neg V}\subseteq \neg V'$.  This contradicts $x\in V'$.

{\em(3)}~First we recall that $X$ is regular, because it is a compact
Hausdorff space; this means that points and closed subsets can be
separated, or equivalently: $U =
\bigcup\setin{V}{\open(X)}{\overline{V}\subseteq U}$, for each
$U\in\open(X)$.

\auxproof{ 
First we recall that a compact Hausdorff space $X$ is regular. If
$x\not\in W\in\closed(X)$. Then $W$ is compact, since any
cover of $W$ can be extended to a cover of $X$ by adding $\neg W$.
Because $X$ is Hausdorff, for each $y\in W$ there are disjoint $U_{y},
V_{y}\in\open(X)$ with $x\in U_{y}$ and $y\in V_{y}$. Hence
$W\subseteq \bigcup_{y\in W} V_{y}$, and thus by compactness
$W\subseteq V_{y_{1}} \cup \cdots \cup V_{y_{n}}\in\open(X)$ for
certain $y_{1},\ldots, y_{n}\in W$. Then also $x\in U_{y_1} \cap
\cdots \cap U_{y_n}\in\open(X)$, where clearly $V_{y_{1}} \cup \cdots
\cup V_{y_{n}}$ and $U_{y_1} \cap \cdots \cap U_{y_n}$ are disjoint.

Next we show that regularity is equivalent to: $U =
\bigcup\setin{V}{\open(X)}{\overline{V}\subseteq U}$, for
$U\in\open(X)$.

In one direction, assume that $X$ is regular. The inclusion
$(\supseteq)$ is obvious, and so we concentrate on
$(\subseteq)$. Suppose $x\in U\in\open(X)$. Then $x\not\in\neg
U\in\closed(X)$. By regularity there are disjoint $V,W\in\open(X)$
with $x\in V$ and $\neg U\subseteq W$.  Thus $V \subseteq \neg W$, and
thus $\overline{\neg V} \subseteq \neg W \subseteq U$.

Conversely, assume we have $x\not\in W\in\closed(X)$. Then $x\in \neg
W = \bigcup\setin{V}{\open(X)}{\overline{V} \subseteq \neg W}$, and so
$x\in V$ for some $V\in\open(X)$ with $\overline{V}\subseteq \neg W$.
Thus $W\subseteq \neg\overline{V}\in\open(X)$.
}

First assume $x = \varepsilon_{X}(\varphi)\in U =
\bigcup\setin{V}{\open(X)}{\overline{V}\subseteq U}$, say $x\in
V\in\open(X)$, with $\overline{V}\subseteq U$. Then $V\in\mathcal{N}(x)
\subseteq \varphi^{-1}(1)$, so that $\varphi(V) = 1$.

In the reverse direction, assume $\varphi(V) = 1$ for
$\overline{V}\subseteq U$. Then $\{a(\varphi)\} =
\bigcap\set{\overline{W}}{W\in\Pow(X), \varphi(W) = 1} \subseteq
\overline{V} \subseteq U$.
\end{myproof}

\begin{theorem}
\label{CHAdjThm}
The forgetful functor $\CH\rightarrow\Sets$ has $\UF$ as left adjoint.
\end{theorem}

\begin{myproof}
We already have (natural) unit~(\ref{UFunitEqn}) and
counit~(\ref{UFcounitEqn}) maps, so we only need to check the
triangular identities. First, for $X\in\CH$ we need to see that
$\varepsilon_{X} \after \eta_{X} = \idmap[X]$, where the forgetful
functor is not written. Suppose $x\in X$ and
$\varepsilon_{X}(\eta_{X}(x)) \neq x$. Then there are disjoint open
$U,V$ with $\varepsilon_{X}(\eta_{X}(x)) \in U$ and $x\in V$. The
former yields $U\in\mathcal{N}\big(\varepsilon_{X}(\eta_{X}(x))\big)
\subseteq \eta_{X}(x)^{-1}(1) = \set{V}{\eta(x)(V) = 1} = \set{V}{x\in
  V}$. Hence $x\in U$, which contradicts disjointness of $U$ and $V$.

Next for $Y\in\Sets$ we need to show $\varepsilon_{\UF(X)} \after
\UF(\eta_{X}) = \idmap$. For $\varphi\in\UF(X)$ write $\psi =
\varepsilon_{\UF(X)}\big(\UF(\eta_{X})(\varphi)\big)$. Suppose
$\psi(U)=1$, for $U\in\Pow(X)$. Notice that:
$$\begin{array}{rcl}
\mathcal{N}(\psi)
& = &
\setin{A}{\Pow(\UF(X))}{\exin{A'}{\open(\UF(X))}{\psi\in A'\subseteq A}} \\
& = &
\setin{A}{\Pow(\UF(X))}{\exin{V}{\Pow(X)}{\psi\in D(V)\subseteq A}} \\
& = &
\setin{A}{\Pow(\UF(X))}{\exin{V}{\Pow(X)}{\psi(V) = 1 \conjun
   D(V)\subseteq A}}.
\end{array}$$

\noindent Thus $D(U)\in\mathcal{N}(\psi)$ and thus: 
$$\begin{array}{rcl}
D(U)\in \big(\UF(\eta_{X})(\varphi)\big)^{-1}(1)
& = &
\set{A}{\varphi(\eta_{X}^{-1}(A)) = 1} \\
& = &
\set{A}{\varphi(\set{x}{\eta(x)\in A})=1}.
\end{array}$$

\noindent Hence: 
$$1 = \varphi(\set{x}{\eta(x)\in D(U)}) =
\varphi(\set{x}{\eta(x)(U) = 1}) = \varphi(\set{x}{x\in U}) =
\varphi(U).$$

\noindent This yields $\psi=\varphi$ by Lemma~\ref{UltraFilterEqLem}. \QED
\end{myproof}

\auxproof{
Earlier text, starting with a monad description.

The multiplication $\mu\colon \UF^{2} \Rightarrow \UF$ is
defined on $\Phi\colon \Pow(\UF(X))\rightarrow 2$ and $U\in \Pow(X)$ as:
$$\begin{array}{rcl}
\mu_{X}(\Phi)(U)
& = &
\Phi\big(\set{\varphi\colon\Pow(X)\rightarrow 2}{\varphi(U) = 1}\big).
\end{array}$$

\noindent It takes some effort but it is essentially straightforward
to show that $\UF$ with this $\eta$ and $\mu$ is indeed a monad.

We first check naturality of $\mu$. For $f\colon X\rightarrow Y$ we have,
for $\Phi\in\UF^{2}(X) = \BA(\Pow(\UF(X)), 2)$ and $V\in\Pow(Y)$,
$$\begin{array}{rcl}
\lefteqn{\big(\mu_{Y} \after \UF^{2}(f)\big)(\Phi)(V) = 1} \\
& \Longleftrightarrow &
\mu_{Y}\big(\UF(\UF(f))(\Phi)\big)(V) = 1 \\
& \Longleftrightarrow &
\mu_{Y}\big(\Phi \after \UF(f)^{-1}\big)(V) = 1 \\
& \Longleftrightarrow &
\big(\Phi \after \UF(f)^{-1}\big)\big(\set{\psi\colon\Pow(Y)\rightarrow 2}{
   \psi(V) = 1}\big) = 1 \\
& \Longleftrightarrow &
\Phi\big(\set{\varphi\colon\Pow(X)\rightarrow 2}{\UF(f)(\varphi)(V) = 1}\big) 
  = 1 \\
& \Longleftrightarrow &
\Phi\big(\set{\varphi\colon\Pow(X)\rightarrow 2}{\varphi(f^{-1}(V)) = 1}\big) 
  = 1 \\
& \Longleftrightarrow &
\mu(\Phi)(f^{-1}(V)) = 1 \\
& \Longleftrightarrow &
\big(\mu_{X}(\Phi) \after f^{-1}\big)(V) = 1 \\
& \Longleftrightarrow &
\UF(f)\big(\mu_{X}(\Phi)\big)(V) = 1 \\
& \Longleftrightarrow &
\big(\UF(f) \after \mu_{X}\big)(\Phi)(V) = 1.
\end{array}$$

\noindent We continue with the monad equations. First, for
$\psi\in\UF(X)$ and $U\in\Pow(X)$,
$$\begin{array}{rcl}
\big(\mu_{X} \after \eta_{\UF(X)}\big)(\psi)(U) = 1
& \Longleftrightarrow &
\mu_{X}(\eta_{\UF(X)}(\psi))(U) = 1 \\
& \Longleftrightarrow &
\eta_{\UF(X)}(\psi)\big(\set{\varphi\colon\Pow(X)\rightarrow 2}
   {\varphi(U) = 1}\big) = 1 \\
& \Longleftrightarrow &
\psi(U) = 1 \\
& \Longleftrightarrow &
\idmap[\UF(X)](\psi)(U) = 1 \\
\big(\mu_{X} \after \UF(\eta_{X})\big)(\psi)(U) = 1
& \Longleftrightarrow &
\mu_{X}(\UF(\eta_{X})(\psi))(U) = 1 \\
& \Longleftrightarrow &
\mu_{X}(\psi \after \eta_{X}^{-1})(U) = 1 \\
& = &
\big(\psi \after \eta_{X}^{-1}\big)\big(\set{\varphi\colon\Pow(X)\rightarrow 2}
   {\varphi(U) = 1}\big) = 1 \\
& \Longleftrightarrow &
\psi\big(\set{x}{\eta(x)(U) = 1}\big) = 1 \\
& \Longleftrightarrow &
\psi\big(\set{x}{x\in U}\big) = 1 \\
& \Longleftrightarrow &
\psi(U) = 1 \\
& \Longleftrightarrow &
\idmap[\UF(X)](\psi)(U) = 1.
\end{array}$$

\noindent Next for $F\in\UF^{3}(X) = \BA(\UF^{2}(X), 2)$ and $U\in\Pow(X)$
we have:
$$\begin{array}{rcl}
\lefteqn{\big(\mu_{X} \after \UF(\mu_{X})\big)(F)(U)} \\
& = &
\mu_{X}\big(\UF(\mu_{X})(F)\big)(U) \\
& = &
\mu_{X}\big(F \after \mu_{X}^{-1}\big)(U) \\
& = &
\big(F \after \mu_{X}^{-1}\big)\big(\set{\varphi\colon\Pow(X)\rightarrow 2}
   {\varphi(U) = 1}\big) \\
& = &
F\big(\set{\Phi\colon\Pow(\UF(X))\rightarrow 2}{\mu_{X}(\Phi)(U) = 1}\big) \\
& = &
F\big(\set{\Phi\colon\Pow(\UF(X))\rightarrow 2}
   {\Phi(\set{\varphi\colon \Pow(X)\rightarrow 2}{\varphi(U)=1}) = 1}\big) \\
& = &
\mu_{\UF(X)}(F)\big(\set{\varphi\colon\Pow(X)\rightarrow 2}
   {\varphi(U) = 1}\big) \\
& = &
\mu_{X}\big(\mu_{\UF(X)}(F)\big)(U) \\
& = &
\big(\mu_{X} \after \mu_{\UF(X)}\big)(F)(U).
\end{array}$$

This $\UF$ is not a commutative monad: recall first that the strength
map $\st\colon X\times\UF(Y) \rightarrow \UF(X\times Y)$ is defined
for $x\in X$, $\psi\in\UF(Y)$ and $W\in\Pow(X\times Y)$ as:
$$\begin{array}{rcl}
\st(x, \psi)(W)
& = &
\UF(\lamin{y}{Y}{\tuple{x,y}})(\psi)(V) \\
& = &
\big(\psi \after (\lamin{y}{Y}{\tuple{x,y}})^{-1}\big)(W) \\
& = &
\psi(\set{y}{\tuple{x,y}\in W}).
\end{array}$$

\noindent Hence the twisted strength $\st'\colon \UF(X)\times Y
\rightarrow \UF(X\times Y)$ is given by $\st'(\varphi,y)(W) =
\varphi(\set{x}{\tuple{x,y}\in W})$. Thus:
$$\begin{array}{rcl}
\lefteqn{\big(\mu \after \UF(\st) \after \st'\big)(\varphi, \psi)(W)} \\
& = &
\mu\big(\UF(\st)(\st'(\varphi,\psi))\big)(W) \\
& = &
\mu\big(\st'(\varphi,\psi) \after \st^{-1}\big)(W) \\
& = &
\big(\st'(\varphi,\psi) \after \st^{-1}\big)\big(
   \set{\chi\colon\Pow(X\times Y)\rightarrow 2}{\chi(W) = 1}\big) \\
& = &
\st'(\varphi,\psi)\big(\set{(x,\psi')}{\st(x,\psi')(W) = 1}\big) \\
& = &
\st'(\varphi,\psi)\big(\set{(x,\psi')}
   {\psi'(\set{y}{\tuple{x,y}\in W}) = 1}\big) \\
& = &
\varphi(\set{x}{\psi(\set{y}{\tuple{x,y}\in W}) = 1}).
\end{array}$$
}

We shall also write $\UF\colon \Sets\rightarrow\Sets$ for the monad
induced by this adjunction, with multiplication $\mu\colon \UF^{2}(X)
\rightarrow \UF(X)$ characterized by:
\begin{equation}
\label{UltraFilterMultiplicationEqn}
\begin{array}{rcl}
\mu(\Phi)(U) = 1
& \Longleftrightarrow &
\Phi(D(U))=1
\end{array}
\end{equation}

\auxproof{
We have $\mu_{X}(\Phi) = \varepsilon_{\UF(X)}(\Phi)$ and thus
$$\qquad\begin{array}{rcl}
\lefteqn{\llap{$\Phi^{-1}(1) \supseteq\;$} \mathcal{N}(\mu_{X}(\Phi))} \\
& = &
\setin{A}{\Pow(\UF(X))}{\exin{A'}{\open(\UF(X))}
   {\mu(\Phi) \in A' \subseteq A}} \\
& = &
\setin{A}{\Pow(\UF(X))}{\exin{U}{\Pow(X)}
   {\mu(\Phi) \in D(U) \subseteq A}} \\
& = &
\setin{A}{\Pow(\UF(X))}{\exin{U}{\Pow(X)}
   {\mu(\Phi)(U) = 1 \conjun D(U) \subseteq A}}.
\end{array}$$

\noindent Hence if $\mu(\Phi)(U) = 1$, then $D(U)\in \Phi^{-1}(1)$ and
thus $\Phi(D(U)) = 1$. This suffices by Lemma~\ref{UltraFilterEqLem},
since $D\colon \Pow(X)\rightarrow \Pow(\UF(X))$ is a map of
Boolean algebras, see~(\ref{UltraFilterOpensEqn}).
}

\noindent The main result is that its category of
algebras is the category of compact Hausdorff spaces. As a first step
we consider the following partial result.

\begin{lemma}
\label{CHAlgLem}
The comparison functor $\CH\rightarrow \Alg(\UF)$ is full and
faithful: a function $f\colon X\rightarrow Y$ between compact
Hausdorff spaces $X,Y$ is continuous if and only if it is a morphism
of algebras $\varepsilon_{X} \rightarrow \varepsilon_{Y}$.
\end{lemma}

\begin{myproof}
Suppose $f\colon X\rightarrow Y$ is continuous. Then for
$\varphi\in\UF(X)$ and $y\in Y$,
$$\begin{array}{rcl}
\lefteqn{\big(\varepsilon_{Y} \after \UF(f)\big)(\varphi) = y} \\
& \Longleftrightarrow &
\varepsilon_{Y}(\UF(f)(\varphi)) = \varepsilon_{Y}(\varphi \after f^{-1}) = y \\
& \Longleftrightarrow &
\mathcal{N}(y) \subseteq (\varphi \after f^{-1})^{-1}(1)
   = \setin{V}{\Pow(Y)}{\varphi(f^{-1}(V)) = 1} \\
& \smash{\stackrel{(*)}{\Longleftrightarrow}} &
f(\varepsilon_{X}(\varphi)) = \big(f \after \varepsilon_{X}\big)(\varphi) = y.
\end{array}$$

\noindent For the marked equivalence
$\smash{\stackrel{(*)}{\Longleftrightarrow}}$ we first do the
direction $(\Rightarrow)$. Assume $\mathcal{N}(y) \subseteq
\setin{V}{\Pow(Y)}{\varphi(f^{-1}(V)) = 1}$, and $\varepsilon_{X}(\varphi) = x$
\textit{i.e.}~$\mathcal{N}(x) \subseteq \varphi^{-1}(1)$; we need to
prove $f(x) = y$. Suppose $f(x) \neq y$, so that we can find disjoint
open $V,W\in\open(X)$ with $f(x)\in V$ and $y\in W$. Then $x\in
f^{-1}(V)$, which is open because $f$ is continuous, so that
$f^{-1}(V)\in\mathcal{N}(x)$ and thus $\varphi(f^{-1}(V)) = 1$. At the
same time we have $W\in\mathcal{N}(y)$ and thus $\varphi(f^{-1}(W)) =
1$. This yields a contradiction: $1 = 1\conjun 1 = \varphi(f^{-1}(V))
\conjun \varphi(f^{-1}(W)) = \varphi(f^{-1}(V) \cap f^{-1}(W)) =
\varphi(f^{-1}(V\cap W)) = \varphi(f^{-1}(\emptyset)) =
\varphi(\emptyset) = 0$.

For the direction $(\Leftarrow)$ write $x=\varepsilon_{X}(\varphi)$ so that
$f(x)=y$ by assumption. We need to prove $\mathcal{N}(y) \subseteq
\setin{V}{\Pow(Y)}{\varphi(f^{-1}(V)) = 1}$. So assume
$V\in\mathcal{N}(y)$, say $y\in V'\subseteq V$ where $V'\in\open(Y)$.
Since $f$ is continuous we get $x\in f^{-1}(V') \in \open(X)$, so that
$f^{-1}(V')\in\mathcal{N}(x)$ and thus $\varphi(f^{-1}(V')) = 1$.
Since $V'\subseteq V$ we get $\varphi(f^{-1}(V)) = 1$, as required.

Next assume $f$ is a map of algebras $\varepsilon_{X} \rightarrow
\varepsilon_{Y}$, \textit{i.e.}~$\varepsilon_{Y} \after \UF(f) = f
\after \varepsilon_{X}$. We have to prove that $f$ is continuous.  So
assume $V\in\open(Y)$; we plan to prove $f^{-1}(V) =
\bigcup\setin{U}{\open(X)}{U\subseteq f^{-1}(V)}$. Take $x\in
f^{-1}(V)$ and form the filter $F = \filter\set{U -
  f^{-1}(V)}{U\in\mathcal{N}(x)}\subseteq \Pow(X)$. Distinguish:
\begin{itemize}
\item $F$ is not disjoint from the ideal $\downset f^{-1}(V)$. This
  means that there is a set $U\in \filter A \cap \downset
  f^{-1}(V)$. Hence $U\subseteq f^{-1}(V)$ and $U\supseteq (U_{1} -
  f^{-1}(V)) \cap \cdots \cap (U_{n}-f^{-1}(V)) = (U_{1}\cap \cdots
  \cap U_{n})-f^{-1}(V)$ for certain $U_{i}\in\mathcal{N}(x)$. Then
  $U_{1}\cap \cdots \cap U_{n} \subseteq f^{-1}(V)$, since if not,
  then there is a $z\in U_{1}\cap \cdots \cap U_{n}$ with $z\not\in
  f^{-1}(V)$, so that $z\in (U_{1}\cap \cdots \cap U_{n})-f^{-1}(V)
  \subseteq U \subseteq f^{-1}(V)$.

Find $U'_{i}\in\open(X)$ with $x\in U'_{i} \subseteq U_{i}$. Then
$x\in U'_{1} \cap \cdots \cap U'_{n} \subseteq U_{1} \cap \cdots \cap
U_{n} \subseteq U \subseteq f^{-1}(V)$, so we are done.

\item $F$ is disjoint from the ideal $\downset f^{-1}(V)$. The
  Ultrafilter Lemma yields a $\varphi\in\UF(X)$ with $\filter
  A\subseteq \varphi^{-1}(1)$ and $\downset f^{-1}(V) \subseteq
  \varphi^{-1}(0)$. In particular, $\varphi(f^{-1}(V)) = 0$. We have
  $\varepsilon_{X}(\varphi) = x$ since $\mathcal{N}(x) \subseteq
  \varphi^{-1}(1)$; indeed, if $U\in\mathcal{N}(x)$ we have
  $U\supseteq U - f^{-1}(V) \in F$, so $\varphi(U) = 1$. Since $f$ is
  a map of algebras we get: $f(x) = f(\varepsilon_{X}(\varphi)) =
  \varepsilon_{Y}(\UF(f)(\varphi)) = \varepsilon_{Y}(\varphi \after
  f^{-1})$. Thus $\mathcal{N}(f(x)) \subseteq (\varphi\after
  f^{-1})(1) = \setin{V}{\Pow(Y)}{\varphi(f^{-1}(V)) = 1}$. Since
  $V\in\open(Y)$ and $f(x)\in V$ we have $V\in\mathcal{N}(f(x))$ and
  thus $\varphi(f^{-1}(V)) = 1$. This is impossible. Hence this second
  option cannot occur. \QED
\end{itemize}
\end{myproof}

\begin{theorem}[Manes]
\label{CHAlgThm}
The category $\CH$ of compact Hausdorff spaces is algebraic over sets
via the ultrafilter monad $\UF$: the comparison functor
$\CH\rightarrow \Alg(\UF)$ is an equivalence of categories.
\end{theorem}

\begin{myproof}
Assume that $a\colon \UF(X)\rightarrow X$ in $\Sets$ is an
(Eilenberg-Moore) algebra of the monad $\UF$. It induces a compact
Hausdorff topology on $X$. To start one defines a closure operation
$\overline{(-)}$ on $\Pow(X)$ by:
$$\begin{array}{rcccl}
\overline{U}
& = &
\set{a(\varphi)}{\varphi\in\UF(X), \varphi(U)=1}
& = &
\set{a(\varphi)}{\varphi\in D(U)}.
\end{array}$$

\noindent Then:
\begin{itemize}
\item $U\subseteq\overline{U}$, since for $x\in U$ we have
$\eta(x)\in\UF(X)$ satisfying $\eta(x)(U) = 1$. From $a(\eta(x))=x$ we
then get $x\in\overline{U}$.

\item $\overline{(-)}$ is monotone: if $U\subseteq V$, then:
$$\setin{\varphi}{\UF(X)}{\varphi(U) = 1} 
\subseteq
\setin{\varphi}{\UF(X)}{\varphi(V) = 1}$$ 

\noindent and thus:
$$\overline{U} 
=
\set{a(\varphi)}{\varphi\in\UF(X), \varphi(U)=1} 
\subseteq
\set{a(\varphi)}{\varphi\in\UF(X), \varphi(V)=1} 
= 
\overline{V}.$$

\item $\overline{U} \cup \overline{V} = \overline{U\cup V}$, of which
  by monotony only the implication $(\supseteq)$ requires proof. So
  assume $x\in\overline{U\cup V}$, say via $x=a(\varphi)$ where $1 =
  \varphi(U\cup V) = \varphi(U) \disjun \varphi(V)$. Hence either
  $\varphi(U) = 1$ or $\varphi(V) = 1$, that is, either
  $x\in\overline{U}$ or $x\in\overline{V}$.
\end{itemize}

Idempotency $\overline{\overline{U}} = \overline{U}$ is more involved,
and requires the $\mu$-requirement for the algebra $a$, where
$\mu(\Phi) = $. Assume therefor $x\in\overline{\overline{U}}$ via $x =
a(\varphi)$ where $\varphi(\overline{U}) = 1$. The aim is to find a
$\psi\in\UF(X)$ with $x=a(\psi)$ and $\psi(U)=1$, so that
$x\in\overline{U}$.

First note that for $V\in\Pow(X)$ with $\varphi(V)=1$ we have $V\cap
\overline{U} \neq \emptyset$, because otherwise $1 = \varphi(V)
\conjun \varphi(\overline{U}) = \varphi(V\cap\overline{U}) =
\varphi(\emptyset) = 0$. Hence also $a^{-1}(V) \cap D(U) \neq
\emptyset$, because if $y\in V\cap\overline{U}$, then
$y\in\overline{U}$ means $y = a(\psi)$ for some $\psi$ with $\psi(U) =
1$; this means $\psi\in D(U)$. But $y = a(\psi) \in V$ means $\psi \in
a^{-1}(V)$. Hence $\psi \in a^{-1}(V) \cap D(U)$. One now forms the
following filter in $\UF^{2}(X)$.
$$\begin{array}{rcl}
A
& = &
\filter\set{a^{-1}(V)\cap D(U)}{V\in\varphi^{-1}(1)}.
\end{array}$$

\noindent This is a proper filter. If not, there are $V_{1},\ldots,V_{n}
\in \varphi^{-1}(1)$ with:
$$\begin{array}{rcl}
\lefteqn{\llap{$\emptyset \supseteq\;$} 
   \big(a^{-1}(V_{1}) \cap D(U)\big) \cap \cdots \cap
   \big(a^{-1}(V_{n}) \cap D(U)\big)} \\
& = &
\big(a^{-1}(V_{1}) \cap \cdots \cap a^{-1}(V_{n})\big) \cap D(U) \\
& = &
a^{-1}\big(V_{1} \cap \cdots \cap V_{n}\big) \cap D(U).
\end{array}$$

\noindent But the latter is non-empty, as we just noted, because
$V_{1} \cap \cdots \cap V_{n} \in\varphi^{-1}(1)$. Hence by the
Ultrafilter Lemma~\ref{UltraFilterLem} in $\Pow(\UF(X))$ there a
$\Phi\in\UF^{2}(X)$ with $A\subseteq \Phi^{-1}(1)$.

One now uses the following two facts.
\begin{enumerate}
\item $\UF(a)(\Phi) = \varphi$. Indeed, if for $V\in\Pow(X)$ we have
  $\varphi(V)=1$, then $V\in\varphi^{-1}(1)$ and thus $a^{-1}(V) \cap
  D(U) \in A$. Hence $\Phi(a^{-1}(V) \cap D(U)) = 1$, and thus also
  $\UF(a)(\Phi)(V) = \Phi(a^{-1}(V)) = 1$. Hence we are done by
  Lemma~\ref{UltraFilterEqLem}.

\item Write $\psi = \mu(\Phi)$. Then $\psi(U) = 1$ since $\Phi(D(U)) =
  1$ by construction, using~(\ref{UltraFilterMultiplicationEqn}),
  since $D(U)\in A$. Now $a(\psi) = a(\mu(\Phi)) = a(\UF(a)(\varphi))
  = a(\varphi)$. Hence we are done showing that
  $\overline{\overline{U}}\subseteq \overline{U}$.
\end{enumerate}

\noindent The set $X$ thus carries a topology in which $U\subseteq X$
is closed if
$$\allin{\varphi}{\UF(X)}{a(\varphi) \in U \Longrightarrow \varphi(U)=1}$$

\noindent That is, $U\subseteq X$ is open if $\neg U$ is closed, so:
$$\allin{\varphi}{\UF(X)}{\varphi(U)=1 \Longrightarrow 
   a(\varphi) \in U}$$

What remains is showing that the convergence map $\varepsilon_{X}
\colon \UF(X)\rightarrow X$ for this topology on $X$ coincides with
$a$, that is $\varepsilon_{X}(\varphi) = x$ iff $a(\varphi) = x$.
First the easy direction: if $x=a(\varphi)$, then for each
$V\in\Pow(X)$ with $\varphi(V) = 1$ we have $x\in\overline{V} =
\set{a(\psi)}{\psi\in\UF(X)\mbox{ with }\psi(V)=1}$, by definition of
closure. This means:
$$\begin{array}{rcccl}
x 
& \in &
\bigcap\set{\overline{V}}{V\in\varphi^{-1}(1)}
& = &
\{\varepsilon_{X}(\varphi)\},
\quad\mbox{by Lemma~\ref{UFcounitLem}}
\end{array}$$

\noindent In the other direction, assume $x=\varepsilon_{X}(\varphi)$,
so that $x\in\{\varepsilon_{X}(\varphi)\} =
\bigcap\set{\overline{V}}{V\in \varphi^{-1}(1)}$. Hence
$x\in\overline{V}$ for each $V\in\varphi^{-1}(1)$. Since by definition
$\overline{V} = \set{a(\psi)}{\psi\in D(V)}$, there is for each
$V\in\varphi^{-1}(1)$ a $\psi\in D(V)$ with $x=a(\psi)$.  Thus
$D(V)\cap a^{-1}(x)\neq\emptyset$, for each $V\in\varphi^{-1}(1)$.
Hence the filter $A = \filter\set{D(V)\cap
  a^{-1}(x)}{V\in\varphi^{-1}(1)} \subseteq \Pow(\UF(X))$ is proper:
if not, then $\emptyset \supseteq (D(V_{1})\cap a^{-1}(x)) \cap \cdots
\cap (D(V_{n})\cap a^{-1}(x)) = (D(V_{1}) \cap \cdots \cap
D(V_{n}))\cap a^{-1}(x) = D(V_{1} \cap \cdots \cap V_{n})\cap
a^{-1}(x)$ for certain $V_{1}, \ldots, V_{n}\in\varphi^{-1}(1)$; since
also $V_{1} \cap \cdots \cap V_{n} \in \varphi^{-1}(1)$ this is
impossible.

Thus there is a $\Phi\in\UF^{2}(X)$ with $A\subseteq \Phi^{-1}(1)$. 
It satisfies two crucial properties.
\begin{itemize}
\item $\UF(a)(\Phi) = \eta(x)$. Since $a^{-1}(x) \in A$ we have
$1 = \Phi(a^{-1}(x)) = \UF(a)(\Phi)(\{x\})$. Thus if $\eta(x)(U) = 1$,
then $x\in U$ and $\{x\}\subseteq U$, so that $\UF(a)(\Phi)(U) = 1$. 
This suffices by Lemma~\ref{UltraFilterEqLem}. 

\item $\mu(\Phi) = \varphi$. If $\varphi(U)=1$, then $\varphi\in D(U)
  \in A$, so that $\Phi(D(U)) = 1$, and thus $\mu(\Phi)(U) = 1$ 
by~(\ref{UltraFilterMultiplicationEqn}). Hence we are done again
by Lemma~\ref{UltraFilterEqLem}. 
\end{itemize}

\noindent Now we have $x = a(\eta(x)) = a(\UF(a)(\Phi)) = a(\mu(\Phi)) =
a(\varphi)$, as required.

In particular, for each $\varphi\in\UF(X)$ there is a unique $x$ with
$\varphi\downset x$, namely $x = \varepsilon_{X}(\varphi) =
a(\varphi)$. Hence $X$ is a compact Hausdorff space by
Proposition~\ref{ConvergenceProp}. \QED
\end{myproof}
}

\subsection{The continuation monad}\label{ContinuationSubsec}

The so-called continuation monad is useful in the context of
programming semantics, where it is employed for a particular style of
evaluation. The monad starts from a fixed set $C$ and takes the
``double dual'' of a set, where $C$ is used as dualizing object.
Hence we first form a functor $\mathcal{C} \colon \Sets \rightarrow
\Sets$ by:
$$\begin{array}{rclcrcl}
\mathcal{C}(X)
& = &
C^{(C^{X})}
& \quad\mbox{and}\quad &
\mathcal{C}\big(\smash{X \stackrel{f}{\rightarrow} Y}\big)
& = &
\lamin{h}{C^{(C^{X})}}{\lamin{g}{C^{Y}}{h(g \after f)}}.
\end{array}$$

\noindent This functor $\mathcal{C}$ forms a monad via:
$$\xymatrix@R-2pc@C-.5pc{
X\ar[r]^-{\eta} & C^{(C^{X})}
& &
C^{\Big(C^{\big(C^{(C^{X})}\big)}\Big)}\ar[r]^-{\mu} & C^{(C^{X})} \\
x\ar@{|->}[r] & \lamin{g}{C^{X}}{g(x)}
& &
H\ar@{|->}[r] & \lamin{g}{C^{X}}{H\big(\lamin{k}{C^{(C^{X})}}{k(g)}\big)}.
}$$

\auxproof{
Naturality of $\eta$ and $\mu$:
$$\begin{array}{rcl}
\big(\mathcal{C}(f) \after \eta_{X}\big)(x)
& = &
\mathcal{C}(f)(\eta(x)) \\
& = &
\lamin{g}{C^Y}{\eta(x)(g \after f)} \\
& = &
\lamin{g}{C^Y}{g(f(x))} \\
& = &
\eta(f(x)) \\
& = &
\big(\eta_{Y} \after f\big)(x) \\
\big(\mu_{Y} \after \mathcal{C}^{2}(f)\big)(H)(g) 
& = &
\mu_{Y}\big(\mathcal{C}^{2}(f)(H)\big)(g) \\
& = &
\mathcal{C}^{2}(f)(H)(\lam{k}{k(g)}) \\
& = &
H\big((\lam{k}{k(g)}) \after \mathcal{C}(f)\big) \\
& = &
H\big(\lam{h}{\mathcal{C}(f)(h)(g)}\big) \\
& = &
H\big(\lam{h}{h(g \after f)}\big) \\
& = &
\mu_{X}(H)(g \after f) \\
& = &
\mathcal{C}(f)(\mu_{X}(H))(g) \\
& = &
\big(\mathcal{C}(f) \after \mu_{X}\big)(H)(g).
\end{array}$$

Next we check the monad equations.
$$\begin{array}{rcl}
\big(\mu_{X} \after \mathcal{C}(\eta_{X})\big)(h)(g)
& = &
\mu\big(\mathcal{C}(\eta)(h)\big)(g) \\
& = &
\mathcal{C}(\eta)(h)(\lam{k}{k(g)}) \\
& = &
h\big((\lam{k}{k(g)}) \after \eta\big) \\
& = &
h\big(\lam{x}{\eta(x)(g)}\big) \\
& = &
h\big(\lam{x}{g(x)}\big) \\
& = &
h(g) \\
& = &
\idmap{}(h)(g) \\
\big(\mu_{X} \after \eta_{\mathcal{C}(X)}\big)(h)(g)
& = &
\mu\big(\eta(h)\big)(g) \\
& = &
\eta(h)(\lam{k}{k(g)}) \\
& = &
h(g).
\end{array}$$

\noindent For the $\mu$-equation assume $K\in\mathcal{C}^{3}(X)$.
Then:
$$\begin{array}{rcl}
\big(\mu_{X} \after \mathcal{C}(\mu_{X})\big)(K)(g)
& = &
\mu\big(\mathcal{C}(\mu)(K)\big)(g) \\
& = &
\mathcal{C}(\mu)(K)(\lam{k}{k(g)}) \\
& = &
K\big((\lam{k}{k(g)}) \after \mu\big) \\
& = &
K\big(\lam{H}{\mu(H)(g)}\big) \\
& = &
K\big(\lam{H}{H(\lam{k}{k(g)})}\big) \\
& = &
\mu(K)(\lam{k}{k(g)}) \\
& = &
\mu\big(\mu(K)\big)(g) \\
& = &
\big(\mu_{X} \after \mu_{\mathcal{C}(X)}\big)(K)(g).
\end{array}$$
}

The following folklore result will be useful in the present 
context.

\begin{lemma}
\label{ContAlgMndMapLem}
Let $T\colon \Sets \rightarrow \Sets$ be an arbitrary monad and 
$\mathcal{C}(X) = C^{(C^{X})}$ be the continuation monad on a set $C$.
Then there is a bijective correspondence between:
$$\begin{prooftree}
{\xymatrix{T(C)\ar[r]^-{\scriptstyle a} & C}
   \rlap{\qquad Eilenberg-Moore algebras}}
\Justifies
{\xymatrix{T \ar@{=>}[r]_-{\scriptstyle \sigma} & \mathcal{C}}
   \rlap{\qquad\hspace*{.5em} maps of monads.}}
\end{prooftree}\qquad\qquad$$
\end{lemma}

\begin{myproof}
First, given an algebra $a\colon T(C)\rightarrow C$ define
$\sigma_{X}\colon T(X) \rightarrow C^{(C^{X})}$ by:
$$\begin{array}{rcl}
\sigma_{X}(u)(g)
& = &
a\big(T(g)(u)\big).
\end{array}$$

\auxproof{
We check the details. First, naturality: for $f\colon X\rightarrow Y$,
$$\begin{array}{rcl}
\big(\sigma_{Y} \after T(f)\big)(u)(g)
& = &
\sigma_{Y}\big(T(f)(u)\big)(g) \\
& = &
a\big(T(g)(T(f)(u))\big) \\
& = &
a\big(T(g \after f)(u)\big) \\
& = &
\sigma_{X}(u)(g \after f) \\
& = &
\mathcal{C}(f)\big(\sigma_{X}(u)\big)(g) \\
& = &
\big(\mathcal{C}(f) \after \sigma_{X}\big)(u)(g).
\end{array}$$

\noindent Further, $\sigma$ commutes with the unit and multiplication:
$$\begin{array}{rcl}
\big(\sigma_{X} \after \eta_{X}\big)(x)(g)
& = &
\sigma_{X}\big(\eta_{X}(x)\big)(g) \\
& = &
a\big(T(g)(\eta_{X}(x))\big) \\
& = &
\big(a \after T(g) \after \eta_{X}\big)(x) \\
& = &
\big(a \after \eta \after g\big)(x) \\
& = &
g(x) \\
& = &
\eta(x)(g) \\
\big(\sigma_{X} \after \mu_{X}\big)(H)(g) 
& = &
\sigma_{X}\big(\mu_{X}(H)\big)(g) \\
& = &
a\big(T(g)(\mu_{X}(H))\big) \\
& = &
\big(a \after T(g) \after \mu_{X}\big)(H) \\
& = &
\big(a \after \mu_{C} \after T^{2}(g)\big)(H) \\
& = &
\big(a \after T(a) \after T^{2}(g)\big)(H) \\
& = &
a\big(T(a \after T(g))(H)\big) \\
& = &
a\big(T(\lam{f}{a\big(T(g)(f)\big)})(H)\big) \\
& = &
a\big(T(\lam{f}{\sigma_{X}(f)(g)})(H)\big) \\
& = &
\sigma_{TX}(H)\big(\lam{f}{\sigma_{X}(f)(g)}\big) \\
& = &
\sigma_{TX}(H)\big((\lam{k}{k(g)}) \after \sigma_{X}\big) \\
& = &
\mathcal{C}(\sigma_{X})(\sigma_{TX}(H))(\lam{k}{k(g)}) \\
& = &
\mu_{X}\Big(\mathcal{C}(\sigma_{X})(\sigma_{TX}(H))\Big)(g) \\
& = &
\big(\mu_{X} \after \mathcal{C}(\sigma_{X}) \after \sigma_{TX}\big)(H)(g).
\end{array}$$
}

\noindent Conversely, given a map of monads $\sigma\colon T
\Rightarrow C^{(C^{(-)})}$, define as algebra $a\colon T(C)
\rightarrow C$,
$$\begin{array}{rcl}
a(u)
& = &
\sigma_{C}(u)(\idmap[C]).
\end{array}\eqno{\QEDbox}$$

\auxproof{
This map $a\colon T(C) \rightarrow C$ is an algebra:
$$\begin{array}{rcl}
\big(a \after \eta_{C}\big)(x)
& = &
\sigma_{C}(\eta(x))(\idmap[C]) \\
& = &
\eta(x)(\idmap[C]) \\
& = &
\idmap[C](x) \\
& = &
x. \\
\big(a \after \mu_{C}\big)(v)
& = &
a(\mu_{C}(v)) \\
& = &
\sigma_{C}(\mu_{C}(v))(\idmap) \\
& = &
\big(\sigma_{C} \after \mu_{C}\big)(v)(\idmap) \\
& = &
\big(\mu_{C} \after \mathcal{C}(\sigma_{C}) \after \sigma_{TC}\big)
   (v)(\idmap) \\
& = &
\mu_{C}\Big(\mathcal{C}(\sigma_{C})(\sigma_{TC}(v))\Big)(\idmap) \\
& = &
\mathcal{C}(\sigma_{C})(\sigma_{TC}(v))(\lam{k}{k(\idmap)}) \\
& = &
\sigma_{TC}(v)\big((\lam{k}{k(\idmap)}) \after \sigma_{C}\big) \\
& = &
\sigma_{TC}(v)\big(\lam{u}{\sigma_{C}(u)(\idmap)}\big) \\
& = &
\sigma_{TC}(v)\big((\lam{u}{a(u)}\big) \\
& = &
\sigma_{TC}(v)(a) \\
& = &
a\big(T(a)(v)\big) \\
& = &
\big(a \after T(a)\big)(v)
\end{array}$$

Moreover, these translations are each others inverses:
$$\begin{array}{rcl}
\overline{\overline{a}}(u)
& = &
\overline{a}_{C}(u)(\idmap) \\
& = &
a\big(T(\idmap)(u)\big) \\
& = &
a(u) \\
\overline{\overline{\sigma}}_{X}(u)(g)
& = &
\overline{\sigma}\big(T(g)(u)\big) \\
& = &
\sigma_{C}\big(T(g)(u)\big)(\idmap) \\
& = &
\mathcal{C}(g)\big(\sigma_{X}(u)\big)(\idmap) \\
& = &
\sigma_{X}(u)(g).
\end{array}$$
}
\end{myproof}

Taking $C=2=\{0,1\}$ to be the two-element set, yields as associated
continuation monad $\mathcal{C}(X) = 2^{(2^{X})} \cong \Pow(\Pow(X))$,
the double-powerset monad. For a function $f\colon X\rightarrow Y$ we
have a map $\Pow^{2}(X) \rightarrow \Pow^{2}(Y)$, by functoriality,
given by double inverse image: $U\subseteq\Pow(X) \longmapsto
(f^{-1})^{-1}(U) = \set{V\subseteq Y}{f^{-1}(V) \in U}$.

It is not hard to see that the inclusion maps: 
$$\xymatrix{
\UF(X) \ar[r]^-{(\ref{UFDefEqn})}_-{\cong} &
   \BA(2^{X}, 2)\; \ar@{^(->}[r] & \; 2^{(2^{X})}
}$$

\noindent form a map of monads, from the ultrafilter monad to the
continuation monad (with constant $C = 2$).

\auxproof{
We describe functoriality in terms of associated characteristic
functions, $\chi_{U} \colon 2^{X} \rightarrow 2$ gives
$\mathcal{C}(f)(\chi_{U}) = \lam{\varphi\colon Y\rightarrow
  2}{\chi_{U}(\varphi \after f)}$.  Then, for $V\subseteq X$,
$$\begin{array}{rcl}
\mathcal{C}(f)(\chi_{U})(\chi_{V}) = 1
& \Longleftrightarrow &
\chi_{U}(\chi_{V} \after f) = 1 \\
& \Longleftrightarrow &
\set{x}{\chi_{V}(f(x)) = 1} \in U \\
& \Longleftrightarrow &
\set{x}{f(x)\in V} \in U \\
& \Longleftrightarrow &
f^{-1}(V) \in U.
\end{array}$$

Explicitly, if we write $\alpha_X$ for the maps $\UF(X) \rightarrow
2^{(2^{X})}$, given by $\alpha_{X}(\calF)(\varphi) = 1$ iff
$\set{x}{\varphi(x)=1}\in\calF$, then we see that $\alpha$ is a map
of monads:
$$\begin{array}{rcl}
\big(\alpha \after \eta\big)(x)(\varphi) = 1
& \Longleftrightarrow &
\alpha(\eta(x))(\varphi) = 1 \\
& \Longleftrightarrow &
\set{y}{\varphi(y)=1} \in \eta(x) \\
& \Longleftrightarrow &
x \in \set{y}{\varphi(y)=1} \\
& \Longleftrightarrow &
\varphi(x)= 1 \\
& \Longleftrightarrow &
\eta(x)(\varphi) = 1 \\
\big(\mu \after \mathcal{C}(\alpha) \after \alpha\big)
   (\mathcal{A})(\varphi) = 1
& \Longleftrightarrow &
\mu\Big(\mathcal{C}(\alpha)\big(\alpha(\mathcal{A})\big)\Big)
   (\varphi) = 1\\
& \Longleftrightarrow &
\mathcal{C}(\alpha)\big(\alpha(\mathcal{A})\big)
   \big(\lam{k}{k(\varphi)}\big) = 1 \\
& \Longleftrightarrow &
\alpha(\mathcal{A})\big((\lam{k}{k(\varphi)}) \after \alpha\big) = 1 \\
& \Longleftrightarrow &
\alpha(\mathcal{A})\big(\lam{\calF}{\alpha(\calF)(\varphi)}\big) = 1 \\
& \Longleftrightarrow &
\set{\calF}{\alpha(\calF)(\varphi)=1} \in \mathcal{A} \\
& \Longleftrightarrow &
\set{\calF}{\set{x}{\varphi(x)=1}\in\calF} \in \mathcal{A} \\
& \Longleftrightarrow &
D(\set{x}{\varphi(x)=1}) \in \mathcal{A} \\
& \Longleftrightarrow &
\set{x}{\varphi(x)=1} \in \mu(\mathcal{A}) \\
& \Longleftrightarrow &
\alpha(\mu(\mathcal{A}))(\varphi) = 1 \\
& \Longleftrightarrow &
\big(\alpha \after \mu\big)(\mathcal{A})(\varphi) = 1
\end{array}$$
}

\subsection{Monads from composable adjunctions}\label{MndAdjSubsec}

It is well-known, see \textit{e.g.}~\cite[Ch.~VI]{MacLane71} that each
adjunction $F\dashv G$ gives rise to a monad $GF$. The expectation
monad arises from a slightly more complicated situation, involving two
composable adjunctions. This situation is captured abstractly in the
following result.

\begin{lemma}
\label{ComposableAdjunctionLem}
Consider two composable adjunctions $F\dashv G$ and $H\dashv K$ in a
situation:
$$\xymatrix{
\cat{A}\ar@/^2ex/[rr]^-{F}
   \ar@(ur,ul)[]_{T = GF}\ar@(dr,dl)[]^{S = GKHF}  & \bot &
   \cat{B}\ar@/^2ex/[ll]^-{G}\ar@/^2ex/[rr]^-{H} & \bot &
   \cat{C}\ar@/^2ex/[ll]^-{K}
}$$

\noindent with monads $T = GF$ induced by the adjunction $F\dashv G$ 
and $S = GKHF$ induced by the (composite) adjunction $HF \dashv GK$.

Then there is a map of monads $T\Rightarrow S$ given by the
unit $\eta$ of the adjunction $H\dashv K$ in:
\begin{equation}
\label{ComposableAdjunctionNatroDiag}
\vcenter{\xymatrix{
T = GF\ar[rr]^-{G\eta^{H\dashv K} F} & & GKHF = S.
}}
\end{equation}

\noindent It gives rise a functor $\Alg(S) \rightarrow \Alg(T)$ between
the associated categories of Eilenberg-Moore algebras, and thus to
a commuting diagram:
\begin{equation}
\label{ComposableAdjunctionComparisonDiag}
\vcenter{\xymatrix@R-.5pc{
\cat{C}\ar[d]_{K}\ar[rr] & & \Alg(S)\ar[d]^{(-) \after G\eta F} \\
\cat{B}\ar[dr]_{G}\ar[rr] & & \Alg(T)\ar[dl]^{U} \\
& \cat{A} &
}}
\end{equation}

\noindent where the horizontal arrows are the so-called comparison
functors. 
\end{lemma}

\begin{myproof}
Easy. We unravel the relevant ingredients for future use. The unit and
counit of the composite adjunction $HF \dashv GK$ are:
$$\begin{array}{rcl}
\eta^{HF \dashv GK}
& = &
G\eta^{H\dashv K}F \after \eta^{F\dashv G} 
  \;\colon\; \idmap \Longrightarrow GKHF = S \\
\varepsilon^{HF \dashv GK}
& = &
\varepsilon^{H\dashv K} \after H\varepsilon^{F\dashv G} K 
  \;\colon\; HFGK \Longrightarrow \idmap.
\end{array}$$

\noindent This means that the monads $T$ and $S$ have multiplications:
$$\begin{array}{rcl}
\mu^{T}
& = &
G\varepsilon^{F\dashv G}F 
  \;\colon\; T^{2} = FGFG \Longrightarrow FG = T \\
\mu^{S}
& = &
GK\varepsilon^{H\dashv K}HF \after GKH\varepsilon^{F\dashv G} KHF 
  \;\colon\; S^{2} = GKHFGKHF \Longrightarrow GKHF = S.
\end{array}$$

\noindent The comparison functor $K_{T} \colon \cat{B} \rightarrow
\Alg(T)$ is:
$$\begin{array}{rcl}
K_{T}(X)
& = &
\big(TGX = GFGX\xrightarrow{G(\varepsilon^{F\dashv G}_{X})} GX\big).
\end{array}$$

\noindent Similarly, $K_{S} \colon \cat{A} \rightarrow \Alg(S)$ is:
$$\begin{array}{rcl}
K_{S}(Y)
& = &
\big(SGKY = GKHFGKY \xrightarrow{GK(\varepsilon^{H\dashv K}_{Y})
   \after GKH(\varepsilon^{F\dashv G}_{KY})} GKY\big).
\end{array}\eqno{\QEDbox}$$
\end{myproof}

\auxproof{ 
We write $\eta,\varepsilon$ for the unit and counit of the adjunction
$F\dashv G$ and $\eta',\varepsilon'$ for $H\dashv K$. The composite
adjunction $HF \dashv GK$ then has unit $G\eta'F \after \eta \colon
\idmap \Rightarrow GKHF$ and counit $\varepsilon' \after H\varepsilon K
\colon HFGK \Rightarrow \idmap$.

The monad $T = GF$ thus has multiplication $T^{2}\Rightarrow T$ given by:
$$\begin{array}{rcl}
\mu
& = &
G\varepsilon F\colon GFGF \Longrightarrow GF.
\end{array}$$

\noindent Similarly, the multiplication $S^{2} \Rightarrow S$ of 
$S$ is given by:
$$\begin{array}{rcl}
v\mu'
& = &
GK(\varepsilon' \after H\varepsilon K)HF 
   \colon GKHFGKHF \Longrightarrow GKHF \\
& = &
GK\varepsilon'HF \after GKH\varepsilon KHF 
\end{array}$$

The claim is that $\sigma = G\eta'F \colon T \Rightarrow S$ is a map
of monads. We check the equations involved.
$$\begin{array}{rcl}
\sigma \after \eta
& = &
G\eta'F \after \eta \\
& = &
\eta' \\
\mu' \after S\sigma \after \sigma T
& = &
GK\varepsilon'HF \after GKH\varepsilon KHF \after 
   GKHFG\eta'F \after G\eta'FGF \\
& = &
GK\varepsilon'HF \after GKH\eta'F \after GKH\varepsilon F \after 
   G\eta'FGF \\
& = &
GK\big(\varepsilon'H \after H\eta'\big)F \after 
   G\eta'F \after G\varepsilon F  \\
& = &
G\eta'F \after G\varepsilon F  \\
& = &
\sigma \after \mu
\end{array}$$

Finally, the comparison functor $K_{T} \colon \cat{B} \rightarrow \Alg(T)$
is:
$$\begin{array}{rcl}
K_{T}(X)
& = &
\big(TGX = GFGX\xrightarrow{G(\varepsilon_{X})} GX\big).
\end{array}$$

\noindent Similarly, $K_{S} \colon \cat{A} \rightarrow \Alg(S)$ is:
$$\begin{array}{rcl}
K_{S}(Y)
& = &
\big(SGKY = GKHFGKY \xrightarrow{GK(\varepsilon'_{Y})
   \after GKH(\varepsilon_{KY})} GKY\big).
\end{array}$$

\noindent The above square of comparison functors commutes, since
for $Y\in\cat{C}$ we get the same $T$-algebra $T(GKY) = GFGKY \rightarrow
GKY$ in:
$$\begin{array}{rcl}
\Big(\big((-)\after G\eta F\big) \after K_{S}\Big)(Y)
& = &
GK(\varepsilon'_{Y}) \after GKH(\varepsilon_{KY}) \after G(\eta'_{FGKY}) \\
& = &
GK(\varepsilon'_{Y}) \after G(\eta'_{KY}) \after G(\varepsilon_{KY}) \\
& = &
G(\varepsilon_{KY}) \\
& = &
\Big(K_{T} \after G\Big)(Y).
\end{array}$$
}

\begin{remark}
Later on in Section~\ref{ExpEModSec} we will construct a left adjoint
to the comparison functor $\cat{C} \rightarrow \Alg(S)$
in~\eqref{ComposableAdjunctionComparisonDiag}. It is already almost
there, in this abstract situation, using the composite adjunction
$HF\dashv KG$. However, suitable restrictions have to used, which
cannot be expressed at this abstract level. In the more concrete
setting described below, the adjunction $H\dashv K$ is of a special
kind, involving a dualizing object.
\end{remark}

The composable adjunctions that form the basis of the expectation
monad are:
\begin{equation}
\label{SetsConvEModDiag}
\vcenter{\xymatrix@R-1.5pc@C+.5pc{
\Sets\ar@/^2ex/[rr]^-{\Dst} & \bot &
   \Alg(\Dst)\ar@/^2ex/[ll]^-{U}
   \ar@/^2ex/[rr]^-{\Conv(-,[0,1])} & \bot &
   \EMod\op\ar@/^2ex/[ll]^-{\EMod(-,[0,1])} \\
& & \Conv\ar@{=}[u] & & 
}}
\end{equation}

\noindent The adjunction on the left is the standard adjunction
between a category of algebras $\Alg(\Dst)$ of the distribution monad
(see Subsection~\ref{DstSubsec}) and its underlying category. The
adjunction on the right will be described in the next section.

\section{Effect modules}\label{EModSec}

This section introduces the essentials of effect modules and refers
to~\cite{Jacobs10e,JacobsM12a} for further details. Intuitively,
effect modules are vector spaces, not with the real or complex numbers
as scalars, but with scalars from the unit interval $[0,1]\subseteq
\mathbb{R}$.  Also, the addition operation $+$ on vectors is only
partial; it is written as $\ovee$. These effect modules occur
in~\cite{PulmannovaG98} under the name `convex effect algebras'.

More precisely, an effect module is an \emph{effect algebra} $E$ with
an action $[0,1] \otimes E \rightarrow E$ for scalar
multiplication. An effect algebra $E$ carries both:
\begin{itemize}
\item a partial commutative monoid structure $(0,\ovee)$; this means
  that $\ovee$ is a partial operation $E\times E\rightarrow E$ which
  is both commutative and associative, taking suitably account of
  partiality, with $0$ as neutral element;

\item an orthosupplement $(-)^{\bot}\colon E \rightarrow E$. One
  writes $x\orthogonal y$ if the sum $x\ovee y$ is defined;
  $x^{\perp}$ is then the unique element with $x \ovee x^{\perp} = 1$,
  where $1 = 0^{\perp}$; further $x \orthogonal 1$ holds only for
  $x=0$.
\end{itemize}

\noindent These effect algebras carry a partial order given by $x\leq
y$ iff $x\ovee z = y$, for some element $z$. Then $x \orthogonal y$
iff $x \leq y^{\perp}$ iff $y\leq x^{\perp}$. The unit interval
$[0,1]$ is the prime example of an effect algebra with partial sum $r
\ovee s = r+s$ if $r+s \leq 1$; then $r^{\perp} = 1 - r$.

A homomorphism $f\colon E\rightarrow D$ of effect algebras satisfies
$f(1) = 1$ and: if $x\orthogonal x'$ in $E$, then $f(x) \orthogonal
f(x')$ in $D$ and $f(x\ovee x') = f(x) \ovee f(x')$. It is easy to
deduce that $f(x^{\perp}) = f(x)^{\perp}$ and $f(0) = 0$. This yields
a category, written as \EA. It carries a symmetric monoidal structure
$\otimes$ with the 2-element effect algebra $\{0,1\}$ as tensor unit
(which is at the same time the initial object). The usual
multiplication of real numbers (probabilities in this case) yields a
monoid structure on $[0,1]$ in the category \EA. An \emph{effect
  module} is then an effect algebra with an $[0,1]$-action
$[0,1]\otimes E \rightarrow E$. Explicitly, it can be described as a
scalar multiplication $(r,x) \mapsto rx$ satisfying:
$$\begin{array}{rclcrcll}
1x & = & x
& \qquad\quad &
(r+s)x & = & rx + sx \quad & \mbox{if }\;  r+s\leq 1 \\
(rs)x & = & r(sx)
& &
r(x\ovee y) & = & rx \ovee ry & \mbox{if }\; x\orthogonal y.
\end{array}$$

\noindent In particular, if $r+s\leq 1$, then a sum $rx \ovee sy$
always exists (see~\cite{PulmannovaG98}).

\auxproof{
Since $(r+s)1 = r1 \ovee s1$ we have $r1 \orthogonal s1$ and
thus $r1 \leq (s1)^{\perp}$. From $x\leq 1$ we get $rx \leq r1$,
and similarly $sy \leq s1$. Hence $rx \leq r1 \leq (s1)^{\perp}
\leq (sy)^{\perp}$. Hence $rx \orthogonal sy$.
}

\begin{example}
\label{EModEx}
The unit interval $[0,1]$ is again the prime example, this time for
effect modules. But also, for an arbitrary set $X$, the set
$[0,1]^{X}$ of all functions $X\rightarrow [0,1]$ is an effect module,
with structure inherited pointwise from $[0,1]$. Another example,
occurring in integration theory, is the set $[X\rightarrow_{s}[0,1]]$
of \textit{simple} functions $X\rightarrow [0,1]$, having only
finitely many output values (also known as `step functions').
\end{example}

A morphism $E\rightarrow D$ in the category \EMod of such effect
modules is a function $f\colon E\rightarrow D$ between the underlying
sets satisfying:
$$\begin{array}{rclcrclcrcll}
f(rx) & = & rf(x)
& \qquad &
f(1) & = & 1
& \qquad &
f(x\ovee y) & = & f(x) \ovee f(y)  & \;\mbox{ if }\; x\orthogonal y.
\end{array}$$

We now come to the dual adjunction mentioned in the previous section (see~\cite{JacobsM12a} 
for more information).

\begin{proposition}
\label{EModAdjProp}
For each effect module $E$ the homset $\EMod(E,[0,1])$ is a
  convex set. In the other direction, each convex set $X$ gives rise
  to an effect module $\Conv(X,[0,1])$. This gives the adjunction on
  the right in~\eqref{SetsConvEModDiag}, with $[0,1]$ as dualizing
  object
\end{proposition}

The effect algebra structure on the set $\Conv(X,[0,1])$ of affine
maps to $[0,1]$ is obtained pointwise: $f\ovee g$ is defined if
$f(x)+g(x) \leq 1$ for all $x\in X$, and in that case $f\ovee g$ at
$x\in X$ is $f(x)+g(x)$. The orthosupplement is also obtained
pointwise: $(f^{\perp})(x) = 1-f(x)$. Scalar multiplication is done
similarly $(rf)(x) = r(f(x))$. In the reverse direction, each effect
module $E$ gives rise to a convex set $\EMod(E,[0,1])$ of
homomorphisms, with pointwise convex sums. The adjunction
$\Conv(-,[0,1]) \dashv \EMod(-,[0,1])$ arises in the standard way,
with unit and counit given by evaluation.

\subsection{Totalization}

In this section we prove that the category of effect modules is
equivalent to the category of certain ordered vector spaces. For this
we extend a result for effect algebras from~\cite{JacobsM12a}. We
recall the basics below but for details and proofs we refer to that
paper. The idea is that the partial operation $\ovee$ of effect
algebras and effect modules is rather difficult to work with;
therefore we develop an embedding into structures with total
operations.

The first result we need is the following one from~\cite{JacobsM12a}.

\begin{proposition}
There is a coreflection
\begin{equation}
\label{EABCMEqn}
\xymatrix{
\EA\ar@/^/[rr]^{\partotscript}_\bot & & \BCM\ar@/^/[ll]^{\totparscript}
}
\end{equation}

\noindent where $\BCM$ is the category of ``barred commutative
monoids'': its objects are pairs $(M,u)$, where $M$ is a commutative
monoid and $u\in M$ is a unit such that $x+y=0$ implies $x=y=0$ and
$x+y=x+z=u$ implies $y=z$. The morphisms in $\BCM$ are monoid
homomorphisms that preserve the unit. As this is a coreflection every
effect algebra $E$ is isomorphic to $\totpar\partot(E)$. \QED
\end{proposition}

The partialization functor $\totpar$ in~\eqref{EABCMEqn} is defined by:
$$\begin{array}{rcl}
\totpar(M,u)
& = &
\setin{x}{M}{x\preceq u},
\end{array}$$

\noindent where $x\preceq y$ if{f} there exists a $z$ such that
$x+z=y$. The operation $\ojoin$ is defined by $x\ojoin y = x+y$ but
this is only defined if $x+y\preceq u$,
\textit{i.e.}~$x+y\in\totpar(M,u)$.

The totalization functor $\partot$ in~\eqref{EABCMEqn} is defined as:
$$\begin{array}{rcl}
\partot(E)
& = &
(\, \Mlt(E)/\!\sim, \; 1\cdot 1_E\,),
\end{array}$$

\noindent where $\Mlt(E)$ is the free commutative monoid on $E$,
consisting of all finite formal sums $n_1\cdot x_1+\cdots + n_m\cdot
x_m$, with $n_i\in \mathbb{N}$ and $x_i\in E$.  Here we identify sums
such as $1\cdot x+2\cdot x$ with $3\cdot x$. And $\sim$ is the
smallest monoid congruence such that $1\cdot x+1\cdot y \sim 1\cdot
(x\ojoin y)$ whenever $x\ojoin y$ is defined.

\begin{example}
Totalization of the truth values $\{0,1\}\in\EA$ and of the
probabilities $[0,1]\in\EA$ yields the natural numbers and the
non-negative reals:
$$\begin{array}{rclcrcl}
\partot(\{0,1\})
& \cong &
\mathbb{N} 
& \qquad\mbox{and}\qquad &
\partot([0,1]) 
& \cong 
&\mathbb{R}_{\geq 0}.
\end{array}$$
\end{example}

Recall that an effect module $E$ is just an effect algebra together
with a scalar product $[0,1]\otimes E\to E$. Now it turns out that
$\partot$ is a strong monoidal functor, and as a result
$\partot(E)\in\BCM$ comes equipped with a scalar product
$\mathbb{R}_{\geq 0}\otimes \partot(E)\to \partot(E)$. This gives the
monoid $\partot(E)$ the structure of a positive cone of some partially
ordered vector space. To make this exact we give the following
definition.

Construct a category $\cat{Coneu}$ as follows: its objects are pairs
$(M,u)$ where $M$ is a commutative monoid equipped with a scalar
product $\bullet:\mathbb{R}_{\geq 0}\times M\to M$ and $u\in M$ such
that the following axioms hold.
$$\begin{array}{rclcrcl}
1\bullet x 
& =  &
x 
& \qquad\qquad & 
(r+s)\bullet x 
& =  &
r\bullet x + s\bullet \\
(rs)\bullet x 
& =  &
r\bullet (s\bullet x)
& & 
r\bullet (x+ y) 
& =  & 
r\bullet x + r\bullet y \\
x+y
& = & 
0 \;\text{ implies }\; x=y=0  
& &
x+y
\hspace*{\arraycolsep} = \hspace*{\arraycolsep} x+z 
& = & 
u \;\text{ implies }\; y=z, 
\end{array}$$

\noindent and for all $x\in M$ there exists an $n\in\mathbb{N}$ such
that $x\preceq n\bullet u$. Because of this last property we call $u$
a \emph{strong} unit.  The morphisms of $\cat{Coneu}$ are monoid
homomorphisms that respect both the scalar multiplication and the
unit.

We can then extend the coreflection $\partot\dashv\totpar$ to the
categories $\EMod$ and $\cat{Coneu}$. This will actually be an
equivalence of categories. To prove this we first need an
auxiliary result.

\begin{lemma} 
\label{L:ConeuCancel}
If $M\in\cat{Coneu}$ then the cancelation law holds in $M$.
\end{lemma}

\begin{myproof}
Let $x,y,z\in M$ and suppose $x+y=x+z$. Since $u$ is a strong unit we
can find an $n$ such that $x+y\preceq nu$. Therefore
$$\begin{array}{rcccl}
\frac{1}{n}\bullet x+ \frac{1}{n} \bullet y 
& = &
\frac{1}{n}\bullet x+\frac{1}{n}\bullet z 
& \preceq &
u.
\end{array}$$

\noindent Hence we can find an element $w\in M$ such that
$\frac{1}{n}\bullet x+\frac{1}{n}\bullet y+w=\frac{1}{n}\bullet
x+\frac{1}{n}\bullet z+w=u$.  Then $\frac{1}{n}\bullet y
=\frac{1}{n}\bullet z$. And thus $y = \sum_{i=1}^n \frac{1}{n}\bullet
y=\sum_{i=1}^n \frac{1}{n}\bullet z = z$. \QED
\end{myproof}

An immediate consequence is that the preorder $\preceq$ is a partial
order; thus we shall write $\leq$ instead of $\preceq$ from now on.

\begin{lemma} \label{L:[0,1]equiv1}
The coreflection $\partot\dashv\totpar$ between $\EMod$ and
$\cat{Coneu}$ is an equivalence of categories.
\end{lemma}

\begin{myproof}
We only need to show that the counit of the adjunction
$\partot\dashv\totpar$ is an isomorphism.  So let $M\in\cat{Coneu}$; a
typical element of $\partot\totpar (M)$ is an equivalence class of
formal sums like $\sum n_ix_i$ where $n_i\in \mathbb{N}$ and $M\ni
x_i\leq u$. The counit $\eps$ sends the class represented by this
formal sum to its interpretation as an actual sum in $M$.

To show that $\eps$ is surjective suppose $x\in M$. We can find a
natural number $n$ such that $x\leq nu$ so that $\frac{1}{n}\bullet x
\leq u$. This gives us:
$$\begin{array}{rcccl}
x 
& = &
n\cdot (\frac{1}{n}\bullet x)
& = &
\eps(n(\frac{1}{n}\bullet x)).
\end{array}$$

To prove injectivity suppose that $\eps(\sum n_ix_i) = \eps(\sum
k_jy_j)$. Define $N=\sum n_i+\sum k_j$, so that:
$$\begin{array}{rcccccccccl}
\sum n_i\cdot (\frac{1}{N}\bullet x_i) 
& = &
\eps (\sum n_i(\frac{1}{N}\bullet x_i)) 
& = &
\eps(\frac{1}{N}\bullet(\sum n_ix_i)) 
& = &
\frac{1}{N}\bullet\eps(\sum n_ix_i) 
& = &
\frac{1}{N}\bullet\eps(\sum k_jy_j) 
& = &
\sum k_j(\frac{1}{N}\bullet y_j).
\end{array}$$

\noindent Because $N$ is sufficiently large, the terms $\ojoin_{i}\,
n_i\cdot (\frac{1}{N}\bullet x_i)$ and $\ojoin_{j}\,
k_j\cdot(\frac{1}{N}\bullet y_j)$ are both defined in $\totpar M$ and
by the previous calculation they are equal. This means that $\sum
n_i(\frac{1}{N}\bullet x_i)$ and $\sum k_j(\frac{1}{N}\bullet y_k)$
represent equal elements of $\partot\totpar M$ and therefore
the equation
$$\begin{array}{rcccccl}
\sum n_ix_i 
& = &
N\bullet (\sum n_i(\frac{1}{N}\bullet x_i)) 
& = &
N\bullet (\sum k_j(\frac{1}{N}\bullet y_j)) 
& = &
\sum k_jy_j.
\end{array}$$

\noindent holds in $\totpar\partot M$. \QED
\end{myproof}

From positive cones it is but a small step to partially ordered vector
spaces.  Define a category \poVectu as follows; the objects are
partially ordered vector spaces over $\mathbb{R}$ with a strong order
unit $u$, \textit{i.e.}~a positive element $u\in V$ such that for any
$x\in V$ there is a natural number $n$ with $x\leq nu$. The morphisms
in \poVectu are linear maps that preserve both the order and the unit.

\begin{theorem}\label{T:[0,1]equiv2}
The category $\EMod$ is equivalent to $\poVectu$.
\end{theorem}

\begin{myproof}
We will prove that $\cat{Coneu}$ is equivalent to $\poVectu$; the
result then follows from lemma~\ref{L:[0,1]equiv1}.

The functor $F\colon\poVectu\to\cat{Coneu}$ takes the positive cone of
a partially ordered vector space. The construction of
$G\colon\cat{Coneu}\to\poVectu$ is essentially just the usual
construction of turning a cancellative monoid into a group.

In somewhat more detail: if $M\in\cat{Coneu}$ then define $G(M) =
(M\times M)/\sim$ where $\sim$ is defined by $(x,y)\sim(x',y')$ if{f}
$x+y'=y+x'$. We write $[x,y]$ for the equivalence class of $(x,y)\in
M\times M$. Addition is defined by $[x,y]+[x',y'] = [x+x',y+y']$. If
$\alpha\in\reals$ we define $\alpha\scalar[x,y]$ as follows. If
$\alpha\geq 0$ then $\alpha[x,y] = [\alpha \scalar x,\alpha\scalar y]$
and if $\alpha<0$ then $\alpha[x,y] = [-\alpha\scalar y,-\alpha\scalar
  x]$.  It's easy to check that $G(M)$ is indeed a vector
space. Moreover, $G(M)$ is partially ordered by $[x,y]\leq [x',y']$
if{f} $x+y'\leq y+x'$, and $[u,0]$ is its strong unit.

It's easy to see that both constructions can be made functorial and
that this gives an equivalence of categories. \QED
\end{myproof}

We write $\pth:\EMod\leftrightarrows \poVectu:\tph$ for this
equivalence. For a partially ordered vector space $V$ with a strong
unit $u$ the effect module $\tph(V)$ consists of all elements $x$ such
that $0\leq x\leq u$. With this equivalence of categories in hands we
can apply techniques from linear algebra to effect modules. Below we
translate some properties of partially ordered vector spaces to the
language of effect modules. We need these results later on.

If $V\in\poVectu$ and the unit $u$ is Archimedean---in the sense that
$x\leq ru$ for all $r>0$ implies $x\leq 0$---then $V$ is called an
\textit{order unit space}. The Archimedean property of the unit can be
used to define a norm $\|x\| = \inf\setin{r}{[0,1]}{-ru\leq x\leq
  ru}$. We denote by $\OUS$ the full subcategory of $\poVectu$
consisting of all order unit spaces.

This Archimedean property can also be expressed on the effect module level
but some caution is required as effect modules contain no elements
less than $0$ and sums may not be defined. The following formulation
works: an effect module is said to be Archimedean if $x\leq y$ follows
from $\frac{1}{2}x\leq\frac{1}{2}y \ovee \frac{r}{2}1$ for all $r\in
(0,1]$.  All Archimedean effect modules form a full subcategory
  $\AEMod\hookrightarrow \EMod$.  Of course with this definition comes
  a theorem.

\begin{proposition}
\label{ArchimideanEquivProp}
The equivalence $\poVectu\simeq\EMod$, between partially ordered vector
spaces with a strong unit and effect modules, restricts to an
equivalence $\OUS\simeq\AEMod$, between order unit spaces and
Archimedean effect modules.
\end{proposition}

\begin{myproof}
We only check that if $E\in\AEMod$ then its totalization satisfies
$\pth(E)\in\OUS$; the rest is left to the reader.  Suppose
$E\in\AEMod$ and $x\in \pth(E)$ is such that $x\leq ru$ for all
$r\in(0,1]$. The trick is to transform $x$ into an element in the unit
  interval $[0,u]\cong E$.  Since $u$ is a strong unit we can find a
  natural number $n$ such that $x+nu\geq 0$, and again using the fact
  that $u$ is a strong unit we can find a positive real number $s<1$
  such that $sx+nsu\leq u$. Hence $sx+nsu\in [0,u]\cong E$. Now, for
  $r\in (0,1]$ we have $sx \leq x\leq ru$ and so
    $\frac{s}{2}x+\frac{ns}{2}u\leq \frac{ns}{2}u+\frac{r}{2}u$. Thus,
    by the Archimedean property of $E$, we get $sx+nsu\leq nsu$. Hence
    $sx\leq 0$ and therefore $x\leq 0$. \QED
\end{myproof}

Since $E\in\AEMod$ is isomorphic to the unit interval of its
totalization $\pth(E)$, $E$ inherits a metric from the normed space
$\pth(E)$. This metric can be described wholly in terms of
$E$. However the partial addition does force us into a somewhat
awkward definition: for $x,y\in E$ their distance $d(x,y)\in[0,1]$ can
be defined as:
\begin{equation}
\label{EModDistEqn}
\begin{array}{rcl}
d(x,y)
& = &
\max\Big(\,\inf\setin{r}{(0,1]}{\frac{1}{2}x \leq 
   \frac{1}{2}y \ovee\frac{r}{2}1}, \;
  \inf\setin{r}{(0,1]}{\frac{1}{2}y \leq \frac{1}{2}x \ovee \frac{r}{2}1}
  \,\Big).
\end{array}
\end{equation}

A trivial consequence is the following lemma.

\begin{lemma}
\label{AEModNonExpLem}
A map of effect modules $f\colon M\rightarrow M'$ between Archimedean
effect modules $M,M'$ is automatically non-expansive: $d'(f(x),f(y))
\leq d(x,y)$, for all $x,y\in M$. \QED
\end{lemma}

Of particular interest later in this paper are Archimedean effect
modules that are complete in their metric. We call these
\textit{Banach effect modules} and denote by $\BEMod$ the full
subcategory of all Banach effect modules. The previous lemma implies
that each map in $\BEMod$ is automatically continuous.

Since an order unit space is complete in its metric if and only if its
unit interval is complete we get the following result.

\begin{proposition}
The equivalences from Proposition~\ref{ArchimideanEquivProp} restrict
further to an equivalence between Banach effect modules and the full
subcategory $\BOUS\hookrightarrow\OUS$ of those order unit spaces that
are also Banach spaces:
$$\xymatrix@R-.5pc{
\BOUS\ar@{^(->}[d]\ar[rr]^-{\simeq} & & \BEMod\ar@{^(->}[d] \\
\OUS\ar@{^(->}[d]\ar[rr]^-{\simeq} & & \AEMod\ar@{^(->}[d] \\
\poVectu\ar[rr]^-{\simeq}_-{V\mapsto([0,u]\subseteq V)} &  & \EMod
}$$
\end{proposition}

\begin{myproof}
Like in the proof of Proposition~\ref{ArchimideanEquivProp} one
transforms a Cauchy sequence in $\pth(E)$ into a sequence
in $[0,u]\cong E$. \QED
\end{myproof}

\auxproof{ 
Suppose $E\in\BEMod$ and let $(x_n)_{n\in\NNO}$ be a Cauchy sequence
in $\pth(E)$.  By adding some multiple of $u$ to this sequence we can
assume that $x_n\geq 0$. Then by dividing by some constant we can
also assume $x_n\in[0,u]\cong E$. But then the sequence is contained 
in the unit interval and thus converges by completeness of $E$.  
}

%

\begin{example}
\label{AEModEx}
We review Example~\ref{EModEx}: both the effect modules $[0,1]$ and
$[0,1]^{X}$ are Archimedean, and also Banach effect modules.  Norms
and distances in $[0,1]$ are the usual ones, but limits in $[0,1]^{X}$
are defined via the supremum (or uniform) norm: for $p\in[0,1]^{X}$,
we have:
$$\begin{array}{rcl}
\|p\|
& = &
\inf\setin{r}{[0,1]}{p \leq r\cdot u} 
 \quad \mbox{where $u$ is the constant function $\lam{x}{1}$} \\
& = &
\inf\setin{r}{[0,1]}{\allin{x}{X}{p(x) \leq r}} \\
& = &
\sup\set{p(x)}{x\in X} \\
& = &
\|p\|_{\infty}.
\end{array}$$

\noindent The latter notation $\|p\|_{\infty}$ is common for this
supremum norm. The associated metric on $[0,1]^{X}$ is 
according to~\eqref{EModDistEqn}:
$$\begin{array}{rcl}
d(p,q)
& = &
\max\Big(\begin{array}[t]{l}
   \inf\setin{r}{(0,1]}{\allin{x}{X}{\frac{1}{2}p(x) \leq 
   \frac{1}{2}q(x)+ \frac{r}{2}}},
  \\
  \inf\setin{r}{(0,1]}{\allin{x}{X}{\frac{1}{2}q(x) \leq 
   \frac{1}{2}p(x)+ \frac{r}{2}}}\,\Big).
  \end{array} \\
& = &
\max\Big(\begin{array}[t]{l}
   \sup\set{p(x) - q(x)}{x\in X\mbox{ with }p(x) \geq q(x)}, \\
   \sup\set{q(x) - p(x)}{x\in X\mbox{ with }p(x) \leq q(x)} \,\Big)
   \end{array} \\
& = &
\sup\set{|p(x)-q(x)|}{x\in X} \\
& = &
\|p-q\|_{\infty}.
\end{array}$$

Recall that the subset $[X\rightarrow_{s}[0,1]]\subseteq [0,1]^{X}$ of
simple functions contains those $p\in[0,1]^{X}$ that take only
finitely many values, \textit{i.e.}~for which the set $\set{p(x)}{x\in
  X}$ is finite. If we write $\set{p(x)}{x\in X} = \{r_{1}, \ldots,
r_{n}\} \subseteq [0,1]$, then we obtain $n$ disjoint non-empty sets
$X_{i} = \setin{x}{X}{p(x) = r_{i}}$ covering $X$. For a subset
$U\subseteq X$, let $\charac{U} \colon X\rightarrow [0,1]$ be the
corresponding ``characteristic'' simple function, with $\charac{U}(x)
= 1$ iff $x\in U$ and $\charac{U}(x) = 0$ iff $x\not\in U$. Hence we
can write such a simple function $p$ in a normal form in the effect
module $[X\rightarrow_{s}[0,1]]$ of simple functions, namely as finite
sum of characteristic functions:
\begin{equation}
\label{SimpleFunNormalFormEqn}
\begin{array}{rcl}
p
& = &
\ovee_{i}\; r_{i}\cdot \charac{X_i}.
\end{array}
\end{equation}

\noindent Hence $\|p\| = \max\{r_{1}, \ldots, r_{n}\}$. These
simple functions do not form a Banach effect module, since simple
functions are not closed under countable suprema.
\end{example}

\begin{lemma}
\label{SimplePredLem}
The inclusion of simple functions on a set $X$ is dense in
the Banach effect module of all fuzzy predicates on $X$:
$$\xymatrix@C+1pc{
[X\rightarrow_{s}[0,1]] \ar@{ >->}[r]^-{\textrm{dense}} & [0,1]^{X}
}$$

\noindent Explicitly, each predicate $p\in [0,1]^{X}$ can be written
as limit $\smash{p = {\displaystyle\lim_{n\rightarrow\infty}}\,p_{n}}$
of simple functions $p_{n}\in [0,1]^{X}$ with $p_{n} \leq p$. 
\end{lemma}

\begin{myproof}
Define for instance:
$$\begin{array}{rclcrcl}
p_{n}(x)
& = &
0.d_{1}d_{2}\cdots d_{n}
\qquad\mbox{where}\qquad
d_{i}
& = &
\mbox{the $i$-th decimal of $p(x)\in[0,1]$.}
\end{array}$$

\noindent Clearly, $p_{n}$ is simple, because it can take at most
$10^{n}$ different values, since $d_{i}\in\{0,1,\ldots,9\}$. Also,
by construction, $p_{n} \leq p$. For each $\epsilon > 0$, take
$N\in\NNO$ such that for all decimals $d_i$ we have:
$$\begin{array}{rcl}
0.\underbrace{00\cdots 00}_{N\text{ times}}d_{1}d_{2}d_{3}\cdots
& < &
\epsilon.
\end{array}$$

\noindent Then for each $n\geq N$ we have $p(x)-p_{n}(x) < \epsilon$,
for all $x\in X$, and thus $d(p,p_{n}) \leq \epsilon$. \QED
\end{myproof}

\subsection{Hahn-Banach style extension for effect modules}

In this subsection we look at a form of Hahn-Banach theorem for effect
modules.  We need the following version of the Hahn-Banach extension
theorem for partially ordered vector spaces.

\begin{proposition}
\label{HBCofinalProp}
Let $E$ be a partially ordered vector space and let $F\subseteq E$ be a
cofinal subspace (\textit{i.e.}~for all $x\in E$, $x\geq 0$ there is
$y\in F$ with $x\leq y$).  Suppose $f\colon F\to\reals$ is a monotonic
linear function. Then there is a monotonic linear function
$g\colon E\to \reals$ with $g|_{F}=f$ in:
$$\xymatrix@R-.5pc{
F\;\ar[dr]_{f}\ar@{^(->}[r] & E\ar@{-->}[d]^{\exists g} \\
& \reals
}$$
\end{proposition}

\begin{myproof}
We define $p\colon E\to \reals$ by:
$$\begin{array}{rcl}
p(x)
& = &
\inf \set{f(y)}{y\in F\text{ and } y\geq x}.
\end{array}$$

\noindent Notice that $p(x)$ is finite since we can find $y,y'\in F$
with $y\leq x\leq y'$ because $F$ is cofinal.  We need to check that
$p$ is sublinear. So let $x,x'\in E$ and $\epsilon>0$ then we can find
$y,y'\in F$ with $y\geq x,y'\geq x'$, such that $f(y)
<p(x)+\epsilon>0$ and $f(y')<p(x')+\epsilon>0$. Therefore:
$$\begin{array}{rcccccccl}
p(x+x') & \leq & f(y+y') & = & f(y)+f(y') & < & p(x)+p(x')+2\epsilon
  & > & 0.
\end{array}$$

\noindent Also for $r>0$ it is obvious that $p(r \cdot x)=r \cdot
p(x)$.

Having established that $p$ is sublinear we note that for $y\in F$ we
have $p(y)=f(y)$ since $f$ is monotonic. Hence we can apply the
standard (dominated) extension version of Hahn-Banach to find a linear
function $g\colon E\to\reals$ with $g<p$ and $g|_F=f$. Since if $x\leq
0$ then $p(x)\leq 0$ because $0\in F$, hence it follows that $g$ is
monotonic. \QED
\end{myproof}

This version translates effortlessly to effect modules

\begin{proposition}
\label{HBEModProp}
Let $E$ be an effect module and $F\subseteq E$ a sub effect module of
$E$.  Suppose $f\colon F\to[0,1]$ is an effect module map, then there
is an effect module map $g\colon E\to[0,1]$ with $g|_F=E$.
\end{proposition}

\begin{myproof}
We translate effect modules to order unit spaces and apply the
previous result.  Since $u\in\pth(F)$ it's clear that $\pth(F)$ is
cofinal in $\pth(E)$.  Hence using the previous proposition we can
extend $\pth(f)$ to $h\colon \pth(E)\to\reals$. Hence by restriction
to unit intervals $[0,u]$, both in $\pth(E)$ and in $\reals$ we get
the map $g\colon E\rightarrow [0,1]$ that we are looking for. \QED
\end{myproof}

Unfortunately the class of effect module morphisms is too limited to
get a full version of the separation theorem. Consider for example
$E=[0,1]^2$ with the two compact convex subsets
$C_1=\set{(r,\frac{1}{2}+r)}{r \in [0,\frac{1}{2}]}$ and
$C_2=\set{(\frac{1}{2}+r,r)}{r \in [0,\frac{1}{2}]}$. If $f\colon
E\to[0,1]$ is an effect module morphism then the image $f(C_1)$ is the
interval $[f(0,\frac{1}{2}),f(\frac{1}{2},1)]$, and since
$f(\frac{1}{2},\frac{1}{2})=\frac{1}{2}$ it follows that this interval
has length $\frac{1}{2}$. Analogously the interval $f(C_2)$ is also an
interval of length $\frac{1}{2}$ so the two must overlap.

\section{The expectation monad}\label{ExpMndSec}

We now apply Lemma~\ref{ComposableAdjunctionLem} to the composable
adjunctions in~\eqref{SetsConvEModDiag} and take a first look at the
results. In particular, we investigate different ways of describing
the expectation monad $\Exp$ that arises in this way.

Of the two monads resulting from applying
Lemma~\ref{ComposableAdjunctionLem} to the composable adjunctions in
Diagram~\eqref{SetsConvEModDiag}, the first one is the well-known
distribution monad $\Dst$ on \Sets, arising from the adjunction $\Sets
\leftrightarrows \Alg(\Dst) = \Conv$. The second monad on \Sets arises
from the composite adjunction $\Sets \leftrightarrows \EMod\op$ is
less familiar (see Section~\ref{SemanticsSec} for more information and
references). It is what we call the \emph{expectation monad}, written
here as $\Exp$. Following the description in
Lemma~\ref{ComposableAdjunctionLem} this monad is:
$$\begin{array}{rcl}
X 
& \longmapsto &
\EMod\Big(\Conv\big(\Dst(X), [0,1]\big), \; [0,1]\Big).
\end{array}$$

\noindent Since $\Dst\colon\Sets\rightarrow\Alg(\Dst) = \Conv$ is the
free algebra functor, the homset $\Conv(\Dst(X), [0,1])$ is isomorphic
to the set $[0,1]^{X}$ of all maps $X \rightarrow [0,1]$ in
\Sets. Elements of this set $[0,1]^{X}$ can be understood as fuzzy
predicates on $X$.  As mentioned, they form a Banach effect module via
pointwise operations. Thus we describe the expectation monad
$\Exp\colon\Sets\rightarrow\Sets$ as:
\begin{equation}
\label{ExpMndEqn}
\begin{array}{rcl}
\Exp(X)
& = &
\EMod\big([0,1]^{X},[0,1]\big) \\
\Exp\big(X\xrightarrow{f} Y\big)
& = &
\lamin{h}{\Exp(X)}{\lamin{p}{[0,1]^{Y}}{h(p \after f)}}.
\end{array}
\end{equation}

\noindent The unit $\eta_{X}\colon X\rightarrow \Exp(X)$ is given by:
$$\begin{array}{rcl}
\eta_{X}(x)
& = &
\lamin{p}{[0,1]^{X}}{p(x)}.
\end{array}$$

\noindent And the multiplication $\mu_{X} \colon \Exp^{2}(X) \rightarrow
\Exp(X)$ is given on $h\colon [0,1]^{\Exp(X)} \rightarrow [0,1]$ in
\EMod by:
$$\begin{array}{rcl}
\mu_{X}(h)
& = &
\lamin{p}{[0,1]^{X}}{h\Big(\lamin{k}{\Exp(X)}{k(p)}\Big)}.
\end{array}$$

\noindent It is not hard to see that $\eta(x)$ and $\mu(h)$ are
homomorphisms of effect modules. We check explicitly that the
$\mu\mbox{-}\eta$ laws hold and leave the remaining verifications to the
reader. For $h\in\Exp(X)$,
$$\begin{array}{rcl}
\big(\mu_{X} \after \eta_{\Exp(X)}\big)(h)
& = &
\mu_{X}\big(\eta_{\Exp(X)}(h)\big) \\
& = &
\lamin{p}{[0,1]^{X}}{\eta_{\Exp(X)}(h)\big(\lamin{k}{\Exp(X)}{k(p)}\big)} \\
& = &
\lamin{p}{[0,1]^{X}}{\big(\lamin{k}{\Exp(X)}{k(p)}\big)(h)} \\
& = &
\lamin{p}{[0,1]^{X}}{h(p)} \\
& = &
h \\
\big(\mu_{X} \after \Exp(\eta_{X})\big)(h)
& = &
\mu_{X}\big(\Exp(\eta_{X})(h)\big) \\
& = &
\lamin{p}{[0,1]^{X}}{\Exp(\eta_{X})(h)\big(\lamin{k}{\Exp(X)}{k(p)}\big)} \\
& = &
\lamin{p}{[0,1]^{X}}{h\big((\lamin{k}{\Exp(X)}{k(p)}) \after \eta_{X}\big)} \\
& = &
\lamin{p}{[0,1]^{X}}{h\big(\lamin{x}{X}{\eta_{X}(x)(p)}\big)} \\
& = &
\lamin{p}{[0,1]^{X}}{h\big(\lamin{x}{X}{p(x)}\big)} \\
& = &
\lamin{p}{[0,1]^{X}}{h(p)} \\
& = &
h.
\end{array}$$

\auxproof{
Just to be sure we check that the unit and multiplication are 
well-defined, natural, and satisfy the monad equations. 

The map $\eta(x) = \lam{p}{p(x)}$ is a obviously map of effect modules
$[0,1]^{X} \rightarrow [0,1]$. Similarly, $\mu(h) =
\lam{p}{h(\lam{k}{k(p)})}$ is a map of effect monoids $[0,1]^{\Exp(X)}
\rightarrow [0,1]$; the action properties are obvious. And if
$p,q\colon X \rightarrow [0,1]$ satisfy $p\orthogonal q$, then $p(x) +
q(x) \leq 1$, for all $x\in X$.  Then for each $k\in\Exp(X)$,
$k(p\ovee q) = k(p) \ovee k(q)$.  Since $[0,1]^{\Exp(X)}$ has
pointwise structure we get $\lam{k}{k(p\ovee q)} = \big(\lam{k}{k(p)}
\ovee \lam{k}{k(q)}\big)$. Thus for a homomorphism $h\colon
[0,1]^{\Exp(X)} \rightarrow [0,1]$,
$$\begin{array}{rcl}
\mu(h)(p\ovee q)
& = &
h\Big(\lamin{k}{\Exp(X)}{k(p\ovee q)}\Big) \\
& = &
h\Big(\lamin{k}{\Exp(X)}{k(p)}\ovee \lamin{k}{\Exp(X)}{k(q)}\Big) \\
& = &
h\Big(\lamin{k}{\Exp(X)}{k(p)}\Big)\ovee h\Big(\lamin{k}{\Exp(X)}{k(q)}\Big)
   \qquad \mbox{since $h$ is a homomorphism} \\
& = &
\mu(h)(p) \ovee \mu(h)(q) \\
\mu(h)(0)
& = &
h\big(\lamin{k}{\Exp(X)}{k(0)}\big) \\
& = &
h\big(\lamin{k}{\Exp(X)}{0}\big) \\
& = &
h(0) \\
& = &
0.
\end{array}$$

For naturality, assume $f\colon X\rightarrow Y$ in $\Sets$; then:
$$\begin{array}{rcl}
\big(\Exp(f) \after \eta_{X}\big)(x)
& = &
\Exp(f)(\eta_{X}(x)) \\
& = &
\lam{p}{\eta_{X}(x)(p\after f)} \\
& = &
\lam{p}{p(f(x))} \\
& = &
\eta_{Y}(f(x)) \\
& = &
\big(\eta_{Y} \after f\big)(x) \\
\big(\mu_{Y} \after \Exp^{2}(f)\big)(h) 
& = &
\mu_{Y}\big(\Exp(\Exp(f))(h)\big) \\
& = &
\lam{p}{\Exp(\Exp(f))(h)(\lam{k}{k(p)})} \\
& = &
\lam{p}{h\Big(\big(\lam{k}{k(p)}\big) \after \Exp(f)\Big)} \\
& = &
\lam{p}{h\Big(\lam{k}{\Exp(f)(k)(p)}\Big)} \\
& = &
\lam{p}{h\Big(\lam{k}{k(p \after f)}\Big)} \\
& = &
\lam{p}{\mu_{X}(h)(p \after f)} \\
& = &
\Exp(f)(\mu_{X}(h)) \\
& = &
\big(\Exp(f) \after \mu_{X}\big)(h).
\end{array}$$

Finally we check the $\mu$-equation. For $H\in \Exp^{2}(X) = 
\Conv([0,1]^{\Exp(X)}, [0,1])$:
$$\begin{array}{rcl}
\big(\mu_{X} \after \Exp(\mu_{X})\big)(H)
& = &
\mu_{X}\big(\Exp(\mu_{X})(H)\big) \\
& = &
\lamin{p}{[0,1]^{X}}{\Exp(\mu_{X})(H)\Big(\lamin{k}{\Exp(X)}{k(p)}\Big)} \\
& = &
\lamin{p}{[0,1]^{X}}{H\Big((\lamin{k}{\Exp(X)}{k(p)}) \after \mu_{X}\Big)} \\
& = &
\lamin{p}{[0,1]^{X}}{H\Big(\lamin{h}{\Exp^{2}(X)}
   {\mu_{X}(h)(p)}\Big)} \\
& = &
\lamin{p}{[0,1]^{X}}{H\Big(\lamin{h}{\Exp^{2}(X)}
   {h(\lamin{k}{\Exp(X)}{k(p)})}\Big)} \\
& = &
\lamin{p}{[0,1]^{X}}{\mu_{\Exp(X)}(H)\Big(\lamin{k}{\Exp(X)}{k(p)}\Big)} \\
& = &
\mu_{X}\big(\mu_{\Exp(X)}(H)\big) \\
& = &
\big(\mu_{X} \after \mu_{\Exp(X)}\big)(H).
\end{array}$$
}

\begin{remark}
\label{UltrafilterRem}
(1)~We think of elements $h\in\Exp(X)$ as measures. Later on, in
  Theorem~\ref{FinAddMeasThm}, it will be proven that $\Exp(X)$ is
  isomorphic to the set of finitely additive measures
  $\Pow(X)\rightarrow [0,1]$ on $X$. The application $h(p)$ of
  $h\in\Exp(X)$ to a function $p\in[0,1]^{X}$ may then be understood
  as integration $\int p\,dh$, giving the expected value of the
  stochastic variable/predicate $p$ for the measure $h$.

(2)~The description $\Exp(X) = \EMod\big([0,1]^{X},[0,1]\big)$ of the
expectation monad in~\eqref{ExpMndEqn} bears a certain formal
resemblance to the ultrafilter monad $\UF$ from
Subsection~\ref{UFSubsec}.  Recall from~\eqref{UFDefEqn} that:
$$\begin{array}{rcl}
\UF(X)
& \cong &
\BA\big(\{0,1\}^{X}, \; \{0,1\}\big).
\end{array}$$

\noindent Thus, the expectation monad $\Exp$ can be seen as a ``fuzzy''
or ``probabilistic'' version of the ultrafilter monad $\UF$, in which
the set of Booleans $\{0,1\}$ is replaced by the set $[0,1]$ of
probabilities. The relation between the two monads is further
investigated in Section~\ref{ExpUFSec}.

(3)~Using the equivalence $\poVectu\simeq\EMod$ from
Proposition~\ref{ArchimideanEquivProp} via totalization we may
equivalently describe the expectation monad as the homset:
$$\begin{array}{rcl}
\Exp(X)
& \cong &
\poVectu\big(\reals^{X}, \, \reals\big).
\end{array}$$

\noindent It contains the linear monotone functions $\reals^{X}
\rightarrow \reals$ that send the unit $\lam{x}{1}\in\reals^{X}$ to
$1\in\reals$.
\end{remark}

The following result is not a surprise, given the resemblance between
the unit and multiplication for the expectation monad and the ones for
the continuation monad (see Subsection~\ref{ContinuationSubsec}).

\begin{lemma}
\label{ExpContinuationLem}
The inclusion maps:
$$\xymatrix{
\Exp(X) \;=\; \EMod\big([0,1]^{X},[0,1]\big) \;
   \ar@{^(->}[r] & \; [0,1]^{([0,1]^{X})}
}$$

\noindent form a map of monads, from the expectation monad to the
continuation monad (with the set $[0,1]$ as constant). \QED
\end{lemma}

We conclude with an alternative description of the sets $\Exp(X)$, in
terms of finitely additive measures, described as effect algebra
homomorphisms. It also occurs as~\cite[Cor.~4.3]{Gudder98}.

\begin{theorem}
\label{FinAddMeasThm}
For each set $X$ there is a bijection:
$$\xymatrix{
\Exp(X) = \EMod\big([0,1]^{X},\,[0,1]\big)\ar[r]^-{\Phi}_-{\cong} &
   \EA\big(\Pow(X), [0,1]\big)
}$$

\noindent given by $\Phi(h)(U) = h(\charac{U})$.
\end{theorem}

\begin{myproof}
We first check that $\Phi$ is injective: assume $\Phi(h) = \Phi(h')$,
for $h,h'\in\Exp(X)$.  We need to show $h(p) = h'(p)$ for an arbitrary
$p\in[0,1]^{X}$.  We first prove $h(q) = h'(q)$ for a simple function
$q\in [0,1]^{X}$.  Recall that such a simple $q$ can be written as $q
= \ovee_{i}\, r_{i}\charac{X_i}$, like
in~\eqref{SimpleFunNormalFormEqn}, where the (disjoint) subsets
$X_{i}\subseteq X$ cover $X$.  Since $h,h'\in\Exp(X)$ are maps
of effect modules we get:
$$\begin{array}{rcccccccccl}
h(q)
& = &
\sum_{i}r_{i}h(\textbf{1}_{X_i})
& = &
\sum_{i}r_{i}\Phi(h)(X_{i})
& = &
\sum_{i}r_{i}\Phi(h')(X_{i})
& = &
\sum_{i}r_{i}h'(\textbf{1}_{X_i})
& = &
h'(q).
\end{array}$$

\noindent For an arbitrary $p\in [0,1]^{X}$ we first write $p =
\lim_{n}p_{n}$ as limit of simple functions $p_n$ like in
Lemma~\ref{SimplePredLem}. Lemma~\ref{AEModNonExpLem} implies
that $h,h'$ are continuous, and so we get $h=h'$ from:
$$\begin{array}{rcccccl}
h(p)
& = &
\lim_{n}h(p_{n})
& = &
\lim_{n}h'(p_{n})
& = &
h'(p).
\end{array}$$

For surjectivity of $\Phi$ assume a finitely additive measure $m\colon
\Pow(X) \rightarrow [0,1]$. We need to define a function $h\in\Exp(X)$
with $\Phi(h) = m$. We define such a $h$ first on a simple function $q
= \ovee_{i}\, r_{i}\charac{X_i}$ as $h(q) =
\sum_{i}r_{i}m(X_{i})$. For an arbitrary $p\in [0,1]^{X}$, written as
$p = \lim_{n}p_{n}$, like in Lemma~\ref{SimplePredLem}, we define
$h(p) = \lim_{n}h(p_{n})$. Then $\Phi(h) = m$, since for $U\subseteq
X$ we have:
$$\begin{array}{rcccl}
\Phi(h)(U)
& = &
h(\charac{U})
& = &
m(U).
\end{array}\eqno{\QEDbox}$$
\end{myproof}

The inverse $h = \Phi^{-1}(m)$ that is constructed in this proof may
be understood as an integral $h(p) = \int pdm$. The precise nature of
the bijection $\Phi$ remains unclear at this stage since we have not
yet identified the (algebraic) structure of the sets $\Exp(X)$. But
via this bijection we can understand mapping a set to its finitely
additive measures, \textit{i.e.}~$X \mapsto \EA(\Pow(X),[0,1])$, as a
monad.

Yet another perspective is useful in this context. The characteristic
function mapping:
$$\xymatrix{
[0,1] \times \Pow(X) \ar[r] & [0,1]^{X}
\qquad\mbox{given by}\qquad (r,U)\ar@{|->}[r] & r\cdot \charac{U}
}$$

\noindent is a bihomomorphism of effect modules. Hence it gives rise
to a map of effect modules $[0,1]\otimes \Pow(X) \rightarrow
[0,1]^{X}$, where the tensor product $[0,1]\otimes\Pow(X)$ forms a
more abstract description of the effect module of simple (step) functions
$[X\rightarrow_{s}[0,1]]$ from Lemma~\ref{SimplePredLem} (see
also~\cite[Thm.~5.6]{Gudder98}). Lemma~\ref{SimplePredLem} says that
this map is dense. This gives a quick proof of
Theorem~\ref{FinAddMeasThm}:
$$\begin{array}{rcl}
\Exp(X)
& = &
\EMod\big([0,1]^{X},\,[0,1]\big) \\
& \cong &
\EMod\big([0,1]\otimes\Pow(X),\,[0,1]\big)
  \qquad\mbox{by denseness} \\
& \cong &
\EA\big(\Pow(X), [0,1]\big).
\end{array}$$

\noindent This last isomorphism is standard, because
$[0,1]\otimes\Pow(X)$ is the free effect module on $\Pow(X)$.

\section{The expectation and ultrafilter monads}\label{ExpUFSec}

In this section we show that the sets $\Exp(X)$ carry a compact
Hausdorff structure and we identify its topology. The unit interval
$[0,1]$ plays an important role. It is a compact Hausdorff space,
which means that it carries an algebra of the ultrafilter monad, see
Subsection~\ref{UFSubsec}. We shall write this algebra as $\ch =
\ch_{[0,1]}\colon\UF([0,1])\rightarrow[0,1]$. What this map precisely
does is described in Example~\ref{UnitIntervalUFEx}; but mostly we use
it abstractly, as an $\UF$-algebra.  The technique we use to define
the following map of monads is copied from
Lemma~\ref{ContAlgMndMapLem}.

\begin{proposition}
\label{UFExpNatroProp}
There is a map of monads $\tau \colon \UF \Longrightarrow \Exp$, given
on an ultrafilter $\calF\in\UF(X)$ and $p\in[0,1]^{X}$ by:
$$\begin{array}{rcl}
\tau_{X}(\calF)(p)
& = &
\ch\big(\UF(p)(\calF)\big) \\
& = &
\inf\setin{s}{[0,1]}{[0,s]\in \UF(p)(\calF)}
   \qquad\qquad\mbox{by~\eqref{UnitIntervalUFEqn}} \\
& = &
\inf\setin{s}{[0,1]}{\setin{x}{X}{p(x)\leq s}\in\calF}.
\end{array}$$

\noindent In this description the functor $\UF$ is applied to $p$, as
function $X\rightarrow [0,1]$, giving $\UF(p) \colon \UF(X)
\rightarrow \UF([0,1])$.
\end{proposition}

\begin{myproof}
We first have to check that $\tau$ is well-defined, \textit{i.e.}~that
$\tau_{X}(\calF) \colon [0,1]^{X} \rightarrow [0,1]$ is a
morphism of effect modules.
\begin{itemize}
\item Preservation $\tau_{X}(\calF)(r\cdot p) = r\cdot p
  \tau_{X}(\calF)$ of multiplication with scalar $r\in[0,1]$. This
  follows by observing that multiplication $r\cdot (-)\colon [0,1]
  \rightarrow [0,1]$ is a continuous function, and thus a morphism of
  algebras in the square below.
$$\xymatrix@R-.5pc{
\UF([0,1])\ar[rr]^-{\UF(r\cdot (-))}\ar[d]_{\ch} & & 
   \UF([0,1])\ar[d]^{\ch} \\
[0,1]\ar[rr]_-{r\cdot(-)} & & [0,1]
}$$

\noindent Thus:
$$\begin{array}{rcl}
\tau(\calF)(r\cdot p)
& = &
\big(\ch \after \UF(r\cdot(-) \after p)\big)(\calF) \\
& = &
\big(r\cdot(-) \after \ch \after \UF(p)\big)(\calF)
\hspace*{\arraycolsep} = \hspace*{\arraycolsep}
r\cdot \tau(\calF)(p).
\end{array}$$

\item Preservation of $\ovee$, is obtained in the same manner, using
  that addition $+\colon [0,1]\times[0,1] \rightarrow [0,1]$ is
  continuous.

\auxproof{
We use the product algebra on $[0,1]\times[0,1]$ on the left in:
$$\xymatrix@R-.5pc{
\UF([0,1]\times[0,1])\ar[rr]^-{\UF(+)}
   \ar[d]_{\tuple{\ch\after\UF(\pi_{1}), \ch\after\UF(\pi_{2})}} & & 
   \UF([0,1])\ar[d]^{\ch} \\
[0,1]\times[0,1]\ar[rr]_-{+} & & [0,1]
}$$

\noindent For $p,q\in[0,1]^{X}$ with $p(x) + q(x)
\leq 1$, for all $x\in X$, we have:
$$\begin{array}{rcl}
\tau_{X}(\calF)(p\ovee q)
& = &
\big(\ch \after \UF(+ \after \tuple{p,q})\big)(\calF) \\
& = &
\big(+ \after \tuple{\ch\after\UF(\pi_{1}), \ch\after\UF(\pi_{2})} \after
   \UF(\tuple{p,q})\big)(\calF) \\
& = &
\big(+ \after \tuple{\ch\after\UF(p), \ch\after\UF(q)}\big)(\calF) \\
& = &
\ch\big(\UF(p)(\calF)\big) + \ch\big(\UF(q)(\calF)\big) \\
& = &
\tau_{X}(\calF)(p) \ovee \tau_{X}(\calF)(q).
\end{array}$$
}

\item Constant functions $\lam{x}{a} \in [0,1]^X$, including
$0$ and $1$, are preserved:
$$\begin{array}{rcl}
\tau_{X}(\calF)(\lam{x}{a})
& = &
\ch\big(\UF(\lam{x}{a})(\calF)\big) \\
& = &
\ch\big(\setin{U}{\Pow([0,1])}{(\lam{x}{a})^{-1}(U) \in \calF}\big) \\
& = &
\ch\big(\setin{U}{\Pow([0,1])}{\setin{x}{X}{a\in U} \in \calF}\big) \\
& = &
\ch\big(\setin{U}{\Pow([0,1])}{a \in U}\big) 
   \qquad \mbox{since }\emptyset \not\in \calF \\
& = &
\ch(\eta(a)) \\
& = &
a.
\end{array}$$
\end{itemize}

\noindent We leave naturality of $\tau$ and commutation with units to
the reader and check that $\tau$ commutes with multiplications
$\mu^{\Exp}$ and $\mu^{\UF}$ of the expectation and ultrafilter
monads. Thus, for $\mathcal{A}\in\UF^{2}(X)$ and $p\in[0,1]^{X}$,
we calculate:
$$\begin{array}[b]{rcl}
\big(\mu^{\Exp} \after \tau \after \UF(\tau)\big)(\mathcal{A})(p)
& = &
\mu\Big(\tau\big(\UF(\tau)(\mathcal{A})\big)\big)\Big)(p) \\
& = &
\tau\big(\UF(\tau)(\mathcal{A})\big)\big(\lam{k}{k(p)}\big) \\
& = &
\ch\Big(\UF(\lam{k}{k(p)})\big(\UF(\tau)(\mathcal{A})\big)\Big) \\
& = &
\ch\Big(\UF(\lam{\calF}{\tau(\calF)(p)})(\mathcal{A})\big)\Big) \\
& = &
\ch\Big(\UF(\lam{\calF}{\ch(\UF(p)(\calF))})(\mathcal{A})\big)\Big) \\
& = &
\ch\Big(\UF(\ch \after \UF(p))(\mathcal{A})\big)\Big) \\
& = &
\big(\ch \after \UF(\ch \after \UF(p))\big)(\mathcal{A}) \\
& = &
\big(\ch \after \UF(\ch) \after \UF^{2}(p)\big)(\mathcal{A}) \\
& = &
\big(\ch \after \mu^{\UF} \after \UF^{2}(p)\big)(\mathcal{A}) \\
& = &
\big(\ch \after \UF(p) \after \mu^{\UF}\big)(\mathcal{A}) \\
& = &
\ch\big(\UF(p)(\mu^{\UF}(\mathcal{A}))\big) \\
& = &
\big(\tau \after \mu^{\UF}\big)(\mathcal{A})(p).
\end{array}\eqno{\QEDbox}$$

\auxproof{
For a function $f\colon X\rightarrow Y$, and for $\calF\in\UF(X)$,
$$\begin{array}{rcl}
\big(\Exp(f) \after \tau_{X}\big)(\calF)
& = &
\Exp(f)\big(\tau_{X}(\calF)\big) \\
& = &
\lamin{p}{[0,1]^{Y}}{\tau_{X}(\calF)(p \after f)} \\
& = &
\lamin{p}{[0,1]^{Y}}{\ch\big(\UF(p \after f)(\calF)\big)} \\
& = &
\lamin{p}{[0,1]^{Y}}{\ch\big(\UF(p)(\UF(f)(\calF))\big)} \\
& = &
\tau_{Y}\big(\UF(f)(\calF)\big) \\
& = &
\big(\tau_{Y} \after \UF(f)\big)(\calF).
\end{array}$$

\noindent Commutation with units:
$$\begin{array}{rcl}
\big(\tau_{X} \after \eta_{X}\big)(x)
& = &
\tau_{X}\big(\eta_{X}(x)\big) \\
& = &
\lamin{p}{[0,1]^{X}}{\ch\big(\UF(p)(\eta_{X}(x))\big)} \\
& = &
\lamin{p}{[0,1]^{X}}{\ch\big(\eta_{[0,1]}(p(x))\big)} \\
& = &
\lamin{p}{[0,1]^{X}}{p(x)} \\
& = &
\eta_{X}(x).
\end{array}$$
}
\end{myproof}

\begin{corollary}
\label{ExpAlgCHCor}
There is a functor $\Alg(\Exp) \rightarrow \Alg(\UF) = \CH$, by
pre-composition: $\big(\Exp(X)\xrightarrow{\alpha} X\big) \longmapsto
\big(\UF(X) \xrightarrow{\alpha\after\tau} X\big)$. This functor has a
left adjoint by Lemma~\ref{MndMapAlgLem}.

In particular, the underlying set $X$ of each $\Exp$-algebra
$\alpha\colon\Exp(X)\rightarrow X$ carries a compact Hausdorff
topology, with $U\subseteq X$ closed iff for each $\calF\in\UF(X)$
with $U\in\calF$ one has $\alpha(\tau(\calF))\in U$, as described
in Subsection~\ref{UFSubsec}. \QED
\end{corollary}

With respect to this topology on $\Exp(X)$, several maps are continuous.

\begin{lemma}
\label{ExpAlgCHContinuousLem}
The following maps are continuous functions.
$$\xymatrix@C+0.2pc{
\UF(X)\ar[r]^-{\tau_X} & \Exp(X)
\qquad
\Exp(X)\ar[r]^-{\alpha}_-{\text{algebra}} & X
\qquad
\Exp(X)\ar[r]^-{\Exp(f)} & \Exp(Y)
\qquad
\Exp(X)\ar[r]^-{\ev_{p}\,= }_-{\lam{h}{h(p)}} & [0,1].
}$$
\end{lemma}

\begin{myproof}
One shows that these maps are morphisms of $\UF$-algebras. For
instance, $\tau_X$ is continuous because it is a map of monads:
commutation with multiplications, as required in~\eqref{MndMapEqn},
precisely says that it is a map of algebras, in the square on the
left below.
$$\xymatrix@R-.5pc@C-.5pc{
\UF^{2}(X)\ar[d]_{\mu_{X}}\ar[rr]^-{\UF(\tau_{X})} & &
   \UF(\Exp(X))\ar[d]^{\mu_{X}\after\tau_{\Exp(X)}}
& \qquad\qquad & 
\UF(\Exp(X))\ar[d]_{\mu_{X}\after\tau_{\Exp(X)}}\ar[rr]^-{\UF(\alpha)} & &
   \UF(X)\ar[d]^{\alpha\after\tau_{X}} \\
\UF(X)\ar[rr]_-{\tau_{X}} & & \Exp(X)
& &
\Exp(X)\ar[rr]_-{\alpha} & & X
}$$

\noindent The rectangle on the right expresses that an Eilenberg-Moore
algebra $\alpha\colon \Exp(X) \rightarrow X$ is a continuous function.
It commutes by naturality of $\tau$:
$$\begin{array}{rcccl}
\alpha \after \tau_{X} \after \UF(\alpha)
& = &
\alpha \after \Exp(\alpha) \after \tau_{\Exp(X)} 
& = &
\alpha \after \mu_{X} \after \tau_{\Exp(X)}.
\end{array}$$

\noindent For $f\colon X\rightarrow Y$, continuity of $\Exp(f)\colon
\Exp(X) \rightarrow \Exp(Y)$ follows directly from naturality of
$\tau$.  Finally, for $p\in[0,1]^{X}$ the map $\ev_{p} =
\lam{h}{h(p)} \colon \Exp(X) \rightarrow [0,1]$ is continuous because
for $\calF\in\UF(\Exp(X))$,
$$\begin{array}[b]{rcl}
\big(\ev_{p} \after \mu_{X} \after \tau_{\Exp(X)}\big)(\calF)
& = &
\mu_{X}\big(\tau_{\Exp(X)}(\calF)\big)(p) \\
& = &
\tau_{\Exp(X)}(\calF)(\lam{k}{k(p)}) \\
& = &
\tau_{\Exp(X)}(\calF)(\ev_{p}) \\
& = &
\ch\big(\UF(\ev_{p})(\calF)\big) \\
& = &
\big(\ch \after \UF(\ev_{p})\big)(\calF).
\end{array}\eqno{\QEDbox}$$
\end{myproof}

The next step is to give a concrete description of this compact
Hausdorff topology on sets $\Exp(X)$, as induced by the algebra
$\UF(\Exp(X))\rightarrow\Exp(X)$.

\begin{proposition}
\label{ExpTopologyProp}
Fix a set $X$. For a predicate $p\in [0,1]^{X}$ and a rational number
$s\in [0,1]\cap\mathbb{Q}$ write:
$$\begin{array}{rcl}
\square_{s}(p)
& = &
\setin{h}{\Exp(X)}{h(p) > s}.
\end{array}$$

\noindent These sets $\square_{s}(p) \subseteq \Exp(X)$ form a
subbasis for the topology on $\Exp(X)$.
\end{proposition}

\begin{myproof}
We reason as follows. The subsets $\square_{s}(p)$ are open in the
compact Hausdorff topology induced on $\Exp(X)$ by the algebra
structure $\UF(\Exp(X))\rightarrow \Exp(X)$. They form a subbasis for
a Hausdorff topology on $\Exp(X)$. Hence by Lemma~\ref{UniqueCHLem}
this topology is the induced one. We now elaborate these steps.

The Eilenberg-Moore algebra $\UF(\Exp(X)) \rightarrow \Exp(X)$ is
given by $\mu_{X} \after \tau_{\Exp(X)}$. Hence the associated closed
sets $U\subseteq\Exp(X)$ are those satisfying $U\in\calF \Rightarrow
\mu_{X}(\tau_{\Exp(X)}(\calF)) \in U$, for each
$\calF\in\UF(\Exp(X))$, see Subsection~\ref{UFSubsec}.  We wish to
show that $\neg\square_{s}(p) = \set{h}{h(p)\leq s}\subseteq \Exp(X)$
is closed. We reason backwards, starting with the required conclusion.
$$\begin{array}{rcl}
\mu(\tau(\calF))\in\neg\square_{s}(p)
& \Longleftrightarrow &
\mu(\tau(\calF))(p) \leq s \\
& \Longleftrightarrow &
\ch\big(\UF(\lam{k}{k(p)})(\calF)\big) \in [0,s] \\
& & \quad \mbox{since }
   \mu(\tau(\calF))(p) = \tau(\calF)(\lam{k}{k(p)}) = 
   \ch\big(\UF(\lam{k}{k(p)})(\calF)\big) \\
& \Longleftarrow &
[0,s] \in \UF(\lam{k}{k(p)})(\calF) \\
& & \quad \mbox{since  $[0,s]\subseteq[0,1]$ is closed} \\
& \Longleftrightarrow &
(\lam{k}{k(p)})^{-1}([0,s]) \in \calF \\
& \Longleftrightarrow &
\setin{h}{\Exp(X)}{h(p)\in[0,s]} = \neg\square_{s}(p)\in \calF.
\end{array}$$

\noindent Hence $\neg\square_{s}(p) \subseteq \Exp(X)$ is closed,
making $\square_{s}(p)$ open.

Next we need to show that these $\square_{s}(p)$'s give rise to a
Hausdorff topology. So assume $h\neq h'\in\Exp(X)$. Then there must be
a $p\in [0,1]^{X}$ with $h(p) \neq h'(p)$. Without loss of generality
we assume $h(p) < h'(p)$. Find an $s\in [0,1]\cap\mathbb{Q}$ with
$h(p) < s < h'(p)$. Then $h'\in\square_{s}(p)$. Also:
$$\begin{array}{rcccccccl}
h(p^{\perp})
& = &
1 - h(p)
& > &
1-s
& > &
1 - h'(p)
& = &
h'(p^{\perp}).
\end{array}$$

\noindent Hence $h\in\square_{1-s}(p^{\perp})$. These sets
$\square_{s}(p)$ and $\square_{1-s}(p^{\perp})$ are disjoint, since:
$k\in\square_{s}(p) \cap \square_{1-s}(p^{\perp})$ iff both $k(p) > s$
and $1-k(p) > 1-s$, which is impossible. \QED
\end{myproof}

As is well-known, ultrafilters on a set $X$ can also be understood as
finitely additive measures $\Pow(X) \rightarrow \{0,1\}$. Using
Theorem~\ref{FinAddMeasThm} we can express more precisely how the
expectation monad $\Exp$ is a probabilistic version of the ultrafilter
monad $\UF$, namely via the descriptions:
$$\begin{array}{rclcrcl}
\Exp(X)
& \cong &
\EA\big(\Pow(X), [0,1]\big)
& \quad\mbox{and}\quad &
\UF(X)
& \cong &
\EA\big(\Pow(X), \{0,1\}\big). 
\end{array}$$

\noindent We have $\EA\big(\Pow(X), \{0,1\}\big) = \BA\big(\Pow(X),
\{0,1\}\big)$ because in general, for Boolean algebras $B,B'$ a
homomorphism of Boolean algebras $B\rightarrow B'$ is the same as an
effect algebra homomorphism $B\rightarrow B'$.

\begin{lemma}
\label{UFExpNatroMonLem}
The components $\tau_{X}\colon \UF(X) \rightarrow \Exp(X)$ are injections.
\end{lemma}

\begin{myproof}
Because there are isomorphisms:
$$\vcenter{\xymatrix@R-1.55pc{
\UF(X)\ar@{=}[d]_{\wr}\ar[rr]^-{\tau_X} & & \Exp(X)\ar@{=}[d]_{\wr} \\
\EA\big(\Pow(X), \{0,1\}\big)\ar@{ >->}[rr] & &
   \EA\big(\Pow(X), [0,1]\big)
}}\eqno{\QEDbox}$$
\end{myproof}

\section{The expectation and distribution monads}\label{ExpDstSec}

We continue with the implications of
Lemma~\ref{ComposableAdjunctionLem} in the current situation,
especially with the natural
transformation~\eqref{ComposableAdjunctionNatroDiag}. This leads
to convex structure on sets $\Exp(X)$.

\begin{lemma}
\label{DstExpNatroLem}
There is a map of monads:
\begin{equation}
\label{DstExpNatroEqn}
\sigma\colon \Dst \Longrightarrow \Exp
\quad\mbox{given by}\quad
\begin{array}{rcl}
\sigma_{X}(\varphi)
& = &
\lamin{p}{[0,1]^{X}}{\sum_{x}\varphi(x)\cdot p(x)},
\end{array}
\end{equation}

\noindent where the dot $\cdot$ describes multiplication in $[0,1]$.


All components $\sigma_{X} \colon \Dst(X) \rightarrow \Exp(X)$ are
injections. And for \emph{finite} sets $X$ the component at $X$ is an
isomorphism $\Dst(X) \conglongrightarrow \Exp(X)$.
\end{lemma}

With this result we have completed the positioning of the expectation
monad in Diagram~(\ref{MonadsOverviewDiag}), in between the
distribution and ultrafilter monad on the hand, and the continuation
monad on the other.

\begin{myproof}
By construction via~\eqref{ComposableAdjunctionNatroDiag} the natural
transformation $\sigma\colon \Dst \Rightarrow \Exp$ is a map of
monads. Next, assume $X$ is finite, say $X = \{x_{1}, \ldots,
x_{n}\}$.  Each $p\in [0,1]^{X}$ is determined by the values
$p(x_{i})\in[0,1]$.  Using the effect module structure of $[0,1]^{X}$,
this $p$ can be written as sum of scalar multiplications:
$$\begin{array}{rcl}
p
& = &
p(x_{1})\cdot \charac{x_{1}} \ovee \cdots \ovee \,
   p(x_{n})\cdot \charac{x_{n}},
\end{array}$$

\noindent where $\charac{x_{i}}\colon X \rightarrow [0,1]$ is the
characteristic function of the singleton $\{x_{i}\}\subseteq X$. A map
of effect modules $h\in\Exp(X) = \EMod([0,1]^{X},[0,1])$ will thus
send such a predicate $p$ to:
$$\begin{array}{rcl}
h(p)
& = &
h\big(p(x_{1})\cdot \charac{x_{1}} \ovee \cdots \ovee \,
   p(x_{n})\cdot \charac{x_{n}}\big) \\
& = &
p(x_{1})\cdot h(\charac{x_{1}}) + \cdots + 
   p(x_{n})\cdot h(\charac{x_{n}}), \\
\end{array}$$

\noindent since $\ovee$ is $+$ in $[0,1]$. Hence $h$ is completely
determined by these values $h(\charac{x_{i}})\in [0,1]$. But since
$\ovee_{i}\,\charac{x_{i}} = 1$ in $[0,1]^{X}$ we also have
$\sum_{i}h(\charac{x_{i}}) = 1$. Hence $h$ can be described by the
convex sum $\varphi\in\Dst(X)$ given by $\varphi(x) =
h(\charac{x})$. Thus we have a bijection $\Exp(X) \cong \Dst(X)$. In
fact, $\sigma_{X}$ describes (the inverse of) this bijection, since:
$$\begin{array}[b]{rcl}
\sigma_{X}(\varphi)(p)
& = &
\sum_{i}\varphi(x_{i})\cdot p(x_{i}) \\
& = &
\sum_{i}p(x_{i}) \cdot h(\charac{x_{i}}) \\
& = &
h\big(\ovee_{i} p(x_{i})\cdot \charac{x_{i}}\big) \\
& = &
h(p).
\end{array}\eqno{\QEDbox}$$

\auxproof{ 
First we check that $\sigma$ results from
~\eqref{ComposableAdjunctionNatroDiag}, via the isomorphism
$\Conv(\Dst(X),[0,1]) \cong [0,1]^{X}$.  The
formula~\eqref{ComposableAdjunctionNatroDiag} gives in the present
situation:
$$\begin{array}{rcl}
\varphi\in\Dst(X)
& \longmapsto &
\eta(\varphi) = \lam{f}{f(\varphi)}
   \in \EMod\big(\Conv(\Dst(X),[0,1]), [0,1]\big) \\
& \longmapsto &
\lam{p}{\big(\lam{f}{f(\varphi)}\big)(
   \lam{\psi}{\sum_{x}\psi(x)\cdot p(x)})}
   \in \EMod([0,1]^{X},[0,1]) \\
& = &
\lamin{p}{[0,1]^{X}}{\sum_{x}\varphi(x)\cdot p(x)} \\
& = &
\sigma_{X}(\varphi).
\end{array}$$

Clearly, this $\sigma_{X}(\varphi)$ is a map of effect modules. 
We explicitly check the requirements of monad maps:
$$\begin{array}{rcl}
\big(\Exp(f) \after \sigma_{X}\big)(\varphi)
& = &
\Exp(f)\big(\sigma_{X}(\varphi)\big) \\
& = &
\lamin{p}{[0,1]^{Y}}{\sigma_{X}(\varphi)(p \after f)} \\
& = &
\lamin{p}{[0,1]^{Y}}{\sum_{x}\varphi(x)\cdot p(f(x))} \\
& = &
\lamin{p}{[0,1]^{Y}}{\sum_{y}\sum_{x\in f^{-1}(y)}\varphi(x)\cdot p(y)} \\
& = &
\lamin{p}{[0,1]^{Y}}{\sum_{y}\Dst(f)(\varphi)(y) \cdot p(y)} \\
& = &
\sigma_{Y}\big(\Dst(f)(\varphi)\big) \\
& = &
\big(\sigma_{Y} \after \Dst(f)\big)(\varphi) \\
\big(\sigma_{X} \after \eta_{X}\big)(x)
& = &
\sigma_{X}(1x) \\
& = &
\lamin{p}{[0,1]^{X}}{\sum_{x'}(1x)x'\cdot p(x')} \\
& = &
\lamin{p}{[0,1]^{X}}{p(x)} \\
& = &
\eta_{X}(x) \\
\big(\mu_{X} \after \sigma_{\Exp(X)} \after \Dst(\sigma_{X})\big)
   (\sum_{i}r_{i}\varphi_{i})
& = &
\mu_{X}\Big(\sigma_{\Exp(X)}\big(\Dst(\sigma_{X})
   (\sum_{i}r_{i}\varphi_{i})\big)\Big) \\
& = &
\lamin{p}{[0,1]^{X}}{\sigma_{\Exp(X)}\big(\sum_{i}r_{i}\sigma_{X}(\varphi_{i})\big)
   \big(\lamin{k}{\Exp(X)}{k(p)}\big)} \\
& = &
\lamin{p}{[0,1]^{X}}{\sum_{i}r_{i}\sigma_{X}(\varphi_{i})(p)} \\
& = &
\lamin{p}{[0,1]^{X}}{\sum_{i}r_{i}\cdot(\sum_{x}\varphi_{i}(x)\cdot p(x))} \\
& = &
\lamin{p}{[0,1]^{X}}{\sum_{x}\big(\sum_{i}r_{i}\cdot\varphi_{i}(x)\big)\cdot p(x)} \\
& = &
\lamin{p}{[0,1]^{X}}{\sum_{x}\mu_{X}(\sum_{i}r_{i}\varphi_{i})(x)\cdot p(x)} \\
& = &
\sigma_{X}(\mu_{X}(\sum_{i}r_{i}\varphi_{i})) \\
& = &
\big(\sigma_{X} \after \mu_{X}\big)(\sum_{i}r_{i}\varphi_{i}).
\end{array}$$
}
\end{myproof}

\begin{corollary}
\label{ExpAlgConvCor}
There is a functor $\Alg(\Exp) \rightarrow \Alg(\Dst) = \Conv$, by
pre-composition: $\big(\Exp(X)\xrightarrow{\alpha} X\big) \longmapsto
\big(\Dst(X) \xrightarrow{\alpha\after\sigma} X\big)$. It has a
left adjoint by Lemma~\ref{MndMapAlgLem}. \QED
\end{corollary}

Explicitly, for each $\Exp$-algebra $\alpha\colon\Exp(X)\rightarrow
X$, the set $X$ is a convex set, with sum of a formal convex
combination $\sum_{i}r_{i}x_{i}$ given by the element:
$$\begin{array}{rcccl}
\alpha\big(\sigma_{X}(\sum_{i}r_{i}x_{i})\big)
& = &
\alpha\big(\lamin{p}{[0,1]^{X}}{\sum_{i}r_{i}\cdot p(x_{i})}\big) 
& \in &
X.
\end{array}$$

\noindent Lemma~\ref{DstExpNatroLem} implies that if the carrier $X$
is finite, the algebra structure $\alpha$ corresponds precisely to
such convex structure on $X$. If $X$ is non-finite we still have to
find out what $\alpha$ involves.

Here is another (easy) consequence of Lemma~\ref{DstExpNatroLem}.

\begin{corollary}
\label{Exp012Cor}
On the first few finite sets: empty $0$, singleton $1$, and
two-element $2$ one has:
$$\begin{array}{rclcrclcrcl}
\Exp(0)
& \cong &
0
& \qquad &
\Exp(1)
& \cong & 
1
& \qquad &
\Exp(2)
& \cong &
[0,1].
\end{array}$$

\noindent The isomorphism in the middle says that $\Exp$ is an
\emph{affine} functor.

\end{corollary}

\begin{myproof}
The isomorphisms follow easily from $\Exp(X)\cong \Dst(X)$ for finite
$X$.  
\QED
\end{myproof}

\begin{remark}
\label{DstExpNatroRem}
(1)~The natural transformation $\sigma\colon\Dst\Rightarrow\Exp$
from~\eqref{DstExpNatroEqn} implicitly uses that the unit interval
$[0,1]$ is convex. This can be made explicit in the following way.
Describe this convexity via an algebra $\cv\colon \Dst([0,1])
\rightarrow [0,1]$. Then we can equivalently describe $\sigma$ as:
$$\begin{array}{rcl}
\sigma_{X}(\varphi)(p)
& = &
\cv\big(\Dst(p)(\varphi)\big).
\end{array}$$

\noindent This alternative description is similar to the construction
in Proposition~\ref{UFExpNatroProp}, for a natural transformation
$\UF\Rightarrow\Exp$ (see also Lemma~\ref{ContAlgMndMapLem}).


(2)~From Corollaries~\ref{ExpAlgCHCor} and~\ref{ExpAlgConvCor} we
  know that the sets $\Exp(X)$ are both compact Hausdorff and
  convex. This means that we can take free extensions of the maps
  $\tau\colon \UF(X)\rightarrow \Exp(X)$ and $\sigma\colon
  \Dst(X)\rightarrow \Exp(X)$, giving maps $\Dst(\UF(X))\rightarrow
  \Exp(X)$ and $\UF(\Dst(X))\rightarrow \Exp(X)$, \textit{etc.} The
  latter map is the composite:
$$\xymatrix{
\UF(\Dst(X))\ar[r]^-{\UF(\sigma)} &
  \UF(\Exp(X))\ar[r]^-{\tau} & \Exp^{2}(X)\ar[r]^-{\mu} & \Exp(X).
}$$

\noindent Using Example~\ref{UnitIntervalUFEx}, it can be described
more concretely on $\calF\in\UF(\Dst(X))$ and $p\in[0,1]^{X}$ as:
$$\textstyle\inf\setin{s}{[0,1]}{\,\setin{\varphi}{\Dst(X)}
   {\sum_{x}\varphi(x)\cdot p(x) \leq s}\in\calF}.$$

\auxproof{
$$\begin{array}{rcl}
\lefteqn{\big(\mu \after \tau \after \UF(\sigma)\big)(\calF)(p)} \\
& = &
\mu\big(\tau(\UF(\sigma)(\calF))\big)(p) \\
& = &
\tau(\UF(\sigma)(\calF))(\lam{k}{k(p)}) \\
& = &
\ch\big(\UF(\lam{k}{k(p)})(\UF(\sigma)(\calF))\big) \\
& = &
\ch\big(\UF(\lam{\varphi}{\sigma(\varphi)(p)})(\calF)\big) \\
& = &
\inf\set{s}{[0,s]\in\UF(\lam{\varphi}{\sigma(\varphi)(p)})(\calF)} 
   \qquad \mbox{by Example~\ref{UnitIntervalUFEx}} \\
& = &
\inf\set{s}{(\lam{\varphi}{\sigma(\varphi)(p)})^{-1}([0,s])\in \calF} \\
& = &
\inf\set{s}{\set{\varphi}{\sigma(\varphi)(p)\in[0,s]}\in \calF} \\
& = &
\inf\set{s}{\set{\varphi}{\sum_{x}\varphi(x)\cdot p(x)\leq s}\in \calF}.
\end{array}$$
}



\end{remark}

\auxproof{
For later usage we record the following preservation result for the
expectation functor $\Exp$.  Other preservation properties require
further investigation, for instance, whether $\Exp$ preserves weak
pullbacks (like the distribution functor $\Dst$,
see~\cite{Moss99,VinkR99}).

\begin{lemma}
\label{ExpMonoPresLem}
The functor $\Exp\colon\Sets \rightarrow\Sets$ preserves injections.
\end{lemma}

\begin{myproof}
Assume an injection $f\colon X\rightarrowtail Y$ in $\Sets$. If $X=0$,
we obtain an injection $\Exp(0) \cong 0 \rightarrowtail \Exp(Y)$ by
Corollary~\ref{Exp012Cor}. So we may assume $X$ is non-empty, and thus
a function $s\colon Y\rightarrow X$ with $s\after f = \idmap[X]$.
This makes $f$ a split mono. Hence $\Exp(f)$ is also a split
mono. \QED

\auxproof{
Explicitly, let $h,h'\in\Exp(X)$ satisfy $\Exp(f)(h) = \Exp(f)(h')$. In order
to prove $h=h'$, assume $p\in [0,1]^{X}$. Then $q = p \after s \in
[0,1]^{Y}$ yields:
$$\begin{array}[b]{rcl}
h(p)
\hspace*{\arraycolsep} = \hspace*{\arraycolsep}
h(q \after f)
& = &
\Exp(f)(h)(q) \\
& = &
\Exp(f)(h')(q) 
\hspace*{\arraycolsep} = \hspace*{\arraycolsep}
h'(q \after f)
\hspace*{\arraycolsep} = \hspace*{\arraycolsep}
h'(p).
\end{array}\eqno{\QEDbox}$$
}
\end{myproof}
}

The next result is the affine analogue of
Lemma~\ref{ExpAlgCHContinuousLem}.

\begin{lemma}
\label{ExpAlgDstAffineLem}
The following maps are affine functions.
$$\xymatrix@C+0.2pc{
\Dst(X)\ar@{ >->}[r]^-{\sigma_X} & \Exp(X)
\qquad
\Exp(X)\ar[r]^-{\alpha}_-{\text{algebra}} & X
\qquad
\Exp(X)\ar[r]^-{\Exp(f)} & \Exp(Y)
\qquad
\Exp(X)\ar[r]^-{\ev_{p}\,= }_-{\lam{h}{h(p)}} & [0,1].
}$$
\end{lemma}

\begin{myproof}
Verifications are done like in the proof of
Lemma~\ref{ExpAlgCHContinuousLem}. We only do the last one. We need
to prove that the following diagram commutes,
$$\xymatrix@R-.5pc{
\Dst(\Exp(X))\ar[d]_{\mu_{X}\after\sigma_{X}}\ar[rr]^-{\Dst(\ev_{p})} & &
   \Dst([0,1])\ar[d]^{\cv} \\
\Exp(X)\ar[rr]_{\ev_p} & & [0,1]
}$$

\noindent where the algebra \cv interprets formal convex combinations
as actual combinations. For a distribution $\Phi =
\sum_{i}r_{i}h_{i}\in\Dst(\Exp(X))$ we have:
$$\begin{array}[b]{rcl}
\big(\ev_{p} \after \mu \after \sigma\big)(\Phi)
& = &
\mu\big(\sigma(\Phi)\big)(p) \\
& = &
\sigma(\Phi)(\ev_{p}) \\
& = &
\sum_{i}r_{i}\cdot \ev_{p}(h_{i}) \\
& = &
\cv\big(\sum_{i}r_{i}\ev_{p}(h_{i}) \\
& = &
\cv\big(\Dst(\ev_{p})(\sum_{i}r_{i}h_{i})\big) \\
& = &
\big(\cv \after \Dst(\ev_{p})\big)(\Phi).
\end{array}\eqno{\QEDbox}$$

\auxproof{
This follows from---or is actually equivalent to---the requirement that
monad maps commute with multiplication, in:
$$\xymatrix@R-.5pc{
\Dst^{2}(X)\ar[d]_{\mu}\ar[rr]^-{\Dst(\sigma_{X})} & &
   \Dst(\Exp(X))\ar[d]^{\mu\after\sigma} \\
\Dst(X)\ar[rr]_-{\sigma_{X}} & & \Exp(X)
}$$
}

\end{myproof}

The $\Dst$-algebras obtained from $\Exp$-algebras turn out to
be continuous functions. This connects the convex and topological
structures in such algebras.

\begin{lemma}
\label{ExpAlgDstContinuousLem}
The maps $\sigma_{X}\colon \Dst(X) \rightarrowtail \Exp(X)$ are
(trivially) continuous when we provide $\Dst(X)$ with the subspace
topology with basic opens $\square_{s}(p) \subseteq \Dst(X)$ given by
restriction: $\square_{s}(p) =
\setin{\varphi}{\Dst(X)}{\sum_{x}\varphi(x) \cdot p(x) > s}$, for
$p\in [0,1]^{X}$ and $s\in[0,1]\cap\mathbb{Q}$.

For each $\Exp$-algebra $\alpha\colon\Exp(X)\rightarrow X$ the
associated $\Dst$-algebra $\alpha\after\sigma\colon\Dst(X) \rightarrow
X$ is then also continuous.
\end{lemma}

\begin{myproof}
Lemma~\ref{ExpAlgCHContinuousLem} states that $\Exp$-algebras
$\alpha\colon\Exp(X)\rightarrow X$ are continuous. Hence
$\alpha\after\sigma \colon \Dst(X) \rightarrow X$, as composition of
continuous maps, is also continuous. \QED
\end{myproof}

The following property of the map of monads $\Dst\Rightarrow\Exp$ will
play a crucial role.

\begin{proposition}
\label{DstExpDenseProp}
The inclusions $\sigma_{X} \colon \Dst(X) \rightarrowtail \Exp(X)$ are
dense: the topological closure of $\Dst(X)$ is the whole of $\Exp(X)$.
\end{proposition}

\begin{myproof}
We need to show that for each non-empty open $U\subseteq \Exp(X)$
there is a distribution $\varphi\in\Dst(X)$ with $\sigma(\varphi)\in
U$. By Proposition~\ref{ExpTopologyProp} we may assume $U$ is of the
form $U = \square_{s_{1}}(p_{1}) \cap \cdots \cap
\square_{s_{m}}(p_{m})$, for certain $s_{i} \in [0,1]\cap \mathbb{Q}$
and $p_{i}\in [0,1]^{X}$. For convenience we do the proof for
$m=2$. Since $U$ is non-empty there is some inhabitant
$h\in\square_{s_{1}}(p_{2}) \cap \square_{s_2}(p_{2})$. Thus $h(p_{i})
> s_{i}$.  We claim there are simple functions $q_{i} \leq p_{i}$ with
$h(q_{i}) > s_{i}$.

In general, this works as follows. If $h(p) > s$, write $p =
\lim_{n}p_{n}$ for simple functions $p_{n} \leq p$, like in
Lemma~\ref{SimplePredLem}. Then $h(p) = \lim_{n} h(p_{n}) > s$. Hence
$h(p_{n})>s$ for some simple $p_{n} \leq p$.

In a next step we write the simple functions as weighted sum of
characteristic functions, like in~\eqref{SimpleFunNormalFormEqn}.
Thus, let
$$\begin{array}{rclcrcl}
q_{1}
& = &
\ovee_{j}\,r_{j}\charac{U_j}
& \quad\mbox{and}\quad &
q_{2}
& = &
\ovee_{k}\,t_{k}\charac{V_k},
\end{array}$$

\noindent where these $U_{j}\subseteq X$ and $V_{k}\subseteq X$ form
non-empty partitions, each covering $X$. We form a new, refined
partition $(W_{\ell}\subseteq X)_{\ell\in L}$ consisting of the
non-empty intersections $U_{j}\cap V_{j}$, and choose $x_{\ell}\in
W_{\ell}$. Then:
\begin{itemize}
\item $\begin{array}{rcccccl}
\sum_{\ell}h(\charac{W_{\ell}})
& = &
h(\ovee_{\ell}\,\charac{W_{\ell}})
& = &
h(\charac{X})
& = &
1.
\end{array}$

\item There are subsets $L_{j}\subseteq L$ so that each $U_{j}\subseteq X$
can be written as disjoint union $U_{j} = \bigcup_{\ell\in L_{j}}W_{\ell}$.

\item Similarly, $V_{k} = \bigcup_{\ell\in L_{k}}W_{\ell}$ for subsets
$L_{k}\subseteq L$.
\end{itemize}

\noindent We take as distribution $\varphi =
\sum_{\ell\in L}h(\charac{W_{\ell}})\,x_{\ell} \in \Dst(X)$. Then
$\sigma(\varphi) \in \square_{s_i}(p_{i})$. We do the proof for $i=1$.
$$\begin{array}[b]{rcl}
\sigma(\varphi)(p_{1})
& = &
\sum_{\ell\in L}\varphi(x_{\ell})\cdot p_{1}(x_{\ell}) \\
& \geq & 
\sum_{\ell\in L}h(\charac{W_{\ell}})\cdot q_{1}(x_{\ell}) \\
& = &
\sum_{j}\sum_{\ell\in L_{j}}h(\charac{W_{\ell}})\cdot q_{1}(x_{\ell}) \\
& = &
\sum_{j}\sum_{\ell\in L_{j}}h(\charac{W_{\ell}})\cdot r_{j} \\
& = &
\sum_{j}h(\ovee_{\ell\in L_{j}}\charac{W_{\ell}})\cdot r_{j} \\
& = &
\sum_{j}h(\charac{U_{j}})\cdot r_{j} \\
& = &
h(\ovee_{j}\,r_{j}\cdot\charac{U_{j}}) \\
& = &
h(p_{1}) \\
& > &
s_{1}.
\end{array}\eqno{\QEDbox}$$
\end{myproof}

\begin{corollary}
\label{UFDstExpOntoCor}
Each map $\UF(\Dst(X)) \rightarrow \Exp(X)$, described in
Example~\ref{DstExpNatroRem}.(3), is onto (surjective).
\end{corollary}

\begin{myproof}
Since $\Dst(X) \rightarrowtail \Exp(X)$ is dense, each $h\in\Exp(X)$
is a limit of elements in $\Dst(X)$. Such limits can be described for
instance via nets or via ultrafilters. In the present context we
choose the latter approach. Thus there is an ultrafilter
$\calF\in\UF(\Dst(X))$ such that $h$ is the limit of this ultrafilter
$\UF(\sigma)(\calF)\in\UF(\Exp(X))$, when mapped to $\Exp(X)$. The
limit is expressed via the ultrafilter algebra $\mu \after \tau \colon
\UF(\Exp(X)) \rightarrow \Exp(X)$. This means that $(\mu \after \tau
\after \UF(\sigma))(\calF) = h$.  \QED
\end{myproof}

\auxproof{
Each monad on \Sets is strong. Hence the expectation monad also comes
with a strength map $\st\colon \Exp(X)\times Y \rightarrow
\Exp(X\times Y)$. It is defined as:
$$\begin{array}{rcl}
\st(h,y)
& = &
\lamin{p}{[0,1]^{X\times Y}}{h\big(p(\tuple{-,y})\big)}.
\end{array}$$

\noindent The expectation monad is not commutative: the two resulting
maps $\Exp(X)\times\Exp(Y) \rightrightarrows \Exp(X\times Y)$ are not
the same. Also the ultrafilter monad is not commutative.

There is a swapped version $\st'\colon X\times\Exp(Y) \rightarrow
\Exp(X\times Y)$ defined via the swap map $\gamma = \tuple{\pi_{2},
\pi_{1}}$, namely as:
$$\begin{array}{rcl}
\st'(x,k)
\hspace*{\arraycolsep} = \hspace*{\arraycolsep}
\big(\Exp(\gamma) \after \st \after \gamma\big)(x,k)
& = &
\Exp(\gamma)\big(\st(k,x)\big) \\
& = &
\lamin{p}{[0,1]^{X\times Y}}{\st(k,x)(p \after \gamma)} \\
& = &
\lamin{p}{[0,1]^{X\times Y}}{k\big(p \after \gamma \after 
   \tuple{-,x}\big)} \\
& = &
\lamin{p}{[0,1]^{X\times Y}}{k\big(p(\tuple{x,-})\big)}.
\end{array}$$

\noindent There are now two maps $\Exp(X)\times\Exp(Y)
\rightrightarrows \Exp(X\times Y)$. We show that they do not coincide.
$$\begin{array}{rcl}
\big(\mu \after \Exp(\st') \after \st\big)(h,k)
& = &
\mu\big(\Exp(\st')(\st(h,k))\big) \\
& = &
\lamin{p}{[0,1]^{X\times Y}}{\Exp(\st')(\st(h,k))
   \big(\lamin{g}{\Exp(X\times Y)}{g(p)}\big)} \\
& = &
\lamin{p}{[0,1]^{X\times Y}}{\st(h,k)\big(
   (\lamin{g}{\Exp(X\times Y)}{g(p)}) \after \st'\big)} \\
& = &
\lamin{p}{[0,1]^{X\times Y}}{h\big(
   (\lamin{g}{\Exp(X\times Y)}{g(p)}) \after \st' \after
   \tuple{-,k}\big)} \\
& = &
\lamin{p}{[0,1]^{X\times Y}}{h\big(\lamin{x}{X}{\st'(x,k)(p)}\big)} \\
& = &
\lamin{p}{[0,1]^{X\times Y}}{h\big(\lamin{x}{X}{k(p(\tuple{x,-}))}\big)} \\
\end{array}$$

\noindent In the other order we get:
$$\begin{array}{rcl}
\big(\mu \after \Exp(\st) \after \st'\big)(h,k)
& = &
\mu\big(\Exp(\st)(\st'(h,k))\big) \\
& = &
\lamin{p}{[0,1]^{X\times Y}}{\Exp(\st)(\st'(h,k))
   \big(\lamin{g}{\Exp(X\times Y)}{g(p)}\big)} \\
& = &
\lamin{p}{[0,1]^{X\times Y}}{\st'(h,k)\big(
   (\lamin{g}{\Exp(X\times Y)}{g(p)}) \after \st\big)} \\
& = &
\lamin{p}{[0,1]^{X\times Y}}{k\big(
   (\lamin{g}{\Exp(X\times Y)}{g(p)}) \after \st \after
   \tuple{h,-}\big)} \\
& = &
\lamin{p}{[0,1]^{X\times Y}}{k\big(\lamin{y}{Y}{\st(h,y)(p)}\big)} \\
& = &
\lamin{p}{[0,1]^{X\times Y}}{k\big(\lamin{y}{Y}{h(p(\tuple{-,y}))}\big)} \\
\end{array}$$
}

\section{Algebras of the expectation monad}\label{ExpAlgSec}

This section describes algebras of the expectation monad via
barycenters of measures. It leads to an equivalence of categories
between `observable' algebras and `observable' convex compact
Hausdorff spaces. We shall write $\CCH$ for the category of these
convex compact Hausdorff spaces, with affine continuous maps between
them.

We start with the unit interval $[0,1]$. It is both compact Hausdorff
and convex.  Hence it carries algebras $\UF([0,1])\rightarrow[0,1]$
and $\Dst([0,1])\rightarrow [0,1]$. This interval also carries an
algebra structure for the expectation monad.

\begin{lemma}
\label{UnitIntExpAlgLem}
The unit interval $[0,1]$ carries an $\Exp$-algebra structure:
$$\xymatrix{
\Exp([0,1])\ar[r]^-{\ev_{\idmap}} & [0,1]
& \quad\mbox{by}\quad &
h\ar@{|->}[r] & h(\idmap[{[0,1]}]).
}$$

\noindent More generally, for an arbitrary set $A$ the 
set of (all) functions $[0,1]^{A}$ carries an $\Exp$-algebra
structure:
$$\xymatrix@C-.5pc{
\Exp([0,1]^{A})\ar[r] & [0,1]^{A}
& \mbox{namely} &
h\ar@{|->}[r] & \lamin{a}{A}{h\big(\lamin{f}{[0,1]^{A}}{f(a)}\big)}.
}$$
\end{lemma}

\begin{myproof}
It is easy to see that the evaluation-at-identity map
$\ev_{\idmap}\colon\Exp([0,1]) \rightarrow [0,1]$ is an algebra. We
explicitly check the details:
$$\begin{array}{ccc}
\begin{array}{rcl}
\big(\ev_{\idmap} \after \eta\big)(x)
& = &
\ev_{\idmap}\big(\eta(x)\big) \\
& = &
\eta(x)(\idmap) \\
& = &
\idmap(x) \\
& = &
x
\end{array}
& \qquad &
\begin{array}{rcl}
\big(\ev_{\idmap} \after \Exp(\ev_{\idmap})\big)(H)
& = &
\ev_{\idmap}\big(\Exp(\ev_{\idmap})(H)\big) \\
& = &
\Exp(\ev_{\idmap})(H)(\idmap) \\
& = &
H(\idmap \after \ev_{\idmap}) \\
& = &
H\big(\lamin{k}{\Exp([0,1])}{k(\idmap)}\big) \\
& = &
\mu(H)(\idmap) \\
& = &
\ev_{\idmap}\big(\mu(H)\big) \\
& = &
\big(\ev_{\idmap} \after \mu\big)(H).
\end{array}
\end{array}$$

\noindent Since Eilenberg-Moore algebras are closed under products,
there is also an $\Exp$-algebra on $[0,1]^{A}$. \QED

\auxproof{ 
Let's write this map as $\cc\colon\Exp([0,1]^{A}) \rightarrow [0,1]^{A}$.
We explicitly check that it is an algebra:
$$\begin{array}{ccc}
\begin{array}{rcl}
\big(\cc \after \eta\big)(f)(a)
& = &
\cc\big(\eta(f)\big)(a) \\
& = &
\eta(f)(\lam{g}{g(a)}) \\
& = &
f(a) \\
& = &
\idmap(f)(a)
\end{array}
& \quad &
\begin{array}{rcl}
\big(\cc \after \Exp(\cc)\big)(F)(a)
& = &
\cc\big(\Exp(\cc)(F)\big)(a) \\
& = &
\Exp(\cc)(F)(\lam{g}{g(a)}) \\
& = &
F((\lam{g}{g(a)}) \after \cc) \\
& = &
F\big(\lam{h}{\cc(h)(a)}\big) \\
& = &
F\big(\lam{h}{h(\lam{f}{f(a)})}\big) \\
& = &
\mu(F)(\lam{f}{f(a)}) \\
& = &
\cc\big(\mu(F)\big)(a) \\
& = &
\big(\cc \after \mu\big)(F)(a).
\end{array}
\end{array}$$
}
\end{myproof}

From Corollaries~\ref{ExpAlgCHCor} and~\ref{ExpAlgConvCor} we know
that the underlying set $X$ of an algebra $\Exp(X) \rightarrow X$ is
both compact Hausdorff and convex. Additionally,
Lemma~\ref{ExpAlgDstContinuousLem} says that the algebra $\Dst(X)
\rightarrow X$ is continuous. 

We first characterize algebra maps.

\begin{lemma}
\label{ExpAlgHomLem}
Consider Eilenberg-Moore algebras $\smash{(\Exp(X)
  \xrightarrow{\alpha} X)}$ and $\smash{(\Exp(Y) \xrightarrow{\beta}
  Y)}$. A function $f\colon X\rightarrow Y$ is an algebra homomorphism
if and only if it is both continuous and affine, that is, iff the
following two diagrams commute.
$$\xymatrix@R-.5pc@C-.5pc{
\UF(X)\ar[d]_{\alpha\after\tau}\ar[rr]^-{\UF(f)} & &
   \UF(Y)\ar[d]^{\beta\after\tau}
& &
\Dst(X)\ar[d]_{\alpha\after\sigma}\ar[rr]^-{\Dst(f)} & &
   \Dst(Y)\ar[d]^{\beta\after\sigma} \\
X\ar[rr]_-{f} & & Y
& &
X\ar[rr]_-{f} & & Y
}$$

\noindent Thus, the functor $\Alg(\Exp) \rightarrow \CCH$ is
full and faithful.
\end{lemma}

\begin{myproof}
If $f$ is an algebra homomorphism, then $f \after \alpha = \beta
\after \Exp(f)$. Hence the two rectangles above commute by naturality
of $\tau$ and $\sigma$.

For the (if) part we use the property from
Proposition~\ref{DstExpDenseProp} that the maps $\sigma_{X}\colon
\Dst(X) \rightarrowtail \Exp(X)$ are dense monos. This means that for
each map $g\colon \Dst(X) \rightarrow Z$ into a Hausdorff space $Z$
there is at most one continuous $h\colon \Exp(X) \rightarrow Z$ with
$h \after \sigma = g$.  We use this property as follows.
$$\xymatrix@R+1pc{
\Dst(X)\ar@{ >->}[rr]^-{\sigma}_-{\text{dense}}\ar[drr] & &
   \Exp(X)\ar@<+1ex>[d]^{\beta\after\Exp(f)}\ar@<-1ex>[d]_{f\after \alpha} \\
& & Y
}$$

\noindent The triangle commutes for both maps since $f$ is affine:
$$\begin{array}{rcccl}
f \after \alpha \after \sigma
& = &
\beta \after \sigma \after \Dst(f)
& = &
\beta \after \Exp(f) \after \sigma.
\end{array}$$

\noindent Also, both vertical maps are continuous, by
Lemma~\ref{ExpAlgCHContinuousLem}. Hence $f\after \alpha = \beta
\after \Exp(f)$, so that $f$ is an algebra homomorphism. \QED

\auxproof{
Or by the following two diagram chases.
$$\xymatrix@R-.5pc@C-.5pc{
\UF(\Exp(X))\ar[d]_{\tau}\ar[r]^-{\UF(\alpha)} &
   \UF(X)\ar[d]^{\tau}\ar[r]^-{\UF(f)} & \UF(Y)\ar[d]^{\tau} 
& 
\UF(\Exp(X))\ar[d]_{\tau}\ar[r]^-{\UF(\Exp(f))} &
   \UF(\Exp(Y))\ar[d]^{\tau}\ar[r]^-{\UF(\beta)} & \UF(Y)\ar[d]^{\tau} 
\\
\Exp^{2}(X)\ar[d]_{\mu}\ar[r]^-{\Exp(\alpha)} &
   \Exp(X)\ar[d]^{\alpha} & \Exp(Y)\ar[d]^{\beta}
& 
\Exp^{2}(X)\ar[d]_{\mu}\ar[r]^-{\Exp^{2}(f)} &
   \Exp^{2}(Y)\ar[d]^{\mu}\ar[r]^-{\Exp(\beta)} & \Exp(Y)\ar[d]^{\beta}
\\
\Exp(X)\ar[r]_-{\alpha} & X\ar[r]_-{f} & Y
&
\Exp(X)\ar[r]_-{\Exp(f)} & \Exp(Y)\ar[r]_-{\beta} & Y
}$$
}
\end{myproof}

For convex compact Hausdorff spaces $X,Y\in\CCH$ one (standardly)
writes $\mathcal{A}(X,Y) = \CCH(X,Y)$ for the homset of affine
continuous functions $X\rightarrow Y$.  In light of the previous
result, we shall also use this notation $\mathcal{A}(X,Y)$ when $X,Y$
are carriers of $\Exp$-algebras, in case the algebra structure is
clear from the context.

The next result gives a better understanding of $\Exp$-algebras: it
shows that such algebras send measures to barycenters (like for
instance in~\cite{Keimel09}).

\begin{proposition}
\label{ExpAlgBarycenterProp}
Assume an $\Exp$-algebra $\smash{\Exp(X)\xrightarrow{\alpha}X}$.
For each (algebra) map $q\in\mathcal{A}(X,[0,1])$ the following
diagram commutes.
$$\xymatrix@R-.5pc{
\Exp(X)\ar[d]_{\alpha}\ar[drr]^-{\quad\ev_{q} = \lam{h}{h(q)}} \\
X\ar[rr]_-{q} & & [0,1]
}$$

\noindent This says that $x=\alpha(h)\in X$ is a \emph{barycenter} for
$h\in\Exp(X)$, in the sense that $q(x) = h(q)$ for all affine
continuous $q\colon X \rightarrow [0,1]$.
\end{proposition}

\begin{myproof}
Since $\ev_{q} = \ev_{\idmap} \after \Exp(q)$ the above triangle can
be morphed into a rectangle expressing that $q$ is a map of algebras:
$$\xymatrix@R-.5pc{
\Exp(X)\ar[d]_{\alpha}\ar[drr]^-{\ev_{q}}\ar[rr]^-{\Exp(q)} & & 
   \Exp([0,1])\ar[d]^{\ev_{\idmap}} \\
X\ar[rr]_-{q} & & [0,1]
}$$

\noindent where $\ev_{\idmap}$ is the $\Exp$-algebra on $[0,1]$ from
Lemma~\ref{UnitIntExpAlgLem}.
\end{myproof}

Now that we have a reasonable grasp of $\Exp$-algebras, namely as
convex compact Hausdorff spaces with a barycentric operation, we wish
to comprehend how such algebras arise. We first observe that measures
in $\Exp(X)$ in the images of $\Dst(X)\rightarrowtail \Exp(X)$ and
$\UF(X)\rightarrowtail\Exp(X)$ have barycenters, if $X$ carries
appropriate structure.

\begin{lemma}
\label{SimpleBarycenterLem}
Assume $X$ is a convex compact Hausdorff space, described via $\Dst$-
and $\UF$-algebra structures $\cv\colon\Dst(X)\rightarrow X$ and
$\ch\colon\UF(X)\rightarrow X$. Then:
\begin{enumerate}
\item $\cv(\varphi)\in X$ is a barycenter of $\sigma(\varphi)\in
  \Exp(X)$, for $\varphi\in\Dst(X)$;

\item $\ch(\calF)\in X$ is a barycenter of $\tau(\calF)\in\Exp(X)$,
  for $\calF\in\UF(X)$.
\end{enumerate}
\end{lemma}

\begin{myproof}
We write $\cv_{[0,1]}\colon\Dst([0,1])\rightarrow [0,1]$ and
$\ch_{[0,1]}\colon\UF([0,1])\rightarrow [0,1]$ for the convex and compact
Hausdorff structure on the unit interval. Then for
$q\in\mathcal{A}(X,[0,1])$,
$$\begin{array}[b]{rcll}
q\big(\cv(\varphi)\big)
& = &
\cv_{[0,1]}\big(\Dst(q)(\varphi)\big)
   \qquad & \mbox{since $q$ is affine} \\
& = &
\cv_{[0,1]}\big(\sum_{i}r_{i}q(x_{i})\big)
   & \mbox{if }\varphi=\sum_{i}r_{i}x_{i} \\
& = &
\sum_{i}r_{i}\cdot q(x_{i}) \\
& = &
\sigma(\varphi)(q) \\
q\big(\ch(\calF)\big)
& = &
\ch_{[0,1]}\big(\UF(q)(\calF)\big)
   & \mbox{since $q$ is continuous} \\
& = &
\tau(\calF)(q).
\end{array}\eqno{\QEDbox}$$
\end{myproof}

We call a convex compact Hausdorff space $X$ \emph{observable} if the
collection of affine continuous maps $X\rightarrow[0,1]$ is jointly
monic. This means that $x=x'$ holds if $q(x)=q(x')$ for all
$q\in\mathcal{A}(X,[0,1])$. In a similar manner we call an
$\Exp$-algebra observable if its underlying convex compact Hausdorff
space is observable. This yields full subcategories $\CCHobs
\hookrightarrow \CCH$ and $\Algobs(\Exp) \hookrightarrow \Alg(\Exp)$.
By definition, $[0,1]$ is a cogenerator in these categories
$\CCHobs$ and $\Algobs(\Exp)$.


\begin{proposition}
\label{BarycenterExistProp}
Assume $X$ is a convex compact Hausdorff space, described via $\Dst$-
and $\UF$-algebra structures $\cv\colon\Dst(X)\rightarrow X$ and
$\ch\colon\UF(X)\rightarrow X$.
\begin{enumerate}
\item Via the Axiom of Choice one obtains a function
  $\alpha\colon\Exp(X)\rightarrow X$ such that $\alpha(h)\in X$ is a
  barycenter for $h\in\Exp(X)$; that is, $q(\alpha(h)) = h(q)$ for
  each $q\in\mathcal{A}(X,[0,1])$.

\item If $X$ is observable, there is precisely one such $\alpha\colon
  \Exp(X) \rightarrow X$; moreover, it is an $\Exp$-algebra; and its
  induced convex and topological structures are the original ones on
  $X$, as expressed via the commuting triangles:
$$\xymatrix@R-.5pc{
\Dst(X)\ar@{ >->}[r]^-{\sigma}\ar[dr]_{\cv} & 
   \Exp(X)\ar[d]^{\alpha} & \UF(X)\ar@{ >->}[l]_-{\tau}\ar[dl]^{\ch} \\
& X &
}$$

\noindent This yields a functor $\CCHobs \rightarrow \Algobs(\Exp)$.
\end{enumerate}
\end{proposition}

\begin{myproof}
Recall from Corollary~\ref{UFDstExpOntoCor} that the function $\mu
\after \tau \after \UF(\sigma) \colon \UF(\Dst(X)) \rightarrow
\Exp(X)$ is surjective. Using the Axiom of Choice we choose a section
$s\colon \Exp(X) \rightarrow \UF(\Dst(X))$ with $\mu \after \tau
\after \UF(\sigma) \after s = \idmap[\Exp(X)]$. We now obtain, via the
choice of $s$, a map $\alpha\colon \Exp(X)\rightarrow X$ in:
$$\xymatrix{
\UF(\Dst(X))\ar@{->>}[rr]_-{\mu \after \tau \after \UF(\sigma)}
   \ar@/_1ex/[drr]_{\ch \after \UF(\cv)\quad} & & 
   \Exp(X)\ar[d]^{\alpha=\ch \after \UF(\cv)\after s}
      \ar@/_3ex/[ll]_{s} \\
& & X
}$$

\noindent We show that $\alpha(h)\in X$ is a barycenter for the measure
$h\in\Exp(X)$. For each $q\in\mathcal{A}(X,[0,1])$ one has:
$$\begin{array}{rcll}
h(q)
& = &
\big(\mu \after \tau \after \UF(\sigma) \after s\big)(h)(q) \\
& = &
\mu\big((\tau \after \UF(\sigma) \after s)(h)\big)(q) \\
& = &
\big(\tau \after \UF(\sigma) \after s\big)(h)(\ev_{q}) \\
& = &
\big(\ch_{[0,1]} \after \UF(\ev_{q}) \after \UF(\sigma) \after s\big)(h) \\
& = &
\big(\ch_{[0,1]} \after \UF(\lam{\varphi}{\ev_{q}(\sigma(\varphi))}) 
   \after s\big)(h) \\
& = &
\big(\ch_{[0,1]} \after \UF(\lam{\varphi}{\cv_{[0,1]}(\Dst(q)(\varphi))}) 
   \after s\big)(h) \quad
   & \mbox{see Remark~\ref{DstExpNatroRem}.(1)} \\
& = &
\big(\ch_{[0,1]} \after \UF(\cv_{[0,1]} \after \Dst(q))
   \after s\big)(h) \\
& = &
\big(\ch_{[0,1]} \after \UF(q \after \cv) \after s\big)(h)
   & \mbox{since $q$ is affine} \\
& = &
\big(q \after \ch \after \UF(\cv) \after s\big)(h)
   & \mbox{since $q$ is continuous} \\
& = &
\big(q \after \alpha\big)(h) \\
& = &
q(\alpha(h)).
\end{array}$$

\auxproof{
Earlier, partial proofs in terms of nets.

Next we make crucial use of the fact that
$\sigma\colon\Dst(X)\rightarrowtail \Exp(X)$ is dense, see
Proposition~\ref{DstExpDenseProp}. Each $h\in\Exp(X)$ may thus be
written as limit of a net $\psi\colon J\rightarrow \Dst(X)$, where $J$
is a directed set. Thus $h = \lim_{j\in J}\sigma(\psi_{j})$. As just
noted, each $\sigma(\psi_{j})\in\Exp(X)$ has a barycenter $x_{j}$.
This yields a net $x\colon J \rightarrow X$. Since $X$ is compact,
there is a subnet $y\colon K \rightarrow X$, where $K \subseteq J$ is
cofinal, with limit $x_{h} = \lim_{k\in K}y_{k}$. This $x_{h}\in X$
yields a barycenter for $h$ since:
$$\begin{array}{rcll}
p(x_{h})
& = &
\lim_{k}p(y_{k}) 
   & \mbox{because $p$ is continuous} \\
& = &
\lim_{k}\sigma(\psi_{k})(p)
   & \mbox{since $y_{k}$ is barycenter of $\sigma(\psi_{k})$} \\
& = &
\lim_{k} \ev_{p}(\sigma(\psi_{k})) 
   & \mbox{where $\ev_{p} = \lam{g}{g(p)}$} \\
& = &
\ev_{p}(\lim_{k}\sigma(\psi_{k})) \qquad
   & \mbox{by continuity of $\ev_p$, see 
   Lemma~\ref{ExpAlgCHContinuousLem}} \\
& = &
\ev_{p}(\lim_{j}\sigma(\psi_{j})) 
   & \mbox{a limit of a subnet yields the net's limit} \\
& = &
\ev_{p}(h) \\
& = &
h(p).
\end{array}$$
}

For the second point, assume the collection of maps
$q\in\mathcal{A}(X,[0,1])$ is jointly monic. Barycenters are then
unique, since if both $x,x'\in X$ satisfy $q(x) = h(q) = q(x')$ for
all $q\in\mathcal{A}(X,[0,1])$, then $x=x'$. Hence the function
$\alpha\colon\Exp(X)\rightarrow X$ picks barycenters, in a unique
manner. We need to prove the algebra equations (see the beginning of
Section~\ref{MonadSec}). They are obtained via the barycentric
property $q(\alpha(h)) = h(q)$ and observability. First, the equation
$\alpha \after \eta = \idmap$ holds, since for each $x\in X$ and
$q\in\mathcal{A}(X,[0,1])$,
$$\begin{array}{rcccccccl}
q\big((\alpha \after \eta)(x))\big)
& = &
q\big(\alpha(\eta(x))\big)
& = &
\eta(x)(q)
& = &
q(x)
& = &
q\big(\idmap(x)\big).
\end{array}$$

\noindent In the same way we obtain the equation $\alpha \after \mu = 
\alpha \after \Exp(\alpha)$. For $H\in\Exp^{2}(X)$ we have:
$$\begin{array}[b]{rcl}
\big(q \after \alpha \after \mu\big)(H)
& = &
q\big(\alpha(\mu(H))\big) \\
& = &
\mu(H)(q) \\
& = &
H\big(\lamin{k}{\Exp(X)}{k(q)}\big) \\
& = &
H\big(\lamin{k}{\Exp(X)}{q(\alpha(k))}\big) \\
& = &
H\big(q \after \alpha\big) \\
& = &
\Exp(\alpha)(H)(q) \\
& = &
q\big(\alpha(\Exp(\alpha)(H))\big) \\
& = &
\big(q \after \alpha \after \Exp(\alpha)\big)(H).
\end{array}$$

We need to show that $\alpha$ induces the original convexity
and topological structures. Since barycenters are unique, the equations
$\alpha(\sigma(\varphi)) = \cv(\varphi)$ and $\alpha(\tau(\calF)) =
\ch(\calF)$ follow directly from Lemma~\ref{ExpAlgHomLem}. 

Finally, we need to check functoriality. So assume $f\colon X
\rightarrow Y$ is a map in $\CCHobs$, and let $\alpha\colon\Exp(X)
\rightarrow X$ and $\beta\colon\Exp(Y)\rightarrow Y$ be the 
induced algebras obtained by picking barycenters. We need to prove
$\beta \after \Exp(f) = f \after \alpha$. Of course we use
that $Y$ is observable. For $h\in\Exp(X)$, one has for all
$q\in\mathcal{A}(Y,[0,1])$,
$$\begin{array}[b]{rcl}
q\big(\beta(\Exp(f)(h))\big)
& = &
\Exp(f)(h)(q) \\
& = &
h\big(q \after f\big) \\
& = &
(q \after f)(\alpha(h)) \\
& = &
q\big(f(\alpha(h))\big).
\end{array}\eqno{\QEDbox}$$
\end{myproof}

In the approach followed above barycenters are obtained via the Axiom
of Choice. Alternatively, they can be obtained via the Hahn-Banach
theorem, see for instance~\cite[Prop.~I.2.1]{Alfsen71}.

\begin{theorem}
\label{AlgobsCCHobsThm}
There is an isomorphism $\Algobs(\Exp) \cong \CCHobs$ between the
categories of observable $\Exp$-algebras and observable convex compact
Hausdorff spaces in a situation:
$$\xymatrix@R-0.5pc{
\raisebox{.4em}{$\Algobs(\Exp)$}\ar@{^(->}[d]\ar@/_1.3ex/[rr]^-{\cong} & & 
   \CCHobs\ar@/_1.3ex/[ll]\ar@{^(->}[d] \\
\Alg(\Exp)\ar[rr]_-{\text{full \& faithful}} & & \CCH
}$$
\end{theorem}

\begin{myproof}
Obviously the full and faithful functor $\Alg(\Exp) \rightarrow \CCH$
from Lemma~\ref{ExpAlgHomLem} restricts to $\Algobs(\Exp) \rightarrow
\CCHobs$. We need to show that it is an inverse to the functor
$\CCHobs\rightarrow\Algobs(\Exp)$ from
Proposition~\ref{BarycenterExistProp}.(2). 
\begin{itemize}
\item Starting from an algebra $\alpha\colon\Exp(X)\rightarrow X$, we
  know by Proposition~\ref{ExpAlgBarycenterProp} that $\alpha(h)$ is a
  barycenter for $h\in\Exp(X)$. The underlying set $X$ is an
  observable convex compact Hausdorff space. This structure gives by
  Proposition~\ref{BarycenterExistProp}.(2) rise to an algebra
  $\alpha'\colon \Exp(X)\rightarrow X$ such that $\alpha'(h)$ is
  barycenter for $h$. Since $X$ is observable, barycenters are unique,
  and so $\alpha'(h) = \alpha(h)$.

\item Starting from an observable convex compact Hausdorff space $X$,
  we obtain an algebra $\alpha\colon\Exp(X)\rightarrow X$ by
  Proposition~\ref{BarycenterExistProp}.(2), whose induced convex and
  topological structure is the original one. \QED
\end{itemize}
\end{myproof}

Thus we have characterized \emph{observable} $\Exp$-algebras. The
characterization of arbitrary $\Exp$-algebras remains open. Possibly
the functor $\Alg(\Exp)\rightarrow\CCH$ is (also) an isomorphism.
For the duality in the next section the characterization of
observable algebras is sufficient.

\auxproof{
Older approach, more closely following~\cite{Alfsen71}.
\begin{lemma}
Assume $X$ is a non-empty convex compact Hausdorff space. Then:
\begin{enumerate}
\item For each $h\in\Exp(X)$ and $p\in\mathcal{A}(X,[0,1])$ there is
  an $x\in X$ with $p(x) = h(p)$.

\item For each $h\in\Exp(X)$, there is a `barycenter' $x_{h}\in X$
  such that $p(x_{h}) = h(p)$ holds for all
  $p\in\mathcal{A}(X,[0,1])$.
\end{enumerate}
\end{lemma}

\begin{myproof}
For the first point define top and bottom elements $t,b\in[0,1]$ for
$p$ as:
$$\begin{array}{rclcrcl}
t
& = &
\sup\set{p(x)}{x\in X}
& \quad\mbox{and}\quad &
b
& = &
\inf\set{p(x)}{x\in X}.
\end{array}$$

\noindent Then $b\cdot \mathbf{1} \leq p \leq t\cdot\mathbf{1}$, where
$\mathbf{1} = \lam{x}{1} \in [0,1]^{X}$. Since $h(b\cdot\mathbf{1}) =
b\cdot h(\mathbf{1}) = b\cdot 1 = b$, and similarly for $t$, we get:
$b \leq h(p) \leq t$. If $b=t$, $p$ is a constant function $p = b\cdot
\mathbf{1}$ and any $x\in X$ satisfies $p(x) = b = h(p)$. Hence assume
$b<t$. Then we can write $h(p)$ as convex combination $h(p) = \lambda
b + (1-\lambda)t$, for $\lambda = \frac{t-h(p)}{t-b}$.

The image $\set{p(x)}{x\in X} \subseteq [0,1]$ is compact, since $p$
is continuous, and thus closed. Hence $p$ actually reaches its bottom
and top: there are $x_{b}, x_{t}\in X$ with $p(x_{b}) = b$ and
$p(x_{t}) = t$. Now put $x = \lambda x_{b} + (1-\lambda)x_{t}\in X$,
using that $X$ is convex. But then we are done since $p$ is affine:
$$\begin{array}{rcccccl}
p(x)
& = &
\lambda p(x_{b}) + (1-\lambda)p(x_{t})
& = &
\lambda b + (1-\lambda)t
& = &
h(p).
\end{array}$$

For the second point we follow the argument in the proof
of~\cite[Prop.~I.2.1]{Alfsen71}. For $p\in\mathcal{A}(X,[0,1])$ write
$X_{p} = \setin{x}{X}{p(x) = h(p)}$. By the previous point each subset
$X_{p}\subseteq X$ is non-empty. We claim that the collection
$(X_{p})$ has the finite intersection property. Thus, assume $p_{1},
\ldots, p_{n}\in \mathcal{A}(X,[0,1])$. Define a continuous function
$P\colon X \rightarrow \mathbb{R}^{n}$ by $P(x) = \tuple{p_{1}(x),
  \ldots, p_{n}(x)}$. Then:
$$\begin{array}{rcl}
X_{p_{1}} \cap \cdots \cap X_{p_n} \neq \emptyset
& \Longleftrightarrow &
\tuple{h(p_{1}), \ldots, h(p_{n})} \in 
\mathit{im}(P) = \set{P(x)}{x\in X}.
\end{array}$$

\noindent Assume the intersection on the left is empty. The singleton
$\{\tuple{h(p_{1}), \ldots, h(p_{n})}\} \subseteq \mathbb{R}^{n}$ is a
closed subset that is then disjoint from the compact subset
$\mathit{im}(P)\subseteq \mathbb{R}^{n}$. Hence the Hahn-Banach
separation theorem yields a continuous linear function $g\colon
\mathbb{R}^{n} \rightarrow \mathbb{R}$ with:
$$\begin{array}{rcl}
\sup\set{g(v)}{v\in\mathit{im}(P)}
& < &
g(h(p_{1}), \ldots, h(p_{n})).
\end{array}$$

\noindent Since $g\colon \mathbb{R}^{n} \rightarrow \mathbb{R}$ is
linear, we can write it as $g(x_{1}, \ldots, x_{n}) = a_{1}x_{1} +
\cdots + a_{n}x_{n}$ for certain $a_{i}\in\mathbb{R}$. We may assume
$a_{i} \geq 0$. 


\end{myproof}
}

We conclude this section with some further results on observability.
We show that \emph{observable} convex compact Hausdorff spaces can be
considered as part of an enveloping locally convex topological vector
space. This is the more common way of describing such structures,
see \textit{e.g.}~\cite{Alfsen71,AsimovE80}.

\begin{lemma}
Let $X$ be a convex compact Hausdorff space; write $A =
\mathcal{A}(X,[0,1])$. If $X$ is observable, there is (by definition)
an injection:
$$\xymatrix{
X\ar@{ >->}[rr]^-{x\mapsto \ev_{x}} & & [0,1]^{A} 
   \qquad\mbox{where}\qquad\ev_{x} = \lamin{q}{A}{q(x)}.
}$$

\begin{enumerate}
\item This map is both affine and continuous---where $[0,1]^{A}$ carries
the product topology.

\item Hence if $X$ is the carrier of an $\Exp$-algebra, this
map is a homomorphism of algebras---where $[0,1]^{A}$ carries the
$\Exp$-algebra structure from Lemma~\ref{UnitIntExpAlgLem}.
\end{enumerate}
\end{lemma}

\begin{myproof}
The second point follows from the first one via Lemma~\ref{ExpAlgHomLem},
so we only do point~1. Obviously, $x\mapsto\ev_{x}$ is affine. In order
to see that it is also continuous, assume we have a basic open set
$U\subseteq [0,1]^{A}$. The product topology says that $U$  is of the
form $U = \prod_{q\in A}U_{q}$, with $U_{q}\subseteq [0,1]$ open and 
$U_{q} \neq [0,1]$ for only finitely many $q$'s, say $q_{1}, \ldots,
q_{n}$. Thus:
$$\begin{array}{rcccl}
\ev^{-1}(U)
& = &
\set{x}{q_{1}(x) \in U_{q_1} \conjun \cdots \conjun q_{n}(x)\in U_{q_n}}
& = &
\bigcap_{i}q^{-1}_{i}(U_{q_i}).
\end{array}$$

\noindent This intersection of opens is clearly an open set of $X$. \QED

\auxproof{
We check explicitly that $\ev\colon X \rightarrow [0,1]^{A}$ is a map
of algebras in:
$$\xymatrix{
\Exp(X)\ar[d]_{\alpha}\ar[rr]^-{\Exp(\ev)} & &
   \Exp([0,1]^{A})\ar[d]^{\cc = \lam{h}{\lam{f}{h(\ev_{f})}}} \\
X\ar[rr]_-{\ev} & & [0,1]^{A}
}$$

\noindent Indeed for $h\in\Exp(X)$ and $f\in A = \mathcal{A}(X,[0,1])$,
$$\begin{array}{rcl}
\big(\cc \after \Exp(\ev)\big)(h)(f)
& = &
\cc\big(\Exp(\ev)(h)\big)(f) \\
& = &
\Exp(\ev)(h)(\ev_{f}) \\
& = &
h(\ev_{f} \after \ev) \\
& = &
h(\lamin{x}{X}{\ev_{f}(\ev_{x})}) \\
& = &
h(\lamin{x}{X}{\ev_{x}(f)}) \\
& = &
h(\lamin{x}{X}{f(x)}) \\
& = &
h(f) \\
& = &
h(\idmap \after f) \\
& = &
\Exp(f)(h)(\idmap) \\
& = &
\cc_{[0,1]}\big(\Exp(f)(h)\big) \\
& = &
f\big(\alpha(h)\big) \qquad \mbox{since }f\in\mathcal{A}(X,[0,1]) \\
& = &
\ev_{\alpha(h)}(f) \\
& = &
\big(\ev \after \alpha\big)(h)(f).
\end{array}$$
}
\end{myproof}

\begin{proposition}
\label{CCHobsLocConvProp}
Each \emph{observable} convex compact Hausdorff space $X\in\CCHobs$
occurs as subspace of a locally convex topological vector space,
namely via:
$$\xymatrix{
X \ar@{ >->}[rr] & & [0,1]^{A}\; \ar@{^(->}[rr] & & \mathbb{R}^{A}
}$$

\noindent where $A = \mathcal{A}(X,[0,1])$ like in the previous lemma.
\end{proposition}

\begin{myproof}
It is standard that the vector space $\mathbb{R}^{A}$ with product
topology is locally convex. We write $\mathcal{O}(X)$ for the original
compact Hausdorff topology on $X$ and $\mathcal{O}_{i}(X)$ for the
topology induced by the injection $X \rightarrowtail
\mathbb{R}^{A}$. The latter contains basic opens of the form
$q_{1}^{-1}(U_{1}) \cap \cdots \cap q_{n}^{-1}(U_{n})$ for $q_{i}\in
A=\mathcal{A}(X,[0,1])$ and $U_{i}\subseteq\mathbb{R}$ open.  Thus
$\mathcal{O}_{i}(X) \subseteq \mathcal{O}(X)$. We wish to use
Lemma~\ref{UniqueCHLem} to prove the equality $\mathcal{O}_{i}(X) =
\mathcal{O}(X)$. Since $\mathcal{O}(X)$ is compact we only need to
show that the induced topology $\mathcal{O}_{i}(X)$ is Hausdorff. This
is easy since $X$ is observable: if $x\neq x'$ for $x,x'\in X$, then
there is a $q\in\mathcal{A}(X,[0,1])$ with $q(x)\neq q(x')$ in
$[0,1]\subseteq \mathbb{R}$. Hence there are disjoint opens
$U,U'\subseteq\mathbb{R}$ with $q(x)\in U$ and $q(x')\in U'$. Thus
$q^{-1}(U), q^{-1}(U') \in\mathcal{O}_{i}(X)$ are disjoint (induced)
opens containing $x,x'$. \QED
\end{myproof}

\section{Algebras of the expectation monad and effect modules}\label{ExpEModSec}

In this section we relate algebras of the expectation monad to effect
modules via a (dual) adjunction. By suitable restriction this
adjunction gives rise to an equivalence (duality) between observable
$\Exp$-algebras and Banach effect modules. Via a combination with
Theorem~\ref{AlgobsCCHobsThm} we then get our main duality result (see
Theorem~\ref{EModAlgExpDualityThm} below).

We first return to Section~\ref{MndAdjSubsec}. When we apply
Lemma~\ref{ComposableAdjunctionLem} to the adjunctions involving
convex sets and effect modules in Diagram~\eqref{SetsConvEModDiag},
the (upper) comparison functor
in~\eqref{ComposableAdjunctionComparisonDiag} says that each effect
module $M\in\EMod$ gives rise to a $\Exp$-algebra on the homset
$\EMod(M,[0,1])$, namely:
\begin{equation}
\label{EModExpAlgDiag}
\vcenter{\xymatrix@R-2.5pc@C-1pc{
\Exp\Big(\EMod(M,[0,1])\Big)\ar[rr]^-{\alpha_M} & & \EMod(M,[0,1]) \\
h\ar@{|->}[rr] & & 
   \lamin{y}{M}{h\Big(\lamin{k}{\EMod(M,[0,1])}{k(y)}\Big)}
}}
\end{equation}

\noindent In order to simply notation we write:
$$\begin{array}{rcll}
S_{M} 
& = & 
\EMod(M,[0,1]) \quad & \mbox{for the set of ``states'' of $M$} \\
\alpha_{M}(h)(y)
& = &
h(\ev_{y}) & \mbox{where $\ev_{y} = \lam{k}{k(y)}$.}
\end{array}$$

\noindent Thus, Diagram~\eqref{ComposableAdjunctionComparisonDiag}
becomes:
\begin{equation}
\label{ExpComparisonDiag}
\vcenter{\xymatrix@R-.5pc{
\EMod\op\ar[dr]_{S_{(-)}}\ar[rr]^-{S_{(-)}} & &
   \Alg(\Exp)\ar[dl]^(0.4){(-)\after\sigma} \\
& \Alg(\Dst) \rlap{$\;=\Conv$}\ar[d] \\
& \Sets
}}
\end{equation}

\begin{proposition}
\label{EModExpAlgProp}
Consider for an effect module $M$, the $\Exp$-algebra
structure~\eqref{EModExpAlgDiag} on the homset of
states $S_{M} = \EMod(M, [0,1])$.
\begin{enumerate}
\item The induced topology is like the weak-* topology, with subbasic
opens
$$\begin{array}{rcl}
\square_{s}(y)
& = &
\setin{g}{S_{M}}{g(y) < s},
\end{array}$$

\noindent where $y\in M$ and $s\in[0,1]\cap\mathbb{Q}$. It thus
generalizes the topology on $\Exp(X) = \EMod([0,1]^{X},[0,1])$
in Proposition~\ref{ExpTopologyProp}.

\item This convex compact Hausdorff space $\EMod(M, [0,1])$ is
observable.
\end{enumerate}

\noindent Hence the states functor $S_{(-)}$ at the top
of~\eqref{ExpComparisonDiag} restricts to
$\EMod\op\rightarrow\Algobs(\Exp)$.
\end{proposition}

\begin{myproof}
The proof of Proposition~\ref{ExpTopologyProp} generalizes directly
from an effect module of the form $[0,1]^{X}$ to an arbitrary effect
module $M$.

\auxproof{
We show that the subsets $\neq\square_{s}(x)$ are closed in the
topology induced on $X$ by the algebra structure
$\alpha_{M}\after\tau\colon\UF(X)\rightarrow X$. For an ultrafilter
$\calF\in\UF(X)$, we have the implication required
in~\eqref{ClosedInvariantDiag}:
$$\begin{array}{rcll}
\lefteqn{\alpha_{M}(\tau(\calF))\in\neg\square_{s}(x)} \\
& \Longleftrightarrow &
\alpha_{M}(\tau(\calF))(x) = \tau(\calF)(\ev_{x}) \geq s \\
& \Longleftrightarrow &
\ch\big(\UF(\ev_{x})(\calF)\big) \in [s,1] &
   \mbox{where $\ch\colon\UF([0,1])\rightarrow [0,1]$} \\
& \Longleftarrow &
[s,1]\in\UF(\ev_{x})(\calF) &
   \mbox{since $[s,1]\subseteq [0,1]$ is closed} \\
& \Longleftrightarrow &
\neg\square_{s}(x) = \set{g}{g(x)\geq s} 
  = (\ev_{x})^{-1}\rlap{$([s,1])\in\calF.$}
\end{array}$$

The resulting topology is Hausdorff: if $f,g\colon M\rightarrow [0,1]$
in $\EMod$ are different, we may assume $f(x) < g(x)$ for some $x\in M$.
Hence there is an $s\in[0,1]\cap\mathbb{Q}$ with $f(x) < s < g(x)$.
Thus $f\in\square_{s}(x)$ and $g\in\square_{1-s}(x^{\perp})$. The
latter because $g(x^{\perp}) = 1 - g(x) < 1 - s$.

Thus the sets $\square_{s}(x)$ form a subbasis for a Hausdorff
topology on $\Exp(X)$. Hence by Lemma~\ref{UniqueCHLem} the
topology is the same one as the compact topology induced by
the algebra $\alpha_M$.
}

Next suppose $f,g\in S_{M} = \EMod(M,[0,1])$ satisfy $q(f) = q(g)$ for
each affine continuous $q\colon S_{M}\rightarrow [0,1]$. This
applies especially to the functions $\ev_{y} = \lam{k}{k(y)}$, which
are both continuous and affine. Hence $f(y)=\ev_{y}(f) = \ev_{y}(g) =
g(y)$ for each $y\in M$, and thus $f=g$. \QED

\auxproof{
Continuity holds by definition of the weak-* topology. These
$\ev_x$ are affine since:
$$\begin{array}{rcl}
\ev_{x}(\sum_{i}r_{i}h_{i})
& = &
\big(\sum_{i}r_{i}h_{i}\big)(x) \\
& = &
\sum_{i}r_{i}h_{i}(x) \\
& = &
\sum_{i}r_{i}\ev_{x}(h_{i}).
\end{array}$$
}

\end{myproof}

\auxproof{
First we shall see examples of the algebra
structure~\eqref{EModExpAlgDiag} in an indirect manner, via the
following translation result.

\begin{lemma}
\label{EModExpAlgTranslateLem}
Suppose we have a set $X$ for which there is an effect module $M$
with an isomorphism:
$$\xymatrix{
X\ar@/^2ex/[rr]^-{u} & \cong &
   \EMod\rlap{$(M,[0,1])$}\ar@/^2ex/[ll]^-{v}
}$$

\noindent then $X$ carries and $\Exp$-algebra structure given
by composition in:
$$\xymatrix@R-.5pc{
\Exp(X)\ar[d]_{\Exp(u)}^{\cong}\ar[rr]^-{\beta_X} & & X \\
\Exp\Big(\EMod(M,[0,1])\Big)\ar[rr]_-{\alpha_{M}} & & 
   \EMod(M,[0,1])\ar[u]_{v}^{\cong}
}$$

\noindent where $\beta_{X}(h) =
v\big(\lamin{m}{M}{h(\lamin{x}{X}{u(x)(m)})}\big)$.
\end{lemma}

\begin{myproof}
Assume $h\in\Exp(X) = \EMod([0,1]^{X},[0,1])$. Then:
$$\begin{array}[b]{rcl}
\big(v \after \alpha_{M} \after \Exp(u)\big)(h)
& = &
v\Big(\alpha_{M}\big(\Exp(u)(h)\big)\Big) \\
& = &
v\Big(\lamin{m}{M}{\Exp(u)(h)\big(\lam{k}{k(m)}\big)}\Big) \\
& = &
v\Big(\lamin{m}{M}{h\big((\lam{k}{k(m)}) \after u\big)}\Big) \\
& = &
v\Big(\lamin{m}{M}{h\big(\lamin{x}{X}{u(x)(m)}\big)}\Big).
\end{array}\eqno{\QEDbox}$$

\auxproof{
Explicit check that $\beta = \beta_{M}$ is an algebra:
$$\begin{array}{rcl}
\big(\beta \after \eta)(x)
& = &
\beta(\eta(x)) \\
& = &
v\big(\lamin{m}{M}{\eta(x)(\lamin{x}{X}{u(x)(m)})}\big) \\
& = &
v\big(\lamin{m}{M}{u(x)(m)}\big) \\
& = &
v(u(x)) \\
& = &
x \\
\big(\beta \after \Exp(\beta)\big)(H)
& = &
\beta\big(\Exp(\beta)(H)\big) \\
& = &
v\big(\lamin{m}{M}{\Exp(\beta)(H)(\lamin{x}{X}{u(x)(m)})}\big) \\
& = &
v\Big(\lamin{m}{M}{H\big((\lamin{x}{X}{u(x)(m)}) \after \beta\big)}\Big) \\
& = &
v\Big(\lamin{m}{M}{H\big((\lamin{k}{\Exp(X)}{u(\beta(k))(m)}}\Big) \\
& = &
v\Big(\lamin{m}{M}{H\big((\lamin{k}{\Exp(X)}
   {u(v\big(\lamin{m}{M}{k(\lamin{x}{X}{u(x)(m)})}\big))(m)}}\Big) \\
& = &
v\Big(\lamin{m}{M}{H(\lam{k}{k(\lamin{x}{X}{u(x)(m)})})}\Big) \\
& = &
v\Big(\lamin{m}{M}{\mu(H)(\lamin{x}{X}{u(x)(m)})}\Big) \\
& = &
\beta(\mu(H)) \\
& = &
\big(\beta \after \mu\big)(H).
\end{array}$$
}
\end{myproof}

We apply this result twice, in order to obtain and $\Exp$-algebra
structure on the unit interval and on density matrices on a Hilbert
space.

\begin{lemma}
\label{UnitIntExpAlgLem}
The unit interval $[0,1]$ carries an $\Exp$-algebra structure:
$$\xymatrix{
\Exp([0,1])\ar[r]^-{\cc} & [0,1]
& \quad\mbox{by}\quad &
h\ar@{|->}[r] & h(\idmap[{[0,1]}]).
}$$
\end{lemma}

\begin{myproof}
We apply Lemma~\ref{EModExpAlgTranslateLem} for the product effect
module $M = [0,1]^{2}$, since there is an isomorphism (of sets):
$$\xymatrix{
[0,1]\ar@/^2ex/[rr]^-{u} & \cong &
   \EMod\rlap{$([0,1]^{2},[0,1])$}\ar@/^2ex/[ll]^-{v}
}$$

\noindent given by:
$$\begin{array}{rclcrcl}
u(r)(s,t)
& = &
r\cdot s + (1-r)\cdot t
& \quad\mbox{and}\quad &
v(f)
& = &
f(1,0).
\end{array}$$

\noindent This gives an explicit description of the isomorphism
$[0,1]\cong\Exp(2)$ in Lemma~\ref{Exp012Lem}.

\auxproof{
First, $u(r)$ is a map of effect modules:
$$\begin{array}{rcl}
u(r)(0,0)
& = &
r\cdot 0 + (1-r)\cdot 0 \\
& = & 
0 \\
u(r)\big((s_{1},t_{1}) \ovee (s_{2},t_{2})\big)
& = &
u(r)\big(s_{1}+s_{2}, t_{1}+t_{2}\big) \\
& = &
r\cdot (s_{1}+s_{2}) + (1-r)\cdot (t_{1}+t_{2}) \\
& = &
r\cdot s_{1} + (1-r)\cdot t_{1} + r\cdot s_{2} + (1-r)\cdot t_{2} \\
& = &
u(r)(s_{1},t_{1}) \ovee u(r)(s_{2},t_{2}) \\
u(r)\big(a\cdot (s,t)\big)
& = &
u(r)\big(a\cdot s, a\cdot t\big) \\
& = &
r\cdot (a\cdot s) + (1-r)\cdot (a\cdot t) \\
& = &
a\cdot (r\cdot s) + a\cdot ((1-r)\cdot t) \\
& = &
a\cdot u(r)(s,t).
\end{array}$$

\noindent Finally, $u$ and $v$ are each others inverses:
$$\begin{array}{rcl}
\big(v \after u\big)(r)
& = &
v(u(r)) \\
& = &
u(r)(1,0) \\
& = &
r \cdot 1 + (1-r)\cdot 0 \\
& = &
r \\
\big(u \after v)(f)(s,t)
& = &
u(v(f))(s,t) \\
& = &
v(f)\cdot s + (1-v(f))\cdot t \\
& = &
f(1,0)\cdot s + f((1,0)^{\perp})\cdot t \\
& = &
f(s\cdot (1,0)) + f(t\cdot (0,1)) \\
& = &
f((s,0) \ovee (0,t)) \\
& = &
f(s,t).
\end{array}$$
}

\noindent Following Lemma~\ref{EModExpAlgTranslateLem}, the resulting
map $\cc \colon \Exp([0,1]) \rightarrow [0,1]$ is:
$$\begin{array}[b]{rcl}
\cc(h)
& = &
v\big(\lamin{(s,t)}{[0,1]^{2}}{h(\lamin{r}{[0,1]}{u(r)(s,t)})}\big) \\
& = &
h(\lamin{r}{[0,1]}{u(r)(1,0)}) \\
& = &
h(\lamin{r}{[0,1]}{1\cdot r}) \\
& = &
h(\idmap).
\end{array}\eqno{\QEDbox}$$

\auxproof{ 
It is easy to see that the evaluation-at-identity map
$\cc\colon\Exp([0,1]) \rightarrow [0,1]$ is an algebra. We explicitly
check the details:
$$\begin{array}{ccc}
\begin{array}{rcl}
\big(\cc \after \eta\big)(x)
& = &
\cc\big(\eta(x)\big) \\
& = &
\eta(x)(\idmap) \\
& = &
\idmap(x) \\
& = &
x
\end{array}
& \qquad &
\begin{array}{rcl}
\big(\cc \after \Exp(\cc)\big)(H)
& = &
\cc\big(\Exp(\cc)(H)\big) \\
& = &
\Exp(\cc)(H)(\idmap) \\
& = &
H(\idmap \after \cc) \\
& = &
H\big(\lamin{k}{\Exp([0,1])}{k(\idmap)}\big) \\
& = &
\mu(H)(\idmap) \\
& = &
\cc\big(\mu(H)\big) \\
& = &
\big(\cc \after \mu\big)(H).
\end{array}
\end{array}$$
}
\end{myproof}

Since Eilenberg-Moore algebras are closed under products, there is also
an $\Exp$-algebra on $[0,1]^{2}$. Explicitly, it is:
$$\xymatrix{
\Exp([0,1]^{2})\ar[r] & [0,1]^{2}
& \mbox{by} &
h\ar@{|->}[r] & \tuple{h(\pi_{1}), h(\pi_{2})},
}$$

\noindent where $h\colon [0,1]^{([0,1]^{2})} \rightarrow [0,1]$ is a
map of effect modules and the $\pi_{i}\colon [0,1]^{2} \rightarrow
[0,1]$ are the two projections. This extends easily to $[0,1]^{n}$,
and even to the exponent space $[0,1]^{A}$, for an arbitrary set
$A$. Recall that this space $[0,1]^{A}$ is compact (by Tychonoff's
theorem) and convex.  The algebra structure is:
$$\xymatrix@C-.5pc{
\Exp([0,1]^{A})\ar[r] & [0,1]^{A}
& \mbox{namely} &
h\ar@{|->}[r] & \lamin{a}{A}{h\big(\lamin{f}{[0,1]^{A}}{f(a)}\big)}.
}$$

\auxproof{ 
Let's write this map as $\cc\colon\Exp([0,1]^{A}) \rightarrow [0,1]^{A}$.
We explicitly check that it is an algebra:
$$\begin{array}{ccc}
\begin{array}{rcl}
\big(\cc \after \eta\big)(f)(a)
& = &
\cc\big(\eta(f)\big)(a) \\
& = &
\eta(f)(\lam{g}{g(a)}) \\
& = &
f(a) \\
& = &
\idmap(f)(a)
\end{array}
& \quad &
\begin{array}{rcl}
\big(\cc \after \Exp(\cc)\big)(F)(a)
& = &
\cc\big(\Exp(\cc)(F)\big)(a) \\
& = &
\Exp(\cc)(F)(\lam{g}{g(a)}) \\
& = &
F((\lam{g}{g(a)}) \after \cc) \\
& = &
F\big(\lam{h}{\cc(h)(a)}\big) \\
& = &
F\big(\lam{h}{h(\lam{f}{f(a)})}\big) \\
& = &
\mu(F)(\lam{f}{f(a)}) \\
& = &
\cc\big(\mu(F)\big)(a) \\
& = &
\big(\cc \after \mu\big)(F)(a).
\end{array}
\end{array}$$
}

The next application of Lemma~\ref{EModExpAlgTranslateLem} involves
the set $\DM(H)$ of density matrices on a finite-dimensional Hilbert
space $H$. We recall from~\cite{JacobsM??} that there is a
``Hilbert-Schmidt'' isomorphism:
$$\xymatrix{
\DM(H)\ar[rr]^-{\hs}_-{\cong} & & \EMod\big(\Ef(H),[0,1]\big)
}$$

\noindent where $\Ef(H) = \set{A\colon H\rightarrow A}{0 \leq A \leq
  I}$ is the effect module of effects on $H$. This map $\hs$ is given
by $\hs(D)(A) = \tr(DA)$.

On $h\in\Exp(\DM(H)) = \EMod([0,1]^{\DM(H)}, [0,1])$, the algebra
$\beta_{\DM(H)}\colon \Exp(\DM(H)) \rightarrow \DM(H)$ is given by:
$$\begin{array}[b]{rcl}
\beta_{\DM(H)}(h)
& = &
\hs^{-1}\Big(\lamin{A}{\Ef(H)}{h(\lamin{D}{\DM(H)}{\hs(D)(A)})}\Big) \\
& = &
\hs^{-1}\Big(\lamin{A}{\Ef(H)}{h(\lamin{D}{\DM(H)}{\tr(DA)})}\Big) \\
& = &
\hs^{-1}\Big(\lamin{A}{\Ef(H)}{h(\tr(A-))}\Big) \\
\end{array}\eqno{\QEDbox}$$
}

We can also form a functor $\Alg(\Exp)\rightarrow\EMod\op$, in the
reverse direction in~\eqref{ExpComparisonDiag}, by ``homming'' into
the unit interval $[0,1]$. Recall that this interval carries an
$\Exp$-algebra, identified in Lemma~\ref{UnitIntExpAlgLem} as
evaluation-at-identity $\ev_{\idmap}$.  For an algebra
$\alpha\colon\Exp(X)\rightarrow X$ we know that the algebra
homomorphisms $X\rightarrow [0,1]$ are precisely the affine continuous
maps $X\rightarrow [0,1]$, by Lemma~\ref{ExpAlgHomLem}.  We shall be a
bit sloppy in our notation and write the homset
$\Alg(\Exp)(\alpha,\ev_{\idmap}) = \set{q\colon X\rightarrow [0,1]}{ q
  \after \alpha = \ev_{\idmap} \after \Exp(q)}$ of algebras map in
various ways, namely as:
$$\left\{\begin{array}{ll}
\Alg(\Exp)(\alpha,[0,1]) \quad &
   \mbox{leaving the algebra structure $\ev_{\idmap}$ on $[0,1]$
         implicit,} \\
\Alg(\Exp)(X,[0,1]) &
   \mbox{also leaving the algebra structure $\alpha$ on $X$
         implicit,} \\
\mathcal{A}(X,[0,1]) &
   \mbox{as set of affine continuous functions, 
         via Lemma~\ref{ExpAlgHomLem}.}
\end{array}\right.$$

\begin{proposition}
\label{EModAlgExpAdjProp}
The states functor $S_{(-)} = \EMod(-,[0,1]) \colon \EMod\op \rightarrow
\Alg(\Exp)$ from~\eqref{ExpComparisonDiag} has a left adjoint, also
given by ``homming into $[0,1]$'':
$$\xymatrix@C+1pc{
\EMod\op\ar@/^2ex/[rr]^-{S_{(-)} = \EMod(-,[0,1])} & \top &
   \Alg(\Exp)\ar@/^2ex/[ll]^-{\Alg(\Exp)(-,[0,1])}
}$$
\end{proposition}

\begin{myproof}
Assume an $\Exp$-algebra $\alpha\colon\Exp(X)\rightarrow X$.  We
should first check that the set of affine continuous functions is a
sub-effect module: $\Alg(\Exp)(\alpha,[0,1]) = \mathcal{A}(X,[0,1])
\hookrightarrow [0,1]^{X}$. The top and bottom elements $1 =
\lam{y}{1}$ and $0 = \lam{y}{0}$ are clearly in
$\mathcal{A}(X,[0,1])$. Also, $\mathcal{A}(X,[0,1])$ is closed under
(partial) sums $\ovee$ and scalar multiplication with $r\in[0,1]$.
Next, if we have a map of algebras $g\colon X\rightarrow Y$, from
$\smash{\Exp(X)\stackrel{\alpha}{\rightarrow} X}$ to
$\smash{\Exp(Y)\stackrel{\beta}{\rightarrow} Y}$. Then we get a map of
effect modules $g^{*} = (-) \after g \colon \mathcal{A}(X,[0,1])
\rightarrow \mathcal{A}(Y,[0,1])$. This is easy because $g$ is itself
affine and continuous, by Lemma~\ref{ExpAlgHomLem}.

We come to the adjunction $\Alg(\Exp)(-,[0,1]) \dashv
\EMod(-,[0,1])$. For $M\in\EMod$ and
$\smash{(\Exp(X)\stackrel{\alpha}{\rightarrow} X)}\in\Alg(\Exp)$ there
is a bijective correspondence:
$$\begin{prooftree}
\xymatrix@R-.8pc{
\Exp(X)\ar[d]_{\alpha}\ar[r]^-{\Exp(f)} & 
   \Exp(S_{M})\ar[d]^{\alpha_M} \\
X\ar[r]^-{f} & S_{M}
}
\Justifies
\xymatrix{M\ar[r]_-{g} & \Alg(\Exp)(\alpha,[0,1])
   \rlap{$\;=\mathcal{A}(X,[0,1]$}}
\end{prooftree}$$

\noindent We proceed as follows.
\begin{itemize}
\item Given an algebra map $f\colon X\rightarrow S_{M} =
  \EMod(M,[0,1])$ as indicated, define $\overline{f} \colon M
  \rightarrow \mathcal{A}(X,[0,1])$ by $\overline{f}(y)(x) =
  f(x)(y)$. We leave it to the reader to check that $\overline{f}$ is
  a map of effect modules, but we do verify that $\overline{f}(y)$ is
  an algebra map $X \rightarrow [0,1]$; so for $h\in\Exp(X)$,
$$\begin{array}{rcl}
\big(\ev_{\idmap{}} \after \Exp(\overline{f}(y))\big)(h)
& = &
\Exp(\overline{f}(y))(h)(\idmap{}) \\
& = &
h(\idmap{} \after \overline{f}(y)) \\
& = &
h(\overline{f}(y)) \\
& = &
h\big(\lam{x}{f(x)(y)}\big) \\ 
& = &
h\big((\lam{k}{k(y)}) \after f\big) \\
& = &
\Exp(f)(h)\big(\lam{k}{k(y)}\big) \\
& = &
\alpha_{M}\big(\Exp(f)(h)\big)(y) \\
& = &
f\big(\alpha(h)\big)(y) \qquad \mbox{since $f$ is
   an algebra map} \\
& = &
\big(\overline{f}(y) \after \alpha\big)(h).
\end{array}$$

\auxproof{
Verifications that $\overline{f}\colon M \rightarrow \mathcal{A}(X,[0,1])$
is a map of effect modules.
$$\begin{array}{rcl}
\overline{f}(0)
& = &
\lam{x}{\overline{f}(0)(x)} \\
& = &
\lam{x}{f(x)(0)} \\
& = &
\lam{x}{0} \\
& = &
0 \\
\overline{f}(y_{1}\ovee y_{2}) 
& = &
\lam{x}{\overline{f}(y_{1}\ovee y_{2})(x)} \\
& = &
\lam{x}{f(x)(y_{1}\ovee y_{2})} \\
& = &
\lam{x}{f(x)(y_{1}) \ovee f(x)(y_{2})} \\
& = &
\lam{x}{\overline{f}(y_{1})(x) \ovee \overline{f}(y_{2})(x)} \\
& = &
\overline{f}(y_{1}) \ovee \overline{f}(y_{2}) \\
\overline{f}(r\cdot y)
& = &
\lam{x}{\overline{f}(r\cdot y)(x)} \\
& = &
\lam{x}{f(x)(r\cdot y)} \\
& = &
\lam{x}{r\cdot f(x)(y)} \\
& = &
\lam{x}{r\cdot \overline{f}(y)(x)} \\
& = &
r\cdot \overline{f}(y).
\end{array}$$
}

\item Now assume we have a map of effect modules $g\colon M
  \rightarrow \Alg(\Exp)(\alpha,[0,1])=\mathcal{A}(X,[0,1])$. We turn
  it into a map of algebras $\overline{g}\colon X\rightarrow S_{M}$,
  again by twisting arguments: $\overline{g}(x)(y) = g(y)(x)$. Via
  calculations as above one checks that $\overline{g}$ is a map of
  algebras.

\auxproof{
For $h\in \Exp(X)$ and $y\in M$,
$$\begin{array}{rcl}
\big(\alpha_{M} \after \Exp(\overline{g})\big)(h)(y)
& = &
\alpha_{M}\big(\Exp(\overline{g})(h)\big)(y) \\
& = &
\Exp(\overline{g})(h)\big(\lam{k}{k(y)}\big) \\
& = &
h\big((\lam{k}{k(y)}) \after \overline{g}\big) \\
& = &
h\big(\lamin{x}{X}{\overline{g}(x)(y)}\big) \\
& = &
h\big(\lamin{x}{X}{g(y)(x)}\big) \\
& = &
h\big(g(y)\big) \\
& = &
h\big(\idmap \after g(y)\big) \\
& = &
\Exp(g(y))(h)(\idmap) \\
& = &
\ev_{\idmap}\big(\Exp(g(y))(h)\big) \\
& = &
g(y)\big(\alpha(h)\big) \qquad
    \mbox{since $g(y)$ is a map of algebras} \\
& = &
\overline{g}(\alpha(h))(y) \\
& = &
\big(\overline{g} \after \alpha\big)(h)(y).
\end{array}$$
}
\end{itemize}

\noindent Clearly $\smash{\overline{\overline{f}} = f}$
and $\smash{\overline{\overline{g}} = g}$. \QED
\end{myproof}

Let's write the unit and counit of this adjunction as $\eta^{\dashv}$
and $\varepsilon^{\dashv}$, in order make a distinction with the unit
$\eta$ of the monad $\Exp$, see below. The unit and counit are maps:
$$\begin{array}{rcl}
\xymatrix{
X\ar[rr]^-{\eta^{\dashv}} & & \EMod\Big(\Alg(\Exp)(\alpha,[0,1]),[0,1]\Big)
}
& \qquad &
\mbox{in }\Alg(\Exp) \\[-.3em]
\xymatrix{
\Alg(\Exp)\Big(\EMod(M,[0,1]),[0,1]\Big)\ar[rr]_-{\varepsilon^{\dashv}} 
   & & M
}
& & 
\mbox{in }\EMod\op
\end{array}$$

\noindent both given by point evaluation:
\begin{equation}
\label{EModAlgUnitCounitEqn}
\begin{array}{rclcrcl}
\eta^{\dashv}(x)
& = &
\lamin{f}{\mathcal{A}(X,[0,1])}{f(x)}
& \quad &
\varepsilon^{\dashv}(y)
& = &
\lamin{g}{\EMod(M,[0,1])}{g(y)}
\end{array}
\end{equation}

\noindent The unit of this adjunction is related to the unit
of the monad $\Exp$, written explicitly as $\eta^{\Exp}$ in the
following way.
$$\xymatrix@R-2pc{
& & \EMod\big(\mathcal{A}(X,[0,1]),[0,1]\big) \\
X\ar[urr]^-{\eta^\dashv}\ar[drr]_-{\eta^{\Exp}} & & \\
& & \Exp(X)\rlap{$\;=\EMod\big([0,1]^{X},[0,1]\big)$}\ar[uu]
}$$

\begin{lemma}
\label{EModAlgExpAdjUnitLem}
Consider the unit $\eta^{\dashv}$ in~\eqref{EModAlgUnitCounitEqn} of
the adjunction from Proposition~\ref{EModAlgExpAdjProp}, at an algebra
$\smash{\Exp(X)\stackrel{\alpha}{\rightarrow} X}$.
\begin{enumerate}
\item This unit is injective if and only if $X$ is observable.

\item In fact, it is an isomorphism if and only if $X$ is observable.
\end{enumerate}

\noindent Hence the adjunction $\Alg(\Exp) \leftrightarrows \EMod\op$
from Proposition~\ref{EModAlgExpAdjProp} restricts to a coreflection
$\Algobs(\Exp) \leftrightarrows \EMod\op$.
\end{lemma}

\begin{myproof}
The first statement holds by definition of `observable'. So for the
second statement it suffices to assume that $X$ is observable and show
that $\eta^{\dashv}\colon X \rightarrowtail
\EMod\big(\mathcal{A}(X,[0,1]),[0,1]\big)$ is surjective. This unit is
by construction affine and continuous. Hence its image in the space
$\EMod\big(\mathcal{A}(X,[0,1]),[0,1]\big)$ is compact, and thus
closed. 

We are done if $\eta^{\dashv}$ is dense.  Thus we assume a non-empty
open set $U\subseteq \EMod\big(\mathcal{A}(X,[0,1]),[0,1]\big)$ and
need to prove that there is an $y\in X$ with $\eta^{\dashv}(y)\in
U$. By Proposition~\ref{EModExpAlgProp} we may assume $U =
\square_{s_1}(q_{1}) \cap \cdots \cap \square_{s_n}(q_{n})$, for
$q_{i}\in \mathcal{A}(X,[0,1])$ and
$s_{i}\in[0,1]\cap\mathbb{Q}$. Thus we may assume a map of effect
modules $h\colon\mathcal{A}(X,[0,1]) \rightarrow [0,1]$ inhabiting all
these $\square$'s. Hence $h(q_{i}) < s_{i}$.

Since $\iota\colon \mathcal{A}(X,[0,1]) \hookrightarrow [0,1]^{X}$ is
a sub-effect module, by extension, see Proposition~\ref{HBEModProp},
we get a map of effect modules $h'\colon [0,1]^{X}\rightarrow [0,1]$ with
$h' \after \iota = h$. Thus, we can take the inverse image of the open
set $U$ along the continuous map:
$$\xymatrix{
\Exp(X) = \EMod\big([0,1]^{X},[0,1]\big)\ar[rr]^-{(-) \after \iota} & &
   \EMod\big(\mathcal{A}(X,[0,1]),[0,1]\big)
}$$

\noindent The resulting open set is:
$$\begin{array}{rcccl}
V 
& \smash{\stackrel{\textrm{def}}{=}} &
\big((-)\after\iota\big)^{-1}(U)
& = &
\setin{k}{\Exp(X)}{k\after\iota\in U} \\
& & & = &
\setin{k}{\Exp(X)}{\all{i}{k(\iota(q_{i})) < s_{i}}}.
\end{array}$$

\noindent This subset $V\subseteq\Exp(X)$ contains $h'$ and is thus
non-empty. Since $\sigma\colon\Dst(X) \rightarrowtail \Exp(X)$ is
dense, by Proposition~\ref{DstExpDenseProp}, there is a distribution
$\varphi = (\sum_{j}r_{j}x_{j})\in\Dst(X)$ with $\sigma(\varphi) \in
V$. We take $x = \sum_{j}r_{j}x_{j} \in X$ to be the interpretation of
$\varphi$ in $X$, using that $X$ is convex. We claim
$\eta^{\dashv}(x)\in U$.  Indeed, $\eta^{\dashv}(x) \in
\square_{s}(q_{i})$, for each $i$, since:
$$\begin{array}[b]{rcll}
\eta^{\dashv}(x)(q_{i})
\hspace*{\arraycolsep} = \hspace*{\arraycolsep}
q_{i}(x)
& = &
q_{i}(\sum_{j}r_{j}x_{j}) \\
& = &
\sum_{j}r_{j}\cdot q_{i}(x_{j}) \qquad
   & \mbox{since $q_{i}\colon X\rightarrow [0,1]$ is affine} \\
& = &
\sigma(\varphi)(\iota(q_{i})) \\
& < &
s_{i} & \mbox{since $\sigma(\varphi)\in V$.}
\end{array}\eqno{\QEDbox}$$
\end{myproof}

We turn to the counit~\eqref{EModAlgUnitCounitEqn} of the adjunction
in Proposition~\ref{EModAlgExpAdjProp}.

\begin{lemma}
\label{EModAlgExpAdjCounitLem}
For an effect module $M$ consider the counit $\varepsilon^{\dashv}
\colon M \rightarrow \mathcal{A}(S_{M}, [0,1])$, where, as before,
$S_{M} = \EMod(M,[0,1])$ is the convex compact Hausdorff space of
states.
\begin{enumerate}
\item The effect module $\mathcal{A}(S_{M},[0,1])$ is ``Banach'',
\textit{i.e.}~complete.

\item The counit map $\varepsilon^{\dashv}$ is a dense embedding of
  $M$ into this Banach effect module $\mathcal{A}(S_{M},[0,1])$.

\item Hence it is an isomorphism if and only if $M$ is a Banach effect
  module.
\end{enumerate}
\end{lemma}

\begin{myproof}
Completeness of $\mathcal{A}(S_{M},[0,1])$ is inherited from $[0,1]$,
since its norm is the supremum norm, like in Example~\ref{AEModEx}.

For the second point we use the corresponding result for order unit
spaces, via the equivalence $\smash{\pth\colon \AEMod
  \xrightarrow{\simeq}\OUS}$ from
Proposition~\ref{ArchimideanEquivProp}. If $(V,u)$ is an order unit
space then it is well known (see~\cite{AsimovE80}) that the evaluation
map $\theta\colon V\to\mathcal{A}(S,\reals)$ is a dense
embedding. Here $S=\OUS(V,\reals)$ is the state space of $V$. However
if we take $V$ to be the totalization $\pth(M)$ of $M$, then
$\theta$ is precisely $\pth(\varepsilon^{\dashv})$, since:
$$\begin{array}[t]{rcl}
\smash{\xymatrix@C+1pc{
\pth(M) = V\ar@{ >->}[r]_-{\mbox{\small dense}}^-{\theta} &
  \mathcal{A}(S,\reals)}}
& = & 
\mathcal{A}\big(\OUS(V,\reals),\reals\big) \\
& \cong & 
\mathcal{A}\big(\EMod(M,[0,1]),\reals\big) \\
& \cong & 
\pth\Big(\mathcal{A}\big(\EMod(M,[0,1]\big),[0,1])\Big).
\end{array}$$
and both $\theta$ and $\pth(\varepsilon^{\dashv})$ are the evaluation map.

For the third point, one direction is easy: if the counit is an
isomorphism, then $M$ is isomorphic to the complete effect module
$\mathcal{A}\big(S_{M},[0,1])$, and thus complete itself. In the other
direction, denseness of $M\rightarrowtail
\mathcal{A}\big(S_{M},[0,1])$ means that each $h\in
\mathcal{A}\big(S_{M},[0,1])$ can be expressed as limit $h =
\lim_{n}\varepsilon^{\dashv}(x_{n})$ of elements $x_{n}\in M$. But if
$M$ is complete, there is already a limit $x=\lim_{n}x_{n}\in
M$. Hence $\varepsilon^{\dashv}(x) = h$, making $\varepsilon^{\dashv}$
an isomorphism. \QED
\end{myproof}

Combining lemmas~\ref{EModAlgExpAdjCounitLem}
and~\ref{EModAlgExpAdjUnitLem} gives us the main result of this paper.

\begin{theorem}
\label{EModAlgExpDualityThm}
The adjunction $\Alg(\Exp) \leftrightarrows \EMod\op$ from
Proposition~\ref{EModAlgExpAdjProp} restricts to a duality
$\Algobs(\Exp) \simeq \BEMod\op$ between observable $\Exp$-algebras
and Banach effect modules. In combination with
Theorem~\ref{AlgobsCCHobsThm} we obtain:
$$\begin{array}{rcccl}
\CCHobs & \cong & \Algobs(\Exp) & \simeq & \BEMod\op.
\end{array}\eqno{\QEDbox}$$
\end{theorem}

This result can be seen as a probabilistic version of fundamental
results of Manes (Theorem~\ref{ManesThm}) and Gelfand
(Theorem~\ref{GelfandThm}).

\section{A new formulation of Gleason's theorem}\label{GleasonSec}

Gleason's theorem in quantum mechanics says that every state on a
Hilbert space of dimension three or greater corresponds to a density
matrix~\cite{Gleason57}. In this section we introduce a reformulation
of Gleason's theorem, and prove the equivalence via Banach effect
modules (esp.\ Lemma~\ref{EModAlgExpAdjCounitLem}). This reformulation
says that effects are the free effect module on projections. In
formulas: $\Ef(\H) \cong [0,1]\otimes \pr(\H)$, for a Hilbert space
$\H$.

Gleason's theorem is not easy to prove (see
\textit{e.g.}~\cite{Dvurecenskij92}).  Even proofs using elementary
methods are quite involved~\cite{CookeKM85}.  A state on a Hilbert
space $\H$ is a certain probability distribution on the projections
$\pr(\H)$ of $\H$. These projections $\pr(\H)$ form an orthomodular
lattice, and thus an effect
algebra~\cite{DvurecenskijP00,JacobsM12a}. In our current context
these are exactly the effect algebra maps $\pr(\H)\to[0,1]$. So
Gleason's (original) theorem states:
\begin{equation}
\label{GleasonEAEqn}
\begin{array}{rcl}
\EA\big(\pr(\H),[0,1]\big) & \cong & \DM(\H).
\end{array}
\end{equation}

\noindent This isomorphism, from right to left, sends a density matrix
$M$ to the map $p\mapsto \tr(Mp)$---where $\tr$ is the trace map
acting on operators.

Recall that $\Ef(\H)$ is the set of positive operators on $\H$ below
the identity. It is a Banach effect module. One can also consider the
effect module maps $\Ef(\H)\to [0,1]$. For these maps there is a
``lightweight'' version of Gleason's theorem:
\begin{equation}
\label{GleasonEModEqn}
\begin{array}{rcl}
\EMod\big(\Ef(\H),[0,1]\big) & \cong & \DM(\H).
\end{array}
\end{equation}

\noindent This isomorphism involves the same trace computation
as~\eqref{GleasonEAEqn}. This statement is significantly easier to
prove than Gleason's theorem itself, see~\cite{Busch03}.

Because Gleason's original theorem~\eqref{GleasonEAEqn} is so much
harder to prove than the lightweight version~\eqref{GleasonEModEqn}
one could wonder what Gleason's theorem states that Gleason light
doesn't. In Theorem~\ref{T:Gleason} we will show that the difference
amounts exactly to the statement:
\begin{equation}
\label{GleasonFreeEModEqn}
\begin{array}{rcl}
[0,1]\otimes \pr(\H) & \cong & \Ef(\H),
\end{array}
\end{equation}

\noindent where $\otimes$ is the tensor of effect algebras
(see~\cite{JacobsM12a}). A general result,
see~\cite[VII,\S4]{MacLane71}, says that the tensor product
$[0,1]\otimes \pr(\H)$ is the free effect module on $\pr(\H)$.

The following table gives an overview of the various formulations
of Gleason's theorem.
\begin{center}
\begin{tabular}{c|c|c}
\textbf{Description} & \textbf{Formulation} 
   & \quad\textbf{Label}\quad
   \vrule height5mm depth3mm width0mm \\
\hline\hline 
\begin{tabular}{c} original Gleason, \\[-.5em]
for projections \end{tabular} \quad &
   $\EA\big(\pr(\H),[0,1]\big) \cong \DM(\H)$ &
    \eqref{GleasonEAEqn} \ 
   \vrule height5mm depth3mm width0mm \\
\begin{tabular}{c} lightweight version, \\[-.5em]
for effects \end{tabular} &
   $\;\EMod\big(\Ef(\H),[0,1]\big) \cong \DM(\H)\;$ &
   \eqref{GleasonEModEqn} 
   \vrule height5mm depth3mm width0mm \\
\quad\begin{tabular}{c} effects as free \\[-.5em]
module on projections \end{tabular}\quad &
   $[0,1]\otimes \pr(\H) \cong \Ef(\H)$ &
   \eqref{GleasonFreeEModEqn} 
   \vrule height5mm depth3mm width0mm \\
\end{tabular}
\end{center}

\noindent In this section we shall prove~\eqref{GleasonEAEqn}
$\Longleftrightarrow$ \eqref{GleasonFreeEModEqn}, in presence
of~\eqref{GleasonEModEqn}, see Theorem~\ref{T:Gleason}.
Since~\eqref{GleasonEAEqn} is true, for dimension $\geq 3$, the same
then holds for~\eqref{GleasonFreeEModEqn}.

We first prove a general result based on the duality from the previous
section.  There we used the shorthand $S_{M}$ for the algebra of
states $\EMod(M,[0,1])$. We now extend this notation to effect
algebras and write $S_{D} = \EA(D,[0,1])$, where $D$ is an effect
algebra.  We recall from Section~\ref{EModSec} that $S_D$ is a convex
set. We will topologize it via the weakest topology that makes all
point evaluations continuous.

Since the tensor product $[0,1]\otimes D$ of effect algebras is the
free effect module on $D$ it follows that there is an isomorphism:
\begin{equation}
\label{FreeEModMapEqn}
\vcenter{\xymatrix@R-1.5pc{\hspace*{-.5em}
\EA\big(D,[0,1]\big)\ar[r]_-{\cong}^-{\widehat{(-)}} &
   \EMod\big([0,1]\otimes D,[0,1]\big) & \llap{with}\quad
   \widehat{f}(s\sotimes x) = s\cdot f(x) \\
S_{D}\ar@{=}[u] & S_{[0,1]\otimes D}\ar@{=}[u]
}}
\end{equation}

\begin{lemma}
\label{FreeEModMapLem}
The mapping $\smash{\widehat{(-)}\colon S_{D} \conglongrightarrow 
S_{[0,1]\otimes D}}$ in~\eqref{FreeEModMapEqn} is an
affine homeomorphism.
\end{lemma}

\begin{myproof}
We only show that $\smash{\widehat{(-)}}$ is a homeomorphism. For an
arbitrary element $\ovee_{i}\, r_i\sotimes x_i\in [0,1]\otimes D$ we
have in $[0,1]$,
$$\begin{array}{rcl}
\smash{\widehat{f}\big(\ovee_{i}\, r_{i}\sotimes x_{i}\big)}
& = &
\sum_{i} r_{i} \cdot f(x_{i}).
\end{array}$$

\noindent Since the maps $f\mapsto f(x_i)$ are continuous by
definition of the topology on $S_D$, and since addition and
multiplication on $[0,1]$ are continuous, it follows that
$\smash{f\mapsto \widehat{f}\big(\ovee_{i}\, r_i\sotimes x_{i}\big)}$
is continuous. Hence by definition of the topology on $S_{[0,1]\otimes
  D}$ we see that the mapping $\smash{\widehat{(-)}}$ is continuous.

Similarly, the inverse, say written as $\smash{\widetilde{(-)}\colon
  S_{[0,1] \otimes D} \rightarrow S_{D}}$, is continuous. It is given
by $\widetilde{k}(x) = k(1\sotimes x)$. Continuity again follows from
the definition of the topology on $S_{[0,1]\otimes D}$. \QED
\end{myproof}

\begin{lemma}
\label{L:Banachiso}
Suppose $f\colon D\to E$ is an effect algebra map between an effect
algebra $D$ and a Banach effect module $E$ such that the following
hold.
\begin{itemize}
\item The induced map $\widehat{f}\colon [0,1]\otimes D\to E$ is 
surjective---obtained like in~\eqref{FreeEModMapEqn} as
$\widehat{f}(s\sotimes x) = s\scalar f(x)$.

\item The ``precompose with $f$'' map $-\after f\colon S_E\to S_D$ is
  a homeomorphism.
\end{itemize}

\noindent The map $\smash{\widehat{f}}$ is then an isomorphism between
$[0,1]\otimes D$ and $E$.
\end{lemma}

\begin{myproof}
Using Lemma~\ref{EModAlgExpAdjCounitLem}, there are for the
Banach effect module $E$ and for the (free) effect module
$[0,1]\otimes D$, maps $\varepsilon_E$ and $\phi_D$ in:
$$\vcenter{\xymatrix{
E\ar[r]^-{\varepsilon_E}_-{\cong} & \mathcal{A}\big(S_E,[0,1]\big)
}}\qquad
\vcenter{\xymatrix{
[0,1]\otimes D
   \ar@{ >->}[rr]^-{\varepsilon_{[0,1]\otimes D}}_-{\mbox{\small dense}}
   \ar@/_2ex/[drr]_{\phi_{D}} & &
   \mathcal{A}\big(S_{[0,1]\otimes D}, [0,1]\big)
      \ar[d]_{\cong}^{h\mapsto h(\widehat{(-)})} \\
& & \mathcal{A}\big(S_{D}, [0,1]\big)
}}$$

\noindent The operation $\widehat{(-)}$ on the right is as
in~\eqref{FreeEModMapEqn}. We claim that the following diagram
commutes.
$$\xymatrix{
[0,1]\otimes D\ar@{ >->}[rr]_-{\mbox{\small dense}}^-{\phi_D}
   \ar@{->>}[d]_{\widehat{f}} & &
   \mathcal{A}\big(S_{D}, [0,1]\big)
      \ar[d]_{\cong}^{k \mapsto k(- \after f)} \\
E\ar[rr]_-{\varepsilon_E}^-{\cong} & & \mathcal{A}\big(S_E,[0,1]\big)
}$$

\noindent If this is indeed true the map $\smash{\widehat{f}}$ is an
embedding followed by two isomorphism and therefore injective (and
thus an isomorphism). To prove the claim, we assume $\ovee_{i}\,
r_{i}\sotimes x_{i}\in[0,1]\otimes D$ and $g\in S_E$ and compute first
the east-south direction:
$$\begin{array}[b]{rcl}
\Big((\lam{k}{k(-\after f)}) \after \phi_{D}\Big)
   \big(\ovee_{i}\, r_{i}\sotimes x_{i}\big)(g) 
& = &
\phi_{D}\big(\ovee_{i}\, r_{i}\sotimes x_{i}\big)(g \after f) \\
& = &
\varepsilon_{[0,1]\otimes D}\big(\ovee_{i}\, r_{i}\sotimes x_{i}\big)
   \big(\widehat{(g \after f)}\big) \\
& = &
\widehat{(g \after f)}\big(\ovee_{i}\, r_{i}\sotimes x_{i}\big) \\
& = &
\sum_{i} r_{i}\cdot g(f(x_{i})) \\
& = &
g\big(\ovee_{i}\, r_{i}\scalar f(x_{i})\big) 
   \qquad \rlap{since $g$ is affine} \\
& = &
g\big(\widehat{f}\big(\ovee_{i}\, r_{i}\sotimes x_{i}\big)\big) \\
& = &
\varepsilon_{E}\big(\widehat{f}\big(\ovee_{i}\, r_{i}\sotimes x_{i}\big)
   \big)(g) \\
& = &
\big(\varepsilon_{E} \after \widehat{f}\big)
   \big(\ovee_{i}\, r_{i}\sotimes x_{i}\big)(g). \qquad\qquad
\end{array}\eqno{\QEDbox}$$
\end{myproof}

As a consequence we obtain the isomorphism~\eqref{GleasonFreeEModEqn}.
We will show next that it is equivalent to Gleason's (original)
theorem.

\begin{theorem}
\label{T:Gleason}
\eqref{GleasonEAEqn} $\Longleftrightarrow$ \eqref{GleasonFreeEModEqn},
in presence of~\eqref{GleasonEModEqn}.

That is, using Gleason light~\eqref{GleasonEModEqn} the following
statements are equivalent.
\begin{description}
\item[\eqref{GleasonEAEqn}:] $\EA(\pr(\H),[0,1])\cong\DM(\H)$,
  \textit{i.e.}~Gleason's original theorem;

\item[\eqref{GleasonFreeEModEqn}:] The canonical map
  $[0,1]\otimes\pr(\H)\to\Ef(\H)$ is an isomorphism.
\end{description}
\end{theorem}

\begin{myproof}
Assuming $[0,1]\otimes\pr(\H)\conglongrightarrow\Ef(\H)$ we get
Gleason's theorem:
$$\begin{array}{rcll}
\EA\big(\pr(\H),[0,1]\big)
& \cong &
\EMod\big([0,1]\otimes\pr(\H),[0,1]\big) \quad
   & \mbox{by freeness} \\
& \cong &
\EMod\big(\Ef(\H),[0,1]\big) 
   & \mbox{by assumption} \\
& \cong &
\DM(\H) & \mbox{by Gleason light~\eqref{GleasonEModEqn}.}
\end{array}$$


In the other direction assume
$S_{\pr(\H)}=\EA(\pr(\H),[0,1])\cong\DM(\H)$. We apply the previous
lemma to the inclusion $f\colon \pr(\H) \hookrightarrow \Ef(\H)$. Then
indeed:
\begin{itemize}
\item the induced map $\smash{\widehat{f}\colon [0,1]\otimes \pr(\H)
  \rightarrow \Ef(H)}$ is surjective: each effect $A\in\Ef(\H)$ can be
  written as convex combination of projections $A =
  \sum_{i}r_{i}P_{i}$, via the spectral theorem.

\item the precomposition $- \after f \colon S_{\Ef(\H)} \rightarrow
  S_{\pr(\H)}$ is an isomorphism since:
$$\begin{array}{rcccl}
S_{\Ef(\H)} 
& \smash{\stackrel{\eqref{GleasonEModEqn}}{\cong}} &
\DM(\H)
& \smash{\stackrel{\eqref{GleasonEAEqn}}{\cong}} &
S_{\pr(\H)}.
\end{array}$$

\noindent Since both these isomorphisms involve the same trace
computation, this isomorphism is in fact the map induced by the
inclusion $f\colon \pr(\H)\hookrightarrow\Ef(\H)$.
\end{itemize}

\noindent Thus the conditions of Lemma~\ref{L:Banachiso} are met and
so $[0,1]\otimes\pr(\H)\cong\Ef(\H)$. \QED
\end{myproof}

\section{The expectation monad for program semantics}\label{SemanticsSec}

This paper uses the expectation monad $\Exp(X) = \EMod([0,1]^{X},
[0,1])$ in characterization and duality results for convex compact
Hausdorff spaces. Elements of $\Exp(X)$ are characterized as (finitely
additive) measures (see esp.~Theorem~\ref{FinAddMeasThm}). The way the
monad $\Exp$ is defined, via the adjunction $\Sets \leftrightarrows
\EMod\op$, is new. This approach deals effectively with the rather
subtle preservation properties for maps $h\in\Exp(X) =
\EMod([0,1]^{X}, [0,1])$, namely preservation of the structure of
effect modules (with non-expansiveness, and thus continuity, as
consequence, see~Lemma~\ref{AEModNonExpLem}).

Measures have been captured via monads before, first by
Giry~\cite{Giry82} following ideas of Lawvere. Such a description in
terms of monads is useful to provide semantics for probabilistic
programs~\cite{Kozen81,JonesP89,McIverM04,Panangaden09}.  The term
`expectation monad' seems to occur first in~\cite{RamseyP02}, where it
is formalized in Haskell. Such a formalization in a functional
language is only partial, because the relevant equations and
restrictions are omitted, so that there is not really a difference
with the continuation monad $X \mapsto [0,1]^{([0,1]^{X})}$. A
formalization of what is also called `expectation monad' in the
theorem prover Coq occurs in~\cite{AudebaudP09} and is more
informative. It involves maps $h\colon [0,1]^{X} \rightarrow [0,1]$
which are required to be monotone, continuous, linear (preserving
partial sum $\ovee$ and scalar multiplication) and compatible with
inverses---meaning $h(1-p) \leq 1 - h(p)$. This comes very close to
the notion of homomorphism of effect module that is used here, but
effect modules themselves are not mentioned
in~\cite{AudebaudP09}. This Coq formalization is used for instance in
the semantics of game-based programs for the certification of
cryptographic proofs in~\cite{BartheGZ09}
(see~\cite{ZanellaBeguelin10} for an overview of this line of work).
Finally, in~\cite{KeimelRS11} a monad is used of maps $h\colon
[0,1]^{X} \rightarrow [0,1]$ that are (Scott) continuous and
sublinear---\textit{i.e.}~$h(p\ovee q) \leq h(p) \ovee h(q)$, and
$h(r\cdot p) = r\cdot h(p)$.

The definition $\Exp(X) = \EMod([0,1]^{X}, [0,1])$ of the expectation
monad that is used here has good credentials to be the right
definition, because:
\begin{itemize}
\item The monad $\Exp$ arises in a systematic (not \textit{ad hoc})
  manner, namely via the composable
  adjunctions~\eqref{SetsConvEModDiag}.

\item The sets $\Exp(X)$ as defined here form a stable collection, in
  the sense that its elements can be characterized in several other
  ways, namely as finitely additive measures
  (Theorem~\ref{FinAddMeasThm}) or as maps of partially ordered vector
  spaces with strong unit (via Proposition~\ref{EModAdjProp}, see
  Remark~\ref{UltrafilterRem}~(3)).

\item Its (observable) algebras correspond to well-behaved
  mathematical structures (convex compact Hausdorff spaces), via the
  isomorphism $\Algobs(\Exp) \cong \CCHobs$ in
  Theorem~\ref{AlgobsCCHobsThm}.

\item There is a dual equivalence $\Algobs(\Exp) \simeq \BEMod\op$
that can be exploited for program logics, see~\cite{dHondtP06}.
\end{itemize}

It is thus worthwhile to systematically develop a program semantics
and logic based on the expectation monad and its duality. This is a
project on its own. We conclude by sketching some ingredients,
focusing on the program constructs that can be used. 

First we include a small example. Suppose we have a set of states
$S = \{a,b,c\}$ with probabilistic transitions between them as
described on the left below.
$$\vcenter
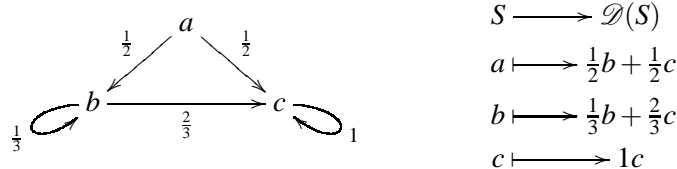
{\xymatrix@R-.5pc{
& a\ar[dl]_{\frac{1}{2}}\ar[dr]^{\frac{1}{2}} & \\
b\ar[rr]_-{\frac{2}{3}}\ar@(l,dl)_{\frac{1}{3}} & & c\ar@(r,dr)^{1}
}}
\qquad\qquad
\vcenter{\xymatrix@R-2pc{
S\ar[r] & \Dst(S) \\
a\ar@{|->}[r] & \frac{1}{2}b + \frac{1}{2}c \\
b\ar@{|->}[r] & \frac{1}{3}b + \frac{2}{3}c \\
c\ar@{|->}[r] & 1c
}}$$

\noindent On the right is the same system described as a function,
namely as coalgebra of the distribution monad $\Dst$. It maps each
state to the corresponding discrete probability distribution. We can
also describe the same system as coalgebra $S\rightarrow \Exp(S)$ of
the expectation monad, via the map $\Dst\rightarrowtail\Exp$. Then it
looks as follows:
$$\xymatrix@R-2pc@C+1pc{
S\ar[r] & \Exp(S) \\
a\ar@{|->}[r] & \lamin{q}{[0,1]^{S}}{\frac{1}{2}q(b) + \frac{1}{2}q(c)} \\
b\ar@{|->}[r] & \lamin{q}{[0,1]^{S}}{\frac{1}{3}q(b) + \frac{2}{3}q(c)} \\
c\ar@{|->}[r] & \lamin{q}{[0,1]^{S}}{q(c)}
}$$

\noindent Thus, via the $\Exp$-monad we obtain a probabilistic
continuation style semantics.

Let's consider this from a more general perspective. Assume we now
have an arbitrary, unspecified set of states $S$, for which we
consider programs as functions $S\rightarrow \Exp(S)$,
\textit{i.e.}~as Kleisli endomaps or as $\Exp$-coalgebras. In a
standard way the monad structure provides a monoid structure on these
maps $S\rightarrow\Exp(S)$ for sequential composition, with the unit
$S\rightarrow\Exp(S)$ as neutral element `skip'. We briefly sketch
some other algebraic structure on such programs (coalgebras), see
also~\cite{Panangaden09}.

Programs $S\rightarrow\Exp(S)$ are closed under convex combinations:
if we have programs $P_{1}, \ldots, P_{n}\colon S \rightarrow \Exp(S)$
and probabilities $r_{i}\in[0,1]$ with $\sum_{i}r_{i}=1$, then we can
form a new program $P = \sum_{i}r_{i}P_{i}\colon S \rightarrow
\Exp(S)$. For $q\in [0,1]^{S}$,
$$\begin{array}{rcl}
P(s)(q)
& = &
\sum_{i}r_{i}\cdot P_{i}(s)(q).
\end{array}$$

Since the sets $S\rightarrow\Exp(S)$ carries a pointwise order with
suprema of $\omega$-chains we can also give meaning to iteration
constructs like `while' and 'for \ldots do'.

Further we can also do ``probabilistic assignment'', written for
instance as $n := \varphi$, where $n$ is a variable, say of integer
type \textsf{int}, and $\varphi$ is a distribution of type
$\Dst(\textsf{int})$. The intended meaning of such an assignment $n :=
\varphi$ is that afterwards the variable $n$ has value
$m\colon\textsf{int}$ with probability $\varphi(m)\in[0,1]$.  In order
to model this we assume an update function $\textsf{upd}_{n}\colon
S\times \textsf{int} \rightarrow S$, which we leave unspecified
(similar functions exist for other variables). The interpretation
$\Scottint{n:=\varphi}$ of the probabilistic assignment is a function
$S\rightarrow\Exp(S)$, defined as follows.
$$\begin{array}{rcl}
\Scottint{n:=\varphi}(s)
& = &
\Exp\big(\textsf{upd}_{n}(s, -)\big)\big(\sigma(\varphi)\big) \\
& = &
\lamin{q}{[0,1]^{S}}{\sum_{i}r_{i}\cdot q(\textsf{upd}_{n}(m_{i}))},
   \qquad\mbox{if }\varphi=\sum_{i}r_{i}m_{i}.
\end{array}$$

\noindent It applies the functor $\Exp$ to the function
$\textsf{upd}_{n}(s,-)\colon \textsf{int} \rightarrow S$ and uses the
natural transformation $\sigma\colon\Dst\Rightarrow\Exp$
from~\eqref{ComposableAdjunctionNatroDiag}.

\bibliographystyle{eptcs}
\bibliography{../../Macros/bib}

\end{document}